\documentclass[10pt]{article}
\usepackage{latexsym,amsfonts,amssymb,amsmath,amsthm}
\usepackage{amscd}
\usepackage{graphicx}
\usepackage{cite}
\usepackage[usenames,dvipsnames]{color}
\usepackage{tocloft}
\usepackage[colorlinks=true,linkcolor=black,urlcolor=blue,citecolor=blue]{hyperref}
\hypersetup{
  linkcolor=blue,
  pdftitle={Robust dynamics of nonautonomous systems with nested invariant cone structures},
  pdfauthor={Dun Zhou},
  pdfsubject={Nested invariant cones, nonautonomous dynamics, and perturbations}
}

\begin{document}
	\setlength{\baselineskip}{13pt}

\parindent 0.5cm
\evensidemargin 0cm \oddsidemargin 0cm \topmargin 0cm \textheight
22cm \textwidth 16cm \footskip 2cm \headsep 0cm
	
\newtheorem{theorem}{Theorem}[section]
\newtheorem{lemma}[theorem]{Lemma}
\newtheorem{proposition}[theorem]{Proposition}
\newtheorem{definition}{Definition}[section]
\newtheorem{example}{Example}[section]
\newtheorem{corollary}[theorem]{Corollary}

\newtheorem{remark}{Remark}[section]
\newtheorem{property}[theorem]{Property}
\numberwithin{equation}{section}
\newtheorem{mainthm}{Theorem}
\newtheorem{mainlem}{Lemma}
	
\def\p{\partial}
\def\I{\textit}
\def\R{\mathbb R}
\def\C{\mathbb C}
\def\u{\underline}
\def\l{\lambda}
\def\a{\alpha}
\def\O{\Omega}
\def\e{\epsilon}
\def\ls{\lambda^*}
\def\D{\displaystyle}
\def\wyx{ \frac{w(y,t)}{w(x,t)}}
\def\imp{\Rightarrow}
\def\tE{\tilde E}
\def\tX{\tilde X}
\def\tH{\tilde H}
\def\tu{\tilde u}
\def\d{\mathcal D}
\def\aa{\mathcal A}
\def\DH{\mathcal D(\tH)}
\def\bE{\bar E}
\def\bH{\bar H}
\def\M{\mathcal M}
\renewcommand{\labelenumi}{(\arabic{enumi})}

\def\disp{\displaystyle}
\def\undertex#1{$\underline{\hbox{#1}}$}
\def\card{\mathop{\hbox{card}}}
\def\sgn{\mathop{\hbox{sgn}}}
\def\exp{\mathop{\hbox{exp}}}
\def\OFP{(\Omega,{\cal F},\PP)}
\newcommand\JM{Mierczy\'nski}
\newcommand\RR{\ensuremath{\mathbb{R}}}
\newcommand\CC{\ensuremath{\mathbb{C}}}
\newcommand\QQ{\ensuremath{\mathbb{Q}}}
\newcommand\ZZ{\ensuremath{\mathbb{Z}}}
\newcommand\NN{\ensuremath{\mathbb{N}}}
\newcommand\PP{\ensuremath{\mathbb{P}}}
\newcommand\abs[1]{\ensuremath{\lvert#1\rvert}}

\newcommand\normf[1]{\ensuremath{\lVert#1\rVert_{f}}}
\newcommand\normfRb[1]{\ensuremath{\lVert#1\rVert_{f,R_b}}}
\newcommand\normfRbone[1]{\ensuremath{\lVert#1\rVert_{f, R_{b_1}}}}
\newcommand\normfRbtwo[1]{\ensuremath{\lVert#1\rVert_{f,R_{b_2}}}}
\newcommand\normtwo[1]{\ensuremath{\lVert#1\rVert_{2}}}
\newcommand\norminfty[1]{\ensuremath{\lVert#1\rVert_{\infty}}}

\title{Robust dynamics of non-autonomous systems with nested invariant cone structures}

\author {Dun Zhou\thanks{Partially supported by the National Natural Science Foundation of China
under Grant Nos. 12671213 and 12331006.}\\
 School of Mathematics and Statistics
 \\ Nanjing University of Science and Technology
 \\ Nanjing, Jiangsu, 210094, P. R. China
 \\
\\
}
\date{}

\maketitle

\begin{abstract}
We develop a perturbation theory for skew-product semiflows with nested
invariant cone (NIC) structures over compact almost periodic minimal flows.
The theory describes how the restrictions imposed by these structures on
asymptotic dynamics persist under small $C^1$ perturbations.
We first establish a trichotomy for differences of perturbed trajectories:
finite-time coincidence, exponential contraction outside the terminal cone,
or eventual separation within a fixed cone layer. On compact invariant
sets near the reference set, the perturbed semiflow extends to a flow.
Under a transversality condition, compact invariant sets whose orbit
differences remain in fixed cone layers admit finite-dimensional fibre
realizations over the original base flow.
In the codimension-one case, every compact minimal set in the perturbative
neighborhood is an almost $1$-cover of the base. For precompact positive
orbits whose $\omega$-limit sets lie in this neighborhood, we prove that
each $\omega$-limit set contains at most two minimal sets and classify the
remaining orbits by their forward and backward limit sets.
Applications include nonlocal and chemotaxis-type perturbations of scalar
parabolic equations with homogeneous Neumann boundary conditions and
perturbations of competitive--cooperative tridiagonal systems. Under
dissipativity assumptions, these applications give global conclusions for
almost periodically forced equations, including perturbations that destroy
the original zero-number or sign-variation mechanism.

\end{abstract}



\tableofcontents

\section{Introduction}

A central question in the qualitative theory of differential equations is
under what conditions a high-dimensional or even infinite-dimensional
dynamical system admits a finite-dimensional description of its asymptotic
dynamics. Invariant cone
structures provide an effective mechanism for such a reduction: they
constrain the evolution of differences between solutions and often force the
recurrent dynamics to exhibit an ordering or oscillatory structure of much
lower dimension than the original phase space. Besides their fundamental role
in monotone dynamical systems, invariant cones are closely related to
hyperbolic structures and finite-dimensional reduction; see, for example,
\cite{KH,LSY}.

An important form of this mechanism is provided by a nested invariant cone
(NIC) structure. Roughly speaking, it consists of a hierarchy
\[
C_0=\{0\}\subset C_1\subset C_2\subset\cdots,
\]
indexed up to $N\in\{1,2,\ldots\}\cup\{\infty\}$, with each $C_i$ of
finite rank and $C_N=X$ when $N<\infty$. The successive cone layers distinguish different
oscillatory modes of orbit differences and, together with exponential
separation and a suitable transversality property, impose strong restrictions
on the asymptotic dynamics. In particular, an NIC structure can reduce the
effective dimension of recurrent dynamics even when the original phase space
is infinite-dimensional.

NIC structures arise naturally in differential equations governed by
sign-variation or ordering mechanisms. A classical precursor is the
sign-change counting method of Gantmacher for oscillatory matrices
\cite[p.~105]{Ga}. For scalar parabolic equations, Matano's zero-number
function \cite{H.MATANO:1982} provides the fundamental mechanism, while
analogous sign-variation properties occur in tridiagonal
competitive--cooperative systems. Further examples include monotone cyclic
feedback systems, systems with time delay, and strongly $k$-cooperative
systems; see
\cite{BSZM,Fiedler,Fusco1987,Fusco1990,Gedeon,Smillie1984,
Ma-Sm,Ma-Se1,Ma-Se2,Ma-Se3,Matano,MS,weiss2021_sub}.

For autonomous and time-periodic systems, these structures have led to
powerful convergence, embedding, and Poincar\'e--Bendixson type results; see
\cite{H.MATANO:1982,Smillie1984,CCH,JR2,Te}. As a classical example,
consider the scalar reaction--diffusion equation with homogeneous Neumann
boundary conditions
\begin{equation}\label{parabolic-eq-01}
\begin{cases}
u_t=u_{xx}+f(t,x,u,u_x), & t>0,\ 0<x<1,\\
u_x(t,0)=u_x(t,1)=0.
\end{cases}
\end{equation}
For autonomous equations, the zero-number method is a basic ingredient in
the convergence of bounded solutions to equilibria
\cite{H.MATANO:1982}. In the time-periodic case, the corresponding
asymptotic objects are periodic solutions; see \cite{BPS,CM}. Within the
abstract NIC framework, Tere\v{s}\v{c}\'ak \cite{Te} obtained
finite-dimensional embeddings of $\omega$-limit sets for autonomous and
time-periodic systems and established the robustness of this structure under
small $C^1$-perturbations. In particular, the codimension-one case yields
convergence-type dynamics, whereas the codimension-two case leads to a
Poincar\'e--Bendixson type description.

Almost periodic forcing generally admits no common return time, so the
recurrent dynamics must be studied together with the motion on the base
flow. The natural recurrent objects are then minimal extensions of the
base, rather than necessarily equilibria or periodic orbits. Almost
$1$-covering properties and restrictions on the number of minimal sets
are already known for several classes of nonautonomous equations; see
\cite{wangyi2007,Fangchun2013,Shen1998,SWZ4}. The question addressed here
is whether such conclusions follow from an abstract NIC structure and
remain valid when a small perturbation destroys the particular
oscillation mechanism of the original equation.

A representative problem is furnished by the nonlocal perturbation
\begin{equation}\label{parabolic-eq-induced0}
\begin{cases}
u_t=u_{xx}+f(t,x,u,u_x)
+\varepsilon c(t,x)\displaystyle\int_0^1\nu(y)u(t,y)\,dy,
& t>0,\ 0<x<1,\\
u_x(t,0)=u_x(t,1)=0,
\end{cases}
\end{equation}
where $f$ and $c$ are uniformly almost periodic in time. Perturbations of
this type are particularly relevant here because the nonlocal term generally
destroys the local scalar equation satisfied by the difference of two
solutions and hence the exact Sturm zero-number mechanism. The same issue
arises for chemotaxis-type couplings and for general perturbations of
competitive--cooperative tridiagonal systems. It is therefore natural to ask
whether the low-dimensional recurrent dynamics induced by the unperturbed NIC
structure persists when the perturbation is sufficiently small. The smallness
condition is essential: even in the autonomous case, nonlocal perturbations
of the form \eqref{parabolic-eq-induced0} may produce complicated dynamics
when $|\varepsilon|$ is not small; see \cite{Fie-Pol}.

Related recent work by Cheng, Wang, and Zhou \cite{ChengWangZhou2025}
uses nested invariant cones to establish persistence of the
Poincar\'e--Bendixson property under $C^1$ perturbations and the generic
Morse--Smale property for dissipative bidirectional cyclic negative
feedback systems.
Here the focus is on minimal extensions and $\omega$-limit sets over
almost periodic base flows, including perturbations that do not retain
the original local or tridiagonal structure.

We work with skew-product semiflows over compact almost periodic minimal
base flows. Following the perturbative cone method of
Tere\v{s}\v{c}\'ak \cite[Propositions~5.1 and~6.1]{Te}, we combine the
unperturbed difference cocycle with uniform finite-time estimates.
Compactness and regularization reduce the relevant cone hierarchy to
finitely many levels near a compact invariant set. The main additional
task is to turn these estimates into information about minimal extensions
and $\omega$-limit sets over an aperiodic base.

\vskip 3mm
\textbf{The main results are as follows.}

First, we establish an abstract perturbation theory for NIC systems
using nonautonomous skew-product semiflows. Theorem~\ref{prop-orbit}
shows that the difference of two perturbed trajectories that remain near the
reference compact invariant set satisfies one of three alternatives: the
trajectories coincide after a finite time, their difference remains outside
the largest perturbed cone and decays exponentially, or their difference is
eventually confined to a single cone layer. On a compact invariant set, the restriction of the perturbed semiflow
extends to a flow, and distinct points in a common fibre are eventually
separated by a fixed cone layer. Thus the essential cone-dynamical structure
is robust under small $C^1$-perturbations in a genuinely nonautonomous
setting.

Second, under a transversality condition involving a codimension-$d$
subspace, we obtain a finite-dimensional fibre realization on suitable
compact invariant sets. More precisely, Corollary~\ref{imbeding-d-plane}
shows that, if each nonzero orbit difference remains in a fixed cone
layer for all time, the restricted flow is topologically conjugate to a
flow on a compact invariant subset of
$
\mathbb R^d\times\Theta.
$
The reduction is finite-dimensional in the state variable while retaining
the original base flow, which need not itself be finite-dimensional.

Third, we obtain a detailed description of the recurrent dynamics in the
codimension-one case. Every compact minimal set contained in the perturbative
neighborhood is an almost $1$-cover, and hence an almost automorphic extension,
of the base flow; see Theorem~\ref{a-a-one}. This is the natural
nonautonomous counterpart of an equilibrium or a periodic orbit. Moreover,
Theorem~\ref{perturbation-thm} shows that the $\omega$-limit set of a
precompact positive orbit contains at most two minimal sets whenever it
lies in the perturbative neighborhood
and admits a three-way classification in terms of the minimal sets
contained in the forward and backward limit sets of the remaining orbits.
The order on each fibre thus controls the number and arrangement of
minimal sets even when the phase space is infinite-dimensional.

We apply the abstract theory to almost-periodically forced scalar parabolic
equations with homogeneous Neumann boundary conditions, including nonlocal
and chemotaxis-type perturbations, and to competitive--cooperative
tridiagonal systems. These applications show that the asymptotic ordering and
low-dimensional recurrent dynamics generated by an NIC structure may persist
even when a perturbation destroys the exact zero-number or sign-variation
property of the unperturbed equation. Under suitable dissipativity and
compactness assumptions, upper semicontinuity of global attractors allows the
local abstract theory to yield global conclusions for the perturbed
equations.

The rest of the paper is organized as follows. Section~2 introduces the
abstract framework and states the main results. Section~3 proves the main
structural theorems. Section~4 applies the theory to scalar parabolic
equations, their nonlocal and chemotaxis-type perturbations, and
competitive--cooperative tridiagonal systems. Section~5 contains further
remarks on the abstract results and their applications. The technical proof
of Proposition~\ref{pertur-pro}, which provides the perturbative cone
estimates used throughout the paper, is deferred to Section~6.

\section{Preliminaries and Main Results}

Throughout this section we assume that $(X,\|\cdot\|)$ is a real Banach space, $Y$ is a metric space, $\Theta$ and $Z$ are compact metric spaces, and $\mathcal{E}=X\times \Theta$. We write
$P:X\times\Theta\to\Theta$ for the natural projection.

\subsection{Basic definitions}

In this subsection we introduce the notions of cones, skew-product semiflows, exponential separation, almost periodic and almost automorphic functions and flows, and recall some basic properties.

A closed subset $C$ of $X$ is called a {\it $k$-cone} if $\lambda C=C$ for all $\lambda\in\RR\setminus\{0\}$ and there exist a $k$-dimensional subspace $V_0$ of $X$ and a subspace $X_0$ of $X$ with codimension $k$ such that $V_0\subset C$ and $X_0\cap C=\{0\}$. A cone is called {\it solid} if it has nonempty interior.

A continuous map
\[
\sigma:Y\times\RR\to Y,\qquad (y,t)\mapsto y\cdot t,
\]
denoted by $(Y,\sigma)$ or $(Y,\RR)$, is called a {\it continuous flow} on $Y$ if $\sigma(y,0)=y$ and $\sigma(\sigma(y,s),t)=\sigma(y,s+t)$ for all $y\in Y$ and $s,t\in\RR$. A subset $\tilde Y\subset Y$ is {\it invariant} if $\sigma_t(\tilde Y)=\tilde Y$ for every $t\in\RR$. A subset $\tilde Y\subset Y$ is called {\it minimal} if it is nonempty, compact, invariant, and the only nonempty compact invariant subset of $\tilde Y$ is $\tilde Y$ itself. Every nonempty compact $\sigma$-invariant set contains a minimal subset, and a compact invariant set $\tilde Y$ is minimal if and only if every trajectory in $\tilde Y$ is dense. The continuous flow $(Y,\sigma)$ is called {\it minimal} if $Y$ is minimal.

The following lemma follows from \cite[Lemma~I.2.16]{Shen1998}.

\begin{lemma}\label{epimorphism-thm}
Let $p:(Z,\RR)\rightarrow(\Theta,\RR)$ be an epimorphism of compact metric flows. Then the set
\[
\Theta'=\Big\{\theta_{0}\in \Theta:\ \text{for any }z_{0}\in p^{-1}(\theta_{0}),\ \theta\in \Theta,\ \text{and any sequence } \{t_{i}\}\subset \RR \text{ with }\theta\cdot t_{i}\rightarrow \theta_{0},
\]
\[
\text{there exists a sequence } \{z_{i}\}\subset p^{-1}(\theta)\ \text{such that}\ z_{i}\cdot t_{i}\rightarrow z_{0}\Big\}
\]
is residual and invariant. In particular, if $(Z,\RR)$ is minimal and distal, then $\Theta'=\Theta$.
\end{lemma}

\begin{definition}
{\rm
Let $p:(Z,\RR)\rightarrow(\Theta,\RR)$ be a homomorphism of minimal flows.
The flow $(Z,\RR)$ is said to be an {\it almost $1$-cover} or an
{\it almost automorphic extension of $(\Theta,\RR)$} if there exists
$\theta_0\in \Theta$ such that ${\rm card}\,p^{-1}(\theta_0)=1$.
}
\end{definition}

\begin{definition}\label{skew-product}
{\rm
Suppose that $\sigma(\theta,t)=\theta\cdot t$ is a compact flow on
$\Theta$. A {\it skew-product semiflow}
$\Lambda_t:\mathcal{E}\rightarrow \mathcal{E}$ is a semiflow of the form
\begin{equation}\label{skew-product-semiflow}
\Lambda_{t}(x,\theta)
=(\phi(t,x,\theta),\theta\cdot t),
\quad t\geq 0,\ (x,\theta)\in \mathcal{E}.
\end{equation}
}
\end{definition}

We always write the skew-product semiflow defined above as
$\Lambda=(\phi,\sigma)$.
A set $E\subset\mathcal E$ is called \emph{positively invariant} if
$\Lambda_t(E)\subset E$ for every $t\ge0$, and \emph{invariant} if
$\Lambda_t(E)=E$ for every $t\ge0$. Thus compact invariant sets of a
semiflow admit complete backward extensions, although these extensions need
not be unique a priori.

Let $L(X,X)$ be the set of all bounded linear operators on $X$. Suppose
$(Z,\sigma)$ is a compact flow and let
$\{T(t,z):t>0,\ z\in Z\}$ be a family of bounded linear operators on $X$
with the following properties:

\medskip

\noindent$\bullet$
$\displaystyle\lim_{t\to 0+}T(t,z)x=x$ for all $z\in Z$ and $x\in X$;

\smallskip

\noindent$\bullet$
the map
$(t,z)\mapsto T(t,z):(0,+\infty)\times Z\to L(X,X)$
is continuous in the operator norm;

\smallskip

\noindent$\bullet$
$T(t_1+t_2,z)=T(t_2,z\cdot t_1)\circ T(t_1,z)$ for all
$t_1,t_2>0$ and $z\in Z$.

\medskip

Together with the base flow, this family defines a linear vector bundle
semiflow for positive times in the sense of \cite{Te}, denoted by
$\bar\Lambda=(T,\sigma)$. We set $T(0,z)=I$ for convenience.
The exponential-separation results below are formulated for these
positive-time operator families.

Let $K=X\times Z$. A subset
$\tilde K=\bigcup_{z\in Z}\tilde K(z)
=\bigcup_{z\in Z}(X_z\times\{z\})$
of $K$ is called a {\it $k$-dimensional vector subbundle of $K$ over $Z$}
if, for each $z\in Z$, $X_z\subset X$ is a $k$-dimensional subspace of $X$.

\begin{definition}\label{ES-define}
{\rm
A linear vector bundle semiflow $\bar\Lambda=(T,\sigma)$ on $X\times Z$
is said to admit a
{\it $k$-dimensional continuous separation along $Z$ associated with a
$k$-cone $C$}
if there exist $k$-dimensional subbundles
$\bigcup_{z\in Z}V_z\times\{z\}$, with $V_z\subset C$, and
$\bigcup_{z\in Z}L_z\times\{z\}$ of $X\times Z$ and $X^*\times Z$,
where $X^*$ is the dual space of $X$, respectively, such that:
\begin{itemize}
\item[\rm(1)]
({\it Continuity})
$\bigcup_{z\in Z}V_z\times\{z\}$ and
$\bigcup_{z\in Z}L_z\times\{z\}$ are continuous bundles: the maps
$z\mapsto V_z$ and $z\mapsto L_z$ are continuous in the gap topology
on the respective Grassmann spaces of $k$-dimensional subspaces.
In particular, the projections onto $V_z$ along $\operatorname{Ann}(L_z)$
depend continuously on $z$ in the operator norm.

\item[\rm(2)]
({\it Invariance})
$T(t,z)V_z=V_{z\cdot t}$ and
$T^*(t,z)L_{z\cdot t}=L_z$
for all $t>0$ and $z\in Z$, where $T^*(t,z)$ denotes the adjoint
operator of $T(t,z)$.

\item[\rm(3)]
({\it Exponential separation})
There exist constants $M,\gamma>0$ such that
\[
\|T(t,z)w\|
\leq
M e^{-\gamma t}\,\|T(t,z)v\|
\]
for all $z\in Z$, $w\in\operatorname{Ann}(L_z)\cap S$,
$v\in V_z\cap S$, and $t\geq1$, where
\[
\operatorname{Ann}(L_z)
=
\{v\in X:\ \ell(v)=0\ \text{for all }\ell\in L_z\}
\]
with $\operatorname{Ann}(L_z)\cap C=\{0\}$, and
$S=\{u\in X:\ \|u\|=1\}$.
\end{itemize}
}
\end{definition}

The following exponential separation lemma is taken from
\cite[Corollary~2.2]{Te}.

\begin{lemma}\label{ES-lm}
Let $(X,\|\cdot\|)$ be a Banach space, let $C$ be a $k$-cone in $X$, and
let $\bar\Lambda=(T,\sigma)$ be a linear vector bundle semiflow on
$X\times Z$. Assume that $T(t,z)$ is compact for each
$(t,z)\in(0,+\infty)\times Z$, and that for any given
$v\in C\setminus\{0\}$ and $(t,z)\in(0,+\infty)\times Z$ there exists a
neighborhood $\mathcal{N}(t,z,v)$ of $v$ in $X$ such that
\[
T(t,z)\,\mathcal{N}(t,z,v)\subset C\setminus\{0\}.
\]
Then the vector bundle semiflow $\bar\Lambda$ admits a $k$-dimensional
continuous separation along $Z$ associated with the $k$-cone $C$.
\end{lemma}

In the following definitions, let $D=\RR^m$. We recall some
definitions and properties related to almost periodic and almost automorphic
functions (and flows).

\begin{definition}\label{admissible}
{\rm
A function $f\in C(\RR\times\RR^n,D)$ is said to be {\it admissible} if,
for any compact subset $K\subset\RR^n$, the function $f$ is bounded and
uniformly continuous on $\RR\times K$. We say that $f$ is
{\it $C^r$ admissible} ($r\geq1$) if $f(t,w)$ is $C^r$ in
$w\in\RR^n$, Lipschitz continuous in $t$ uniformly for $w$ in compact
sets, and $f$ and all its partial derivatives in $w$ up to order $r$
are admissible.
}
\end{definition}

Let $f\in C(\RR\times\RR^n,D)$ be an admissible function. Then
\[
H(f)={\rm cl}\{f\cdot\tau:\ \tau\in\RR\},
\]
called the {\it hull} of $f$, is compact and metrizable under the
compact-open topology (see \cite{Sell,Shen1998}), where
$(f\cdot\tau)(t,\cdot)=f(t+\tau,\cdot)$. Moreover, the time translation
$g\mapsto g\cdot t$ for $g\in H(f)$ induces a natural flow on $H(f)$
(cf.~\cite{Sell}).

\begin{definition}\label{almost}
{\rm
\begin{itemize}
\item[(1)]
A function $f\in C(\RR,D)$ is {\it almost periodic} if, for every sequence
$\{t'_k\}\subset\RR$, there exists a subsequence
$\{t_k\}\subset\{t'_k\}$ such that $\{f(t+t_k)\}$ converges uniformly in
$t\in\RR$.

\item[(2)]
A function $f\in C(\RR,D)$ is {\it almost automorphic} if, for every
sequence $\{t'_k\}\subset\RR$, there exist a subsequence
$\{t_k\}\subset\{t'_k\}$ and a function $g:\RR\to D$ such that
$f(t+t_k)\to g(t)$ and $g(t-t_k)\to f(t)$ for each $t\in\RR$.

\item[(3)]
A function $f\in C(\RR\times\RR^n,D)$ is
{\it uniformly almost periodic in $t$}
(resp.\ {\it uniformly almost automorphic in $t$}) if $f$ is admissible
and, for each fixed $d\in\RR^n$, the function $t\mapsto f(t,d)$ is almost
periodic (resp.\ almost automorphic) in $t\in\RR$.
\end{itemize}
}
\end{definition}

Let $(X,\RR)$ be a compact flow. For $x\in X$ and a net
$\alpha=\{t_n\}$ in $\RR$, define
\[
T_{\alpha}x=\lim_{n\to\infty}x\cdot t_n,
\]
if the limit exists.

\begin{definition}\label{ap-minimal-flow}
{\rm
A compact flow $(X,\RR)$ is called {\it almost periodic} if any nets
$\alpha'$, $\beta'$ in $\RR$ admit subnets
$\alpha=\{t_n\}$, $\beta=\{s_n\}$ such that
\[
T_{\alpha}T_{\beta}x=T_{\alpha+\beta}x
\]
for all $x\in X$, where $\alpha+\beta=\{t_n+s_n\}$.
}
\end{definition}

\begin{definition}
{\rm
Let $(X,\RR)$ be a compact flow. A point $x\in X$ is called an
{\it almost automorphic point} if any net
$\alpha'=\{t'_n\}$ in $\RR$ has a subnet $\alpha=\{t_n\}$ such that both
$T_{\alpha}x$ and $T_{-\alpha}T_{\alpha}x$ exist and
$T_{-\alpha}T_{\alpha}x=x$. The flow $(X,\RR)$ is called
{\it almost automorphic} if there exists an almost automorphic point
$x_0\in X$ with dense orbit.
}
\end{definition}

\begin{remark}\label{a-p-to-minial}
{\rm Assume that $f$ is uniformly almost periodic (resp.\ uniformly almost automorphic) in $t$. Then the following hold.
 \begin{itemize}
   \item[(1)] $H(f)$ is minimal, and $(H(f),\RR)$ is an almost periodic (resp.\ almost automorphic) minimal flow.
   \item[(2)] Every (resp.\ residually many) $g\in H(f)$ is a uniformly almost periodic (resp.\ uniformly almost automorphic) function.
   \item[(3)] There exists a unique $F\in C(H(f)\times \RR^n, D)$ such that $F(g\cdot t,z)\equiv g(t,z)$; moreover, if $f$ is $C^r$ admissible, then $F$ is $C^r$ in $z\in \RR^n$, and the derivatives with respect to $z$ depend continuously on $g\in H(f)$ (see, e.g., \cite{Shen1998}).
 \end{itemize}}
\end{remark}

Before ending this subsection, we introduce the notion of $C^1$-functions on a subset of $X$. Let $B\subset X$ be a bounded open set and define
\[
\begin{aligned}
C^1(B,X):=\{F:B\to X:\;&F \text{ is continuously Fr\'{e}chet
differentiable on }B,\\
&F\text{ and }DF\text{ are bounded on }B\}.
\end{aligned}
\]
with the norm
\[
\|F\|_{C^1(B,X)}=\sup_{x\in B}\|F(x)\|_X+\sup_{x\in B}\|DF(x)\|_{L(X,X)}.
\]

\subsection{Assumptions}

In this subsection, we formulate the standing assumptions on the
skew-product semiflow and its associated linear difference cocycle. These
assumptions provide the abstract NIC framework used throughout the paper.

With $X$ and $L(X,X)$ as above, assume that there exists a nested family of solid
$k_i$-cones $C_i$, indexed by $1\le i\le N$ when $N<+\infty$ and by
$i\ge1$ when $N=+\infty$, where
\[
N\in\NN\cup\{+\infty\}.
\]
Set
\[
C_0=\{0\},\qquad k_0=0,
\]
and assume that
\[
C_i\supset C_{i-1},
\qquad
k_i>k_{i-1},
\]
for every $1\le i\le N$ if $N<+\infty$, and for every $i\ge1$ if
$N=+\infty$. If $N<+\infty$, we further assume that
\[
C_N=X.
\]
No condition of the form
\[
\bigcup_{i\ge1}C_i=X
\]
is imposed when $N=+\infty$. Here and below, an admissible index means
$1\le i\le N$ when $N<+\infty$, and $i\ge1$ when $N=+\infty$.

When $N<+\infty$, the condition $C_N=X$ implies that $X$ is
finite-dimensional, since $C_N$ is assumed to have finite rank. The case
$N=+\infty$ covers the infinite-dimensional applications. The absence of
an exhaustion assumption on $X$ is important for the Sturm cones considered
below, since a general $C^1$ function need not belong to any finite
zero-number cone. Instead, the regularizing part of
{\rm\textbf{(H4)}} requires every nonzero image over the fixed time
interval used below to enter some finite cone level.

Let $(\Theta,\sigma)$ be an almost periodic minimal flow, and let
$\mathcal U_0\subset X\times\Theta$ be an open set that projects onto
$\Theta$. Let
\[
\Lambda_t(x,\theta)
=
\bigl(\phi(t,x,\theta),\theta\cdot t\bigr),
\qquad t\ge0,
\]
be a skew-product semiflow on $\mathcal U_0$ that is $C^1$ in the
$x$-variable.

For each $\theta\in\Theta$, define
\[
\mathcal U_0(\theta)
=
\{x\in X:(x,\theta)\in\mathcal U_0\},
\]
and set
\[
\mathcal U_0^e
=
\left\{
(x,y,\theta)\in X\times X\times\Theta:
x,y\in\mathcal U_0(\theta)
\right\}.
\]
For $t\ge0$ and $(x,y,\theta)\in\mathcal U_0^e$, suppose that there
exists a bounded linear operator
\[
T(t,x,y,\theta)\in L(X,X)
\]
such that
\begin{equation}\label{operator-T-eq1}
T(t,x,y,\theta)(x-y)
=
\phi(t,x,\theta)-\phi(t,y,\theta)
\end{equation}
for every $t\ge0$ and every
$(x,y,\theta)\in\mathcal U_0^e$.

We call $T$ an \emph{associated linear difference cocycle}, with the
cocycle property specified in {\bf(H1)} below. When $[x,y]$ lies in the
fibre domain, the averaged derivative
\[
A(t,x,y,\theta)
:=\int_0^1 D_x\phi\bigl(t,y+s(x-y),\theta\bigr)\,ds
\]
satisfies the difference identity \eqref{operator-T-eq1}, but need not
satisfy the cocycle identity: averaging derivatives does not in general
commute with composing time maps. The usual variational cocycle along a
single orbit, in turn, need not represent finite differences.

In the differential-equation applications, $T$ is obtained from the
linear equation satisfied by the difference of two solutions. Its
coefficients are formed by averaging the derivative of the vector field
along the segment joining the two solutions at each time. The evolution
operator satisfies \eqref{operator-T-eq1}, and uniqueness and time
translation of this linear equation give the cocycle property.

We now state the standing assumptions.

\begin{itemize}

\item[\rm\textbf{(H1)}]
({\rm Cocycle property and continuity})
For every $(x,y,\theta)\in\mathcal U_0^e$, let
\[
T(0,x,y,\theta)=I.
\]
For all $t_1,t_2\ge0$ and all
$(x,y,\theta)\in\mathcal U_0^e$,
\[
\begin{aligned}
T(t_1+t_2,x,y,\theta)
={}&
T\bigl(
t_2,
\phi(t_1,x,\theta),
\phi(t_1,y,\theta),
\theta\cdot t_1
\bigr)\\
&\circ T(t_1,x,y,\theta).
\end{aligned}
\]
Moreover,
\[
\lim_{t\to0^+}T(t,x,y,\theta)v=v
\qquad\text{for every }v\in X,
\]
and the map
\[
(t,x,y,\theta)\longmapsto T(t,x,y,\theta):
(0,1]\times\mathcal U_0^e\longrightarrow L(X,X)
\]
is continuous in the operator norm. Finally, for every compact set
$K\subset\mathcal U_0^e$,
\[
\sup_{\substack{0<t\le1\\ (x,y,\theta)\in K}}
\|T(t,x,y,\theta)\|_{L(X,X)}
<+\infty.
\]

\item[\rm\textbf{(H2)}]
({\rm Compactness})
For every $t>0$ and every
$(x,y,\theta)\in\mathcal U_0^e$, the operator
$T(t,x,y,\theta)$ is compact.

\item[\rm\textbf{(H3)}]
({\rm Injectivity})
For every $t>0$, every
$(x,y,\theta)\in\mathcal U_0^e$, and every
$v\in X\setminus\{0\}$,
\[
T(t,x,y,\theta)v\neq0.
\]

\item[\rm\textbf{(H4)}]
({\rm Strong cone preservation and regularization})
For every $t>0$, every $(x,y,\theta)\in\mathcal U_0^e$, every
admissible index $i$, and every $v\in C_i\setminus\{0\}$, there exists
a neighborhood
\[
\mathcal N(t,x,y,\theta,v)
\]
of $v$ in $X$ such that
\[
T(t,x,y,\theta)
\bigl(\mathcal N(t,x,y,\theta,v)\bigr)
\subset C_i\setminus\{0\}.
\]
If $N=+\infty$, we further assume that
\begin{equation}\label{regularizing-image}
T(t,x,y,\theta)(X\setminus\{0\})
\subset
\bigcup_{j\ge1}(C_j\setminus\{0\})
\qquad\text{for every }\frac12\le t\le1.
\end{equation}

\item[\rm\textbf{(H5)}]
({\rm Transversality with a codimension-$d$ subspace})
There exists a closed linear subspace $\mathcal H\subset X$ of
codimension $d$ such that, for every $t>0$, every
$(x,y,\theta)\in\mathcal U_0^e$, every $v\in X$, and every admissible
index $i$,
\[
T(t,x,y,\theta)v
\in
(C_i\setminus C_{i-1})\cap\mathcal H
\quad\Longrightarrow\quad
v\notin C_i.
\]

\end{itemize}

\begin{remark}\label{new-rk1}
{\rm
\begin{itemize}

\item[(1)]
Assumptions {\bf(H1)}--{\bf(H3)} provide the basic analytic framework:
the difference-cocycle structure, compactness, and injectivity.
Assumptions {\bf(H4)} and {\bf(H5)} encode the geometric properties of
the nested cone family.

By \eqref{operator-T-eq1} and {\bf(H3)}, two distinct points in the same
fibre have distinct positive-time images. Assumption {\bf(H4)} expresses
strong local forward preservation of each cone $C_i$. In the
infinite-dimensional case, \eqref{regularizing-image} is the
continuous-time skew-product analogue of the image-exhaustion condition
in Tere\v{s}\v{c}\'ak's assumption {\rm(H)}; compare also condition
{\rm(iii)} in \cite[Corollary~2.3]{Te}. It does not require every vector
in $X$ to belong to a finite cone level, but only every nonzero image over
the positive-time interval relevant to the perturbation argument.

The role of {\bf(H5)} is to prevent an orbit difference from remaining
in the ``bad'' subspace $\mathcal H$. More precisely, suppose that
\[
\phi(t,x,\theta)-\phi(t,y,\theta)
\in
(C_i\setminus C_{i-1})\cap\mathcal H
\]
for some $t>0$. By \eqref{operator-T-eq1} and {\bf(H5)}, the earlier
difference $x-y$ cannot belong to $C_i$. Thus, whenever the earlier
difference already belongs to some finite cone level, its minimal cone
level is strictly larger than that of its image. In this sense, an
intersection with $\mathcal H$ produces a strict decrease of the cone
level.

By \eqref{regularizing-image}, every nonzero orbit difference enters some
finite cone level after a positive time. Since the cone indices are
positive integers, the strict decrease described above can occur only
finitely many times in the discrete time comparisons used later.
For a perturbed pair whose difference enters
$\mathcal C_{N_0}\setminus\{0\}$, the corresponding conclusion follows
from Theorem~\ref{prop-orbit}(1) and (3): after a finite time, its
difference lies in
\[
\mathcal C_{i_0}\setminus
\bigl(\mathcal H\cup\mathcal C_{i_0-1}\bigr)
\]
for some index $i_0$. This exclusion of $\mathcal H$ is the mechanism
underlying the finite-dimensional reduction established later.

\item[(2)]
Suppose that $(\Theta,\sigma)$ is the hull flow associated with a
$T$-periodic forcing. Then
\[
\theta\cdot T=\theta
\qquad\text{for every }\theta\in\Theta,
\]
so the time-$T$ map of the skew-product semiflow leaves each fibre
invariant and defines the corresponding Poincar\'e map. The resulting
discrete-time system falls within the NIC framework of
Tere\v{s}\v{c}\'ak \cite{Te}.

Under the hypotheses of \cite[Theorem~C]{Te}, the $\omega$-limit set of
every orbit of a sufficiently small $C^1$ perturbation whose closure
remains in the prescribed neighborhood is homeomorphic to a subset of
$\RR^d$. In particular, when $d=1$ and the additional orientation
condition {\rm(A)} in \cite[Corollary~2.1]{Te} is satisfied, every such
$\omega$-limit set consists of a single fixed point of the Poincar\'e map
and hence corresponds to a periodic orbit of the continuous-time system.

\end{itemize}
}
\end{remark}

\subsection{Main results}

We now state the main abstract results. We first establish a trichotomy for
differences of perturbed trajectories near a compact invariant set of the
unperturbed skew-product semiflow. We then show that, on compact invariant
sets, the perturbed semiflow extends to a flow and admits an eventual
cone-layer separation. Under an additional transversality condition,
suitable compact invariant sets admit a finite-dimensional realization
in the fibre variables, provided that their nonzero orbit differences
remain in fixed cone layers for all time; see
Corollary~\ref{imbeding-d-plane}. Finally, we specialize to the
codimension-one case and obtain an almost $1$-covering property for
minimal sets and a trichotomy for $\omega$-limit sets.

The perturbation is measured directly in the fibrewise $C^1$ topology of
the time maps. This is the skew-product counterpart of the perturbation
condition in \cite[Theorem~C and Corollary~2.3]{Te}. In the proof,
auxiliary one-step difference operators are constructed and compared with
the unperturbed difference cocycle, following the idea of
\cite[(5.31)]{Te}.

For later use, associate with a compact invariant set
$K_0\subset\mathcal U_0$ its doubled invariant set
\[
Z(K_0):=\{(x,y,\theta):(x,\theta),(y,\theta)\in K_0\}.
\]

\begin{theorem}\label{prop-orbit}
Suppose that {\rm\textbf{(H1)}--\textbf{(H4)}} hold, and let
$K_0\subset\mathcal U_0$ be a compact invariant set of $\Lambda_t$.
Then there exist $N_0\in\NN$, $\epsilon_0>0$, open bounded neighborhoods
$\mathcal V_0$ and $\mathcal W_0$ of $K_0$ satisfying
\[
\overline{\mathcal V_0}\subset\mathcal W_0\subset\mathcal U_0,
\]
and a finite nested family of $k_i$-cones
$\{\mathcal C_i\}_{i=0}^{N_0}$ such that
\[
\mathcal C_0=\{0\},\qquad
\mathcal C_i\setminus\{0\}\subset\operatorname{int}\mathcal C_{i+1},
\qquad 1\le i<N_0.
\]
For $\theta\in\Theta$, set
\[
\mathcal V_0(\theta)
:=\{x\in X:(x,\theta)\in\mathcal V_0\}.
\]

Let
\[
\widetilde\Lambda_t
=\bigl(\widetilde\phi(t,\cdot,\theta),\theta\cdot t\bigr),
\qquad t\ge0,
\]
be a skew-product semiflow on $\mathcal U_0$ such that
$\widetilde\phi(t,\cdot,\theta)\in
C^1(\mathcal W_0(\theta),X)$ for $t\in[1/2,1]$ and
\begin{equation}\label{main-c1-perturbation}
\sup_{\substack{\frac12\le t\le1\\ \theta\in\Theta}}
\bigl\|\widetilde\phi(t,\cdot,\theta)
-\phi(t,\cdot,\theta)\bigr\|_{C^1(\mathcal W_0(\theta),X)}
<\epsilon_0.
\end{equation}
Then, for any $\theta\in\Theta$ and any distinct
$x,y\in\mathcal V_0(\theta)$ whose trajectories satisfy
\[
\widetilde\phi(t,x,\theta),\ 
\widetilde\phi(t,y,\theta)
\in\mathcal V_0(\theta\cdot t)
\qquad\text{for all }t\ge0,
\]
the following assertions hold.

\begin{itemize}
\item[\rm(1)] Exactly one of the following alternatives occurs:
\begin{itemize}
\item[\rm(a)] The two trajectories meet at a finite time and coincide
thereafter.

\item[\rm(b)] There exists $T^*>0$ such that
\begin{equation}\label{exponential-decay}
\widetilde\phi(t,x,\theta)-\widetilde\phi(t,y,\theta)
\notin\mathcal C_{N_0},
\qquad t\ge T^*,
\end{equation}
and their difference converges to zero exponentially as $t\to+\infty$.

\item[\rm(c)] There exist $T^*>0$ and
$i\in\{1,\dots,N_0\}$ such that
\begin{equation}\label{stayin-cone}
\widetilde\phi(t,x,\theta)-\widetilde\phi(t,y,\theta)
\in\operatorname{int}(\mathcal C_i)\setminus\mathcal C_{i-1},
\qquad t\ge T^*.
\end{equation}
\end{itemize}

\item[\rm(2)] Suppose, in addition, that both trajectories admit complete
extensions in $\mathcal V_0$, denoted by the same symbols, and that their
negative semi-orbits have compact closures contained in $\mathcal V_0$.
If
\[
\|\widetilde\phi(t,x,\theta)-\widetilde\phi(t,y,\theta)\|
\nrightarrow0
\qquad\text{as }t\to-\infty,
\]
then alternative {\rm(c)} in part {\rm(1)} holds. Moreover, there exist
$T^{**}>0$ and integers $1\le i\le i'\le N_0$ such that
\begin{equation}\label{stayin-cone-negative}
\widetilde\phi(t,x,\theta)-\widetilde\phi(t,y,\theta)
\in\operatorname{int}(\mathcal C_{i'})\setminus\mathcal C_{i'-1},
\qquad t\le-T^{**}.
\end{equation}

\item[\rm(3)] Assume, in addition, that {\rm\textbf{(H5)}} holds. If
\[
\widetilde\phi(t_0,x,\theta)-\widetilde\phi(t_0,y,\theta)
\in\mathcal C_{N_0}\setminus\{0\}
\]
for some $t_0>0$, then
\[
\widetilde\phi(t,x,\theta)-\widetilde\phi(t,y,\theta)
\notin\mathcal H
\qquad\text{for all sufficiently large }t.
\]
\end{itemize}
\end{theorem}

We first specialize Theorem~\ref{prop-orbit} to compact invariant sets.
For distinct points in the same fibre, compactness and invariance exclude
alternatives (a) and (b) of Theorem~\ref{prop-orbit}(1).

\begin{corollary}\label{constant-invar-set}
Assume the hypotheses and perturbation setting of
Theorem~\ref{prop-orbit}. Let $E\subset\mathcal V_0$ be a compact invariant
set of $\widetilde\Lambda_t$. Then the restriction of
$\widetilde\Lambda_t$ to $E$ extends uniquely to a flow. Moreover, for any
two distinct points $(x,\theta),(y,\theta)\in E$, there exist $T^*>0$ and
$i\in\{1,\dots,N_0\}$ such that
\[
\widetilde\phi(t,x,\theta)-\widetilde\phi(t,y,\theta)
\in\operatorname{int}(\mathcal C_i)\setminus\mathcal C_{i-1},
\qquad t\ge T^*.
\]
If, in addition,
\[
\|\widetilde\phi(t,x,\theta)-\widetilde\phi(t,y,\theta)\|
\nrightarrow0
\qquad\text{as }t\to-\infty,
\]
where negative times are understood with respect to this flow extension,
then there also exist $T^{**}>0$ and an integer $i'$ with
$i\le i'\le N_0$ such that
\[
\widetilde\phi(t,x,\theta)-\widetilde\phi(t,y,\theta)
\in\operatorname{int}(\mathcal C_{i'})\setminus\mathcal C_{i'-1},
\qquad t\le-T^{**}.
\]
\end{corollary}

\begin{remark}\label{rk-cor21}
{\rm
Corollary~\ref{constant-invar-set} shows that compact invariance reduces
the alternatives in Theorem~\ref{prop-orbit} to eventual cone-layer
separation for distinct points in the same fibre. The relation $i\le i'$
compares the eventual layer indices in the two time directions: the
forward index does not exceed the backward index. This is an asymptotic
comparison and does not assert monotonicity of the layer index at every
time. The resulting separation is used in the finite-dimensional reduction
and the codimension-one analysis below.
}
\end{remark}

The next corollary combines the transversality conclusion of
Theorem~\ref{prop-orbit}(3) with the assumption that each nonzero orbit
difference remains in the same cone layer for all time. Together, these
conditions ensure that the quotient projection separates distinct points
in each fibre.

\begin{corollary}\label{imbeding-d-plane}
Assume the hypotheses and perturbation setting of
Theorem~\ref{prop-orbit}, and assume {\rm\textbf{(H5)}}. Let
$E\subset\mathcal V_0$ be a compact invariant set of
$\widetilde\Lambda_t$. By Corollary~\ref{constant-invar-set}, the
restriction of $\widetilde\Lambda_t$ to $E$ is a flow, and we use the same
notation for its extension to negative times. Assume that, for any two
distinct points $(x,\theta),(y,\theta)\in E$, there exists
$i\in\{1,\dots,N_0\}$ such that
\[
\widetilde\phi(t,x,\theta)-\widetilde\phi(t,y,\theta)
\in\operatorname{int}(\mathcal C_i)\setminus\mathcal C_{i-1}
\qquad\text{for all }t\in\RR.
\]
Then the skew-product flow on $E$ is topologically conjugate to a
skew-product flow on a compact invariant set
$\widehat E\subset\RR^d\times\Theta$.
\end{corollary}

\begin{remark}\label{rk-cor22}
{\rm
The finite-dimensional reduction in
Corollary~\ref{imbeding-d-plane} should be understood as a reduction in the
fibre variables. Since $\mathcal H$ has codimension $d$, the quotient
$X/\mathcal H$ is a $d$-dimensional space and hence isomorphic to $\RR^d$.
The hypotheses of the corollary ensure that the quotient map is injective
on each fibre of $E$. Consequently, the dynamics in each fibre can be
represented by $d$ real coordinates, while the base flow on $\Theta$ is
retained unchanged. In particular, no finite-dimensionality assumption is
imposed on the base space $\Theta$.
}
\end{remark}

We now specialize to the codimension-one case. The resulting fibre order
leads to a precise description of minimal sets and $\omega$-limit sets.

\begin{theorem}[Almost $1$-covering property of minimal sets]
\label{a-a-one}
Assume the hypotheses and perturbation setting of
Theorem~\ref{prop-orbit}, assume {\rm\textbf{(H5)}}, and suppose that
$\mathcal H$ has codimension one in $X$. Then every compact minimal set
$\widetilde{\mathcal M}\subset\mathcal V_0$ of $\widetilde\Lambda_t$ is
an almost $1$-cover of the base flow. More precisely, there exists a
residual subset $\Theta_0\subset\Theta$ such that, for every
$\theta\in\Theta_0$, the fibre
\[
P^{-1}(\theta)\cap\widetilde{\mathcal M}
\]
consists of a single point.
\end{theorem}

\begin{theorem}[Trichotomy for
\texorpdfstring{$\omega$}{omega}-limit sets in the codimension-one case]
\label{perturbation-thm}
Assume the hypotheses and perturbation setting of
Theorem~\ref{prop-orbit}, together with {\rm\textbf{(H5)}}, and
suppose that $\mathcal H$ has codimension one in $X$.
Let $(x_0,\theta_0)\in\mathcal V_0$ have a precompact positive orbit,
and assume that
\[
\widetilde\Omega
:=\widetilde\omega(x_0,\theta_0)
\subset\mathcal V_0.
\]
Then $\widetilde\Omega$ contains at most two minimal sets, and
exactly one of the following alternatives holds:
\begin{itemize}
\item[\rm(1)]
$\widetilde\Omega$ is a minimal set.

\item[\rm(2)]
$\widetilde\Omega$ contains exactly one minimal set
$\widetilde{\mathcal M}_1$, and
\[
\widetilde{\mathcal M}_{11}
:=\widetilde\Omega\setminus\widetilde{\mathcal M}_1
\ne\emptyset.
\]
For every $(x_{11},\theta)\in\widetilde{\mathcal M}_{11}$,
\[
\widetilde{\mathcal M}_1
\subset\widetilde\omega(x_{11},\theta)
\quad\text{and}\quad
\widetilde{\mathcal M}_1
\subset\widetilde\alpha(x_{11},\theta).
\]

\item[\rm(3)]
$\widetilde\Omega$ contains exactly two distinct minimal sets
$\widetilde{\mathcal M}_1$ and $\widetilde{\mathcal M}_2$, and
\[
\widetilde{\mathcal M}_{12}
:=\widetilde\Omega\setminus
\bigl(\widetilde{\mathcal M}_1
\cup\widetilde{\mathcal M}_2\bigr)
\ne\emptyset.
\]
For every $(x_{12},\theta)\in\widetilde{\mathcal M}_{12}$,
exactly one of the following holds:
\[
\widetilde{\mathcal M}_1
\subset\widetilde\omega(x_{12},\theta),
\qquad
\widetilde{\mathcal M}_2
\cap\widetilde\omega(x_{12},\theta)=\emptyset;
\]
\[
\widetilde{\mathcal M}_2
\subset\widetilde\omega(x_{12},\theta),
\qquad
\widetilde{\mathcal M}_1
\cap\widetilde\omega(x_{12},\theta)=\emptyset;
\]
or
\[
\widetilde{\mathcal M}_1
\cup\widetilde{\mathcal M}_2
\subset\widetilde\omega(x_{12},\theta).
\]
The same three alternatives hold with
$\widetilde\omega(x_{12},\theta)$ replaced by
$\widetilde\alpha(x_{12},\theta)$.
\end{itemize}
\end{theorem}

\begin{remark}\label{new-rk2}
{\rm
Theorems~\ref{a-a-one} and \ref{perturbation-thm} provide the
nonautonomous counterparts of the convergence conclusions familiar from
autonomous and time-periodic systems. Here an almost $1$-cover replaces an
equilibrium or a periodic orbit as the natural minimal recurrent object. An
$\omega$-limit set need not itself be minimal, but it can contain at most
two minimal sets. Theorem~\ref{perturbation-thm} specifies which minimal
sets are contained in the forward and backward limit sets of the remaining
points; these containments do not, by themselves, imply that the remaining
orbits are heteroclinic between minimal sets. Taking
$\widetilde\Lambda_t=\Lambda_t$ recovers the corresponding conclusions for
the unperturbed system.
}
\end{remark}

\section{Proofs of the main results}\label{main-proof}

This section is divided into two parts. In Section~\ref{subsec-proof-thm21}, we derive
Theorem~\ref{prop-orbit} and Corollaries~\ref{constant-invar-set}--\ref{imbeding-d-plane}
from a key perturbation proposition. In Section~\ref{subsec-proof-codim1}, we specialize to
the codimension-one case, where the cone structure induces an order on each fibre and leads to
the proofs of Theorems~\ref{a-a-one} and~\ref{perturbation-thm}.

\subsection{Proof of Theorem~\ref{prop-orbit} and its corollaries}
\label{subsec-proof-thm21}

We begin with a key perturbation proposition, adapting the cone estimates
of \cite[Propositions~5.1 and~6.1]{Te} to the present skew-product setting.
It provides the local cone-entry, contraction, and
transversality estimates from which Theorem~\ref{prop-orbit} and its corollaries follow.

\begin{proposition}\label{pertur-pro}
Suppose that {\rm\textbf{(H1)}--\textbf{(H4)}} hold and that $K_0$ is as in
Theorem~\ref{prop-orbit}. Then there exist
$N_0\in\NN$, numbers
$0<\lambda_0<1$, $T_0>0$, $\epsilon_0>0$, open bounded neighborhoods
$\mathcal V_0$ and $\mathcal W_0$ of $K_0$ satisfying
$\overline{\mathcal V_0}\subset\mathcal W_0\subset\mathcal U_0$, and a sequence of $k_i$-cones
$\{\mathcal{C}_i\}_{i=1}^{N_0}$ such that, with
$\mathcal C_0=\{0\}$ and $\mathcal V_0(\theta)=\{x:(x,\theta)\in\mathcal V_0\}$,
\[
\mathcal C_i\setminus\{0\}\subset\operatorname{int}\mathcal C_{i+1},
\qquad 1\le i<N_0.
\]
Moreover, for any skew-product semiflow
\[
\tilde\Lambda_t=(\tilde\phi(t,\cdot,\theta),\theta\cdot t),\qquad t\ge0,
\]
on $\mathcal{U}_0$, with
$\widetilde\phi(t,\cdot,\theta)\in C^1(\mathcal W_0(\theta),X)$ for
$t\in[1/2,1]$, satisfying
\begin{equation}\label{perturb-c1}
\sup_{\substack{\frac12\le t\le1\\ \theta\in\Theta}}
\bigl\|\widetilde\phi(t,\cdot,\theta)
-\phi(t,\cdot,\theta)\bigr\|_{C^1(\mathcal W_0(\theta),X)}
<\epsilon_0,
\end{equation}
the following assertions hold for any $\theta\in\Theta$ and any $x,y\in\mathcal{V}_0(\theta)$:

\begin{itemize}
\item[\rm(1)] If $x-y\in \mathcal{C}_i\setminus\{0\}$ for some $1\le i\le N_0$, and
\[
\tilde\phi(s,x,\theta),\ \tilde\phi(s,y,\theta)\in \mathcal{V}_0(\theta\cdot s)
\quad \text{for all }s\in[0,t]
\]
for some $t\ge T_0$, then
\[
\tilde\phi(t,x,\theta)-\tilde\phi(t,y,\theta)\in \operatorname{int}(\mathcal{C}_i)\setminus\{0\}.
\]

\item[\rm(2)] If
\[
\tilde\phi(s,x,\theta),\ \tilde\phi(s,y,\theta)\in \mathcal{V}_0(\theta\cdot s)
\quad \text{for all }s\in[0,t],
\]
and $\tilde\phi(t,x,\theta)-\tilde\phi(t,y,\theta)\notin \mathcal{C}_{N_0}$ for some $t\ge T_0$,
then
\[
\|\tilde\phi(t,x,\theta)-\tilde\phi(t,y,\theta)\|\le \lambda_0^t\|x-y\|.
\]

\item[\rm(3)] Assume in addition that {\rm\textbf{(H5)}} holds. If
\[
\tilde\phi(s,x,\theta),\ \tilde\phi(s,y,\theta)\in \mathcal{V}_0(\theta\cdot s)
\quad \text{for all }0\le s\le 2T_0,
\]
and
\[
x-y,\ \tilde\phi(2T_0,x,\theta)-\tilde\phi(2T_0,y,\theta)\in
\mathcal{C}_i\setminus \mathcal{C}_{i-1}
\]
for some $i\in\{1,\dots,N_0\}$, then
\[
\tilde\phi(T_0,x,\theta)-\tilde\phi(T_0,y,\theta)\notin \mathcal{H}.
\]
\end{itemize}
\end{proposition}

We now derive Theorem~\ref{prop-orbit} from Proposition~\ref{pertur-pro}. The detailed proof of
Proposition~\ref{pertur-pro} will be given in Section~\ref{proof-section}.

\begin{proof}[Proof of Theorem~\ref{prop-orbit}]
We divide the proof into three steps corresponding to the three parts of the theorem.

\medskip
\noindent\textbf{Step 1. Proof of part {\rm(1)}.}
If the two trajectories meet at some finite time, the semiflow property shows
that they agree thereafter, and the first alternative holds. We may therefore
assume throughout the rest of this step that their difference never vanishes.

If \eqref{exponential-decay} holds, then Proposition~\ref{pertur-pro}(2) immediately yields that
\[
\|\tilde\phi(t,x,\theta)-\tilde\phi(t,y,\theta)\|
\]
decays to zero exponentially fast as $t\to+\infty$.

Assume now that \eqref{exponential-decay} does not hold. Then the
difference enters $\mathcal C_{N_0}$ at some nonnegative time. Let
$i\in\{1,\dots,N_0\}$ be the smallest index for which there exists
$t_i\ge0$ such that
\[
\tilde\phi(t_i,x,\theta)-\tilde\phi(t_i,y,\theta)\in\mathcal C_i.
\]
Proposition~\ref{pertur-pro}(1) places the difference in
$\operatorname{int}\mathcal C_i$ for all $t\ge t_i+T_0$. By the minimal
choice of $i$, it never belongs to $\mathcal C_{i-1}$; for $i=1$, this
follows from the assumed nonvanishing of the difference. Thus
\eqref{stayin-cone} holds with $T^*=t_i+T_0$. The three alternatives are
mutually exclusive, which proves part {\rm(1)}.

\medskip
\noindent\textbf{Step 2. Proof of part {\rm(2)}.}
We first rule out finite-time coincidence. Suppose that the two trajectories
meet at a time $t_c>0$, and fix $r<\min\{0,t_c-T_0\}$. For every $s\le r$,
their difference at time $s$ is nonzero because $x\ne y$ at time zero. It
cannot belong to $\mathcal C_{N_0}$, since
Proposition~\ref{pertur-pro}(1), applied from $s$ to $t_c$, would make the
difference at $t_c$ a nonzero element of
$\operatorname{int}\mathcal C_{N_0}$. The same argument applies at time
$r$. For $s\le r-T_0$, Proposition~\ref{pertur-pro}(2), applied from
$s$ to $r$, therefore gives
\[
\|\tilde\phi(r,x,\theta)-\tilde\phi(r,y,\theta)\|
\le \lambda_0^{r-s}
\|\tilde\phi(s,x,\theta)-\tilde\phi(s,y,\theta)\|.
\]
The right-hand side tends to zero as $s\to-\infty$ by precompactness of the
negative semi-orbits. Thus the two trajectories already agree at time
$r<0$, and hence at time zero, a contradiction.

Suppose, by contradiction, that \eqref{exponential-decay} holds. Then the negative semi-orbits of
$(x,\theta)$ and $(y,\theta)$ are precompact in $\mathcal{V}_0$, and
\[
\|\tilde\phi(t,x,\theta)-\tilde\phi(t,y,\theta)\|\nrightarrow0
\quad\text{as }t\to-\infty.
\]
We claim that \eqref{exponential-decay} can be strengthened to
\begin{equation}\label{exponential-decay2}
\tilde\phi(t,x,\theta)-\tilde\phi(t,y,\theta)\notin \mathcal{C}_{N_0}
\quad\text{for all }t\in\RR.
\end{equation}
Indeed, if there existed some $T_1<T^*$ such that
\[
\tilde\phi(T_1,x,\theta)-\tilde\phi(T_1,y,\theta)\in \mathcal{C}_{N_0},
\]
then Proposition~\ref{pertur-pro}(1) would imply
\[
\tilde\phi(t,x,\theta)-\tilde\phi(t,y,\theta)\in \operatorname{int}\mathcal{C}_{N_0}
\quad\text{for all }t\ge T_1+T_0,
\]
contradicting \eqref{exponential-decay}. Hence \eqref{exponential-decay2} holds.

By compactness of the negative semi-orbits, there exists $c_0>0$ such that
\[
\|\tilde\phi(-t,x,\theta)-\tilde\phi(-t,y,\theta)\|\le 2c_0
\quad\text{for all }t\ge0.
\]
Applying Proposition~\ref{pertur-pro}(2) to the pair
\[
\bigl(\tilde\phi(-t,x,\theta),\tilde\phi(-t,y,\theta),\theta\cdot(-t)\bigr)
\]
and using \eqref{exponential-decay2}, we get
\[
\|x-y\|
\le \lambda_0^t\,\|\tilde\phi(-t,x,\theta)-\tilde\phi(-t,y,\theta)\|
\le 2c_0\lambda_0^t
\qquad\text{for all }t\ge T_0.
\]
Letting $t\to+\infty$ yields $x=y$, a contradiction. Therefore the contracting alternative
\eqref{exponential-decay} cannot occur, and only \eqref{stayin-cone} is possible.

It remains to prove \eqref{stayin-cone-negative}. Since the distance of the two
negative semi-orbits does not converge to zero, there are $c_*>0$ and a
sequence $t_n\to-\infty$ such that
\[
\|\widetilde\phi(t_n,x,\theta)-\widetilde\phi(t_n,y,\theta)\|\ge c_*.
\]
By precompactness, after passing to a subsequence we have
\[
\tilde\Lambda_{t_n}(x,\theta)\to(x^*,\theta^*),\qquad
\tilde\Lambda_{t_n}(y,\theta)\to(y^*,\theta^*)
\]
for some $(x^*,\theta^*),(y^*,\theta^*)\in\mathcal{V}_0$, with
$\|x^*-y^*\|\ge c_*$ and hence $x^*\ne y^*$.
The limiting points lie in the compact invariant alpha-limit sets of the
original complete trajectories, which are contained in $\mathcal V_0$.
They therefore admit complete extensions with compact negative semi-orbits.
The preceding arguments excluding finite-time coincidence and the
contracting alternative use only boundedness of these negative semi-orbits,
not their backward nonconvergence. They apply to the limiting pair as well.
Applying part {\rm(1)} to
$(x^*,\theta^*)$ and $(y^*,\theta^*)$, we find some $j\in\{1,\dots,N_0\}$ and some $T_2^*>0$
such that
\[
\tilde\phi(t,x^*,\theta^*)-\tilde\phi(t,y^*,\theta^*)
\in \operatorname{int}(\mathcal{C}_{j})\setminus \mathcal{C}_{j-1}
\quad\text{for all }t\ge T_2^*.
\]
Since $\operatorname{int}(\mathcal{C}_{j})\setminus \mathcal{C}_{j-1}$ is open, there exists
$N_1\in\NN$ such that $T_2^*+t_{N_1}-T_0<0$ and
\begin{equation}\label{negative-sequence}
\tilde\phi(T_2^*+t_n,x,\theta)-\tilde\phi(T_2^*+t_n,y,\theta)
\in \operatorname{int}(\mathcal{C}_{j})\setminus \mathcal{C}_{j-1}
\quad\text{for all }n\ge N_1.
\end{equation}

Fix any $t\in\RR$. Choose $n\ge N_1$ so large that $t_n+T_2^*+T_0<t$. Then
Proposition~\ref{pertur-pro}(1), applied on the interval $[t_n+T_2^*,\,t]$, yields
\begin{equation}\label{stayin-i'}
\tilde\phi(t,x,\theta)-\tilde\phi(t,y,\theta)\in \operatorname{int}\mathcal{C}_{j}.
\end{equation}
If $j=1$, the difference never belongs to $\mathcal C_0=\{0\}$, as
already shown. For $j\ge2$, suppose that there existed some
$t_0\le T_2^*+t_{N_1}-T_0$ such that
\[
\tilde\phi(t_0,x,\theta)-\tilde\phi(t_0,y,\theta)\in \mathcal{C}_{j-1}.
\]
Applying Proposition~\ref{pertur-pro}(1) on $[t_0,\,T_2^*+t_{N_1}]$ would then imply
\[
\tilde\phi(T_2^*+t_{N_1},x,\theta)-\tilde\phi(T_2^*+t_{N_1},y,\theta)
\in \operatorname{int}\mathcal{C}_{j-1},
\]
contradicting \eqref{negative-sequence}. Hence
\[
\tilde\phi(t,x,\theta)-\tilde\phi(t,y,\theta)\notin \mathcal{C}_{j-1}
\quad\text{for all }t\le T_2^*+t_{N_1}-T_0.
\]
Combining this with \eqref{stayin-i'} gives
\[
\tilde\phi(t,x,\theta)-\tilde\phi(t,y,\theta)\in
\operatorname{int}(\mathcal{C}_{j})\setminus \mathcal{C}_{j-1}
\quad\text{for all }t\le T_2^*+t_{N_1}-T_0.
\]
This is exactly \eqref{stayin-cone-negative}, with
\[
T^{**}=-(T_2^*+t_{N_1}-T_0)>0.
\]
Finally, let $i$ be the forward layer index in \eqref{stayin-cone}.  Since
\eqref{stayin-i'} holds for every $t\in\RR$, every sufficiently large positive
time belongs both to $\operatorname{int}\mathcal C_j$ and to
$\operatorname{int}(\mathcal C_i)\setminus\mathcal C_{i-1}$.  If $j<i$, then
$\mathcal C_j\subset\mathcal C_{i-1}$, a contradiction.  Hence $i\le j$, as
asserted in the theorem.
Thus part {\rm(2)} is proved.

\medskip
\noindent\textbf{Step 3. Proof of part {\rm(3)}.}
The assumption at time $t_0$, together with
Proposition~\ref{pertur-pro}(1), excludes both finite-time coincidence and the contracting alternative. Hence part {\rm(1)} gives some $i_0\in\{1,\dots,N_0\}$ and some $T^*>0$ such that
\[
\tilde\phi(t,x,\theta)-\tilde\phi(t,y,\theta)\in
\operatorname{int}(\mathcal{C}_{i_0})\setminus \mathcal{C}_{i_0-1}
\quad\text{for all }t\ge T^*.
\]
Enlarging $T^*$ if necessary, we may assume $T^*\ge t_0+T_0$. Fix any $t\ge T^*+T_0$. Then
both
\[
\tilde\phi(t-T_0,x,\theta)-\tilde\phi(t-T_0,y,\theta)
\quad\text{and}\quad
\tilde\phi(t+T_0,x,\theta)-\tilde\phi(t+T_0,y,\theta)
\]
belong to $\mathcal{C}_{i_0}\setminus \mathcal{C}_{i_0-1}$. Applying
Proposition~\ref{pertur-pro}(3) to the shifted pair
\[
\bigl(\tilde\phi(t-T_0,x,\theta),\tilde\phi(t-T_0,y,\theta),\theta\cdot(t-T_0)\bigr),
\]
we obtain
\[
\tilde\phi(t,x,\theta)-\tilde\phi(t,y,\theta)\notin \mathcal H.
\]
Since $t\ge T^*+T_0$ is arbitrary, the conclusion follows. This proves part {\rm(3)}.
\end{proof}

\begin{remark}\label{rk-proof-thm21}
{\rm
Theorem~\ref{prop-orbit} should be viewed as the global consequence of
Proposition~\ref{pertur-pro}. Part {\rm(1)} gives three forward alternatives
for orbit differences; part {\rm(2)} excludes finite-time coincidence and
the contracting alternative for the distinct precompact entire orbits under consideration;
and part {\rm(3)} excludes asymptotic motion inside the codimension-$d$ subspace
$\mathcal H$. This last point is the key ingredient in the finite-dimensional reduction below.
}
\end{remark}

The first corollary shows that, on a compact invariant set, the contracting alternative in
Theorem~\ref{prop-orbit} cannot occur.

\begin{proof}[Proof of Corollary~\ref{constant-invar-set}]
We first prove that every pair of distinct points in the same fibre of $E$
has its difference in the terminal cone. Let
$(x,\theta),(y,\theta)\in E$ be distinct. Since
$\widetilde\Lambda_s(E)=E$, for every $s>0$ there exist
$(x_{-s},\theta\cdot(-s)),(y_{-s},\theta\cdot(-s))\in E$ whose time-$s$
images are $(x,\theta)$ and $(y,\theta)$, respectively. If
$x-y\notin\mathcal C_{N_0}$, Proposition~\ref{pertur-pro}(2), applied to
these two preimages for $s\ge T_0$, gives
\[
 \|x-y\|\le \lambda_0^s
 \|x_{-s}-y_{-s}\|.
\]
The last factor is uniformly bounded because $E$ is compact. Letting
$s\to\infty$ yields $x=y$, a contradiction. Therefore
\begin{equation}\label{compact-set-terminal-cone}
x-y\in\mathcal C_{N_0}
\end{equation}
for every two distinct same-fibre points of $E$.

We next prove injectivity. Suppose that two distinct points
$(x,\theta),(y,\theta)\in E$ had the same image at some time $t>0$.
By \eqref{compact-set-terminal-cone},
$x-y\in\mathcal C_{N_0}\setminus\{0\}$. Choose
$\tau\ge\max\{t,T_0\}$. The semiflow property makes their time-$\tau$
images equal, whereas Proposition~\ref{pertur-pro}(1), applied on
$[0,\tau]$, says that their difference at time $\tau$ is a nonzero element
of $\operatorname{int}\mathcal C_{N_0}$, a contradiction. Hence every
positive-time map is injective. Invariance makes it onto $E$, so it is a
homeomorphism of the compact space $E$. The inverses are compatible with
the semigroup law. Their joint continuity in negative time follows from
continuity of the positive-time maps and compactness of $E$ (equivalently,
from the standard compact-space inverse-flow lemma). Thus the restricted
semiflow extends uniquely to a continuous flow.

For the original distinct pair, finite-time coincidence is now impossible,
and \eqref{compact-set-terminal-cone} together with
Proposition~\ref{pertur-pro}(1) rules out the contracting alternative in
Theorem~\ref{prop-orbit}(1). Thus the eventual cone-layer alternative holds.
If the additional backward nonconvergence hypothesis is satisfied,
Theorem~\ref{prop-orbit}(2) gives the stated backward layer and its index
relation.
\end{proof}

The next corollary uses the transversality conclusion in
Theorem~\ref{prop-orbit}(3) to realize the dynamics on compact invariant sets
in finite dimensions. The stronger assumption below guarantees injectivity
of the projection along $\mathcal H$ on the whole invariant set.

\begin{proof}[Proof of Corollary~\ref{imbeding-d-plane}]
Fix two distinct points $(x,\theta),(y,\theta)\in E$ and a time $s\in\RR$.
The hypothesis of the corollary shows that the differences at times
$s-T_0$ and $s+T_0$ lie in the same layer.  Applying
Proposition~\ref{pertur-pro}(3) to the common backward translate at time
$s-T_0$ gives
\[
\tilde\phi(s,x,\theta)-\tilde\phi(s,y,\theta)\notin \mathcal{H}.
\]
Since $s$ was arbitrary, the difference avoids $\mathcal H$ for every
$s\in\RR$.
Let $\mathcal{H}'$ be a $d$-dimensional subspace of $X$ such that
\[
X=\mathcal{H}'\oplus \mathcal{H},
\]
and let
\[
p:X\to \mathcal{H}'
\]
be the projection onto $\mathcal{H}'$ along $\mathcal{H}$. Define
\[
h:X\times\Theta\to \mathcal{H}'\times\Theta,\qquad h(x,\theta)=(p(x),\theta).
\]
Since $E$ is compact and $h$ is continuous and injective on $E$, the restriction
\[
h|_E:E\to \hat E:=h(E)
\]
is a homeomorphism.

By Corollary~\ref{constant-invar-set}, the restriction of
$\tilde\Lambda_t$ to $E$ extends uniquely to a flow. For any
$(\hat x,\theta)\in\hat E$ and $t\in\RR$, define
\[
\hat\phi(t,\hat x,\theta):=p\bigl(\tilde\phi(t,x,\theta)\bigr),
\]
where $x\in X$ is chosen so that $(x,\theta)\in E$ and $p(x)=\hat x$. This is well-defined by
injectivity of $h|_E$. We then obtain a skew-product flow
\begin{equation}\label{induced-skew-product}
\hat\Lambda_t:\hat E\to \hat E,\qquad
\hat\Lambda_t(\hat x,\theta)=(\hat\phi(t,\hat x,\theta),\theta\cdot t).
\end{equation}
Moreover,
\[
h\circ \tilde\Lambda_t(x,\theta)=\hat\Lambda_t\circ h(x,\theta)
\quad\text{for all }t\in\RR,\ (x,\theta)\in E,
\]
so the following diagram commutes:
\[
\begin{CD}
E @>h>> \hat E\\
@VV\tilde\Lambda_tV @VV\hat\Lambda_tV\\
E @>h>> \hat E.
\end{CD}
\]

Equivalently,
$\hat\Lambda_t=h|_E\circ\tilde\Lambda_t|_E\circ(h|_E)^{-1}$, so continuity
and the flow identity follow directly from those of $\tilde\Lambda_t|_E$.
Identifying $\mathcal H'$ with $\RR^d$ proves the corollary.
\end{proof}

\begin{remark}\label{rk-embedding}
{\rm
The realization above is defined on $\hat E=h(E)$. It does not assert
that the induced flow extends to all of $\RR^d\times\Theta$.
}
\end{remark}

\subsection{Proofs of Theorems~\ref{a-a-one} and~\ref{perturbation-thm}}
\label{subsec-proof-codim1}

We now specialize to the codimension-one case. In this setting, the complement of
$\mathcal H$ splits into two ordered half-spaces, and the perturbative cone structure induces an
eventual order on each fibre. This order is the basis for the proofs of
Theorems~\ref{a-a-one} and~\ref{perturbation-thm}.

Throughout this subsection, we assume that $E\subset\mathcal{V}_0$ is a compact invariant set
and that $\tilde{\mathcal M}\subset\mathcal{V}_0$ is a compact minimal invariant set of the
skew-product semiflow $\tilde\Lambda_t$.

By Corollary~\ref{constant-invar-set}, the restriction to any such compact
invariant set extends uniquely to a continuous flow. We use this extension
whenever negative times occur below.

\subsubsection{Proof of Theorem~\ref{a-a-one}}

We first establish several auxiliary lemmas.

\begin{lemma}\label{sequence-limit}
Let $E\subset\mathcal V_0$ be a compact invariant set, and let
$(x_i,\theta),(x_i^*,\theta^*)\in E$ for $i=1,2$, with
$x_1\ne x_2$ and $x_1^*\ne x_2^*$. If there exists a sequence
$t_n\to+\infty$ such that
\[
\widetilde\Lambda_{t_n}(x_i,\theta)\longrightarrow(x_i^*,\theta^*)
\qquad(i=1,2),
\]
then there exists $j\in\{1,\ldots,N_0\}$ such that
\[
\widetilde\phi(t,x_1^*,\theta^*)-\widetilde\phi(t,x_2^*,\theta^*)
\in\operatorname{int}(\mathcal C_j)\setminus\mathcal C_{j-1}
\qquad\text{for every }t\in\RR.
\]
The index $j$ is the eventual forward cone-layer index of the original
pair. The same conclusion holds for $t_n\to-\infty$ if the original
orbit difference belongs to
$\operatorname{int}(\mathcal C_j)\setminus\mathcal C_{j-1}$ for all
sufficiently negative times.
\end{lemma}

\begin{proof}
By Corollary~\ref{constant-invar-set}, the restriction to $E$ is a
continuous flow, and there are $j\in\{1,\ldots,N_0\}$ and $T>0$ such that
\begin{equation}\label{eventual-layer-original-pair}
\widetilde\phi(t,x_1,\theta)-\widetilde\phi(t,x_2,\theta)
\in\operatorname{int}(\mathcal C_j)\setminus\mathcal C_{j-1},
\qquad t\ge T.
\end{equation}
For a sequence tending to $-\infty$, the corresponding membership for
$t\le-T$ is assumed. In either case, for each fixed $s\in\RR$,
continuity of the flow gives
\[
\widetilde\Lambda_{t_n+s}(x_i,\theta)
\longrightarrow\widetilde\Lambda_s(x_i^*,\theta^*),\qquad i=1,2.
\]
Since $t_n+s$ eventually lies in the relevant time interval, closedness
of $\mathcal C_j$ implies
$d^*(s):=\widetilde\phi(s,x_1^*,\theta^*)-
\widetilde\phi(s,x_2^*,\theta^*)\in\mathcal C_j$.
Injectivity of the flow gives $d^*(s)\ne0$, and
Proposition~\ref{pertur-pro}(1), applied from $s-T_0$ to $s$, gives
$d^*(s)\in\operatorname{int}\mathcal C_j$.

If $j\ge2$ and $d^*(s_0)\in\mathcal C_{j-1}$, the same proposition
gives $d^*(s_0+T_0)\in\operatorname{int}\mathcal C_{j-1}$.
Convergence then places the original difference at $t_n+s_0+T_0$ in
$\operatorname{int}\mathcal C_{j-1}$ for all sufficiently large $n$,
contradicting its prescribed layer. If $j=1$, nonvanishing already
excludes $\mathcal C_0$. This proves both assertions.
\end{proof}

Since $\mathcal{H}$ has codimension one in $X$, there exists a one-dimensional subspace
\[
\mathcal{H}'=\operatorname{span}\{u^+\}
\]
with $u^+\in X\setminus\mathcal H$ and $\|u^+\|=1$ such that
\[
X=\mathcal{H}'\oplus\mathcal{H}.
\]
Define the two connected components of $X\setminus\mathcal H$ by
\[
X^+ := \{su^+ + u \mid s>0,\,u\in\mathcal{H}\},\qquad
X^- := \{su^+ + u \mid s<0,\,u\in\mathcal{H}\}.
\]

The next lemma isolates the codimension-one strengthening of
Proposition~\ref{pertur-pro}(3) that will be used repeatedly.

\begin{lemma}\label{halfspace-entry}
Assume that $\mathcal H$ has codimension one, that
$x,y\in\mathcal V_0(\theta)$, and that
$\tilde\phi(s,x,\theta),\tilde\phi(s,y,\theta)\in\mathcal V_0(\theta\cdot s)$
for $0\le s\le2T_0$. Put
$\Delta(t)=\tilde\phi(t,x,\theta)-\tilde\phi(t,y,\theta)$.
If
\[
\Delta(0),\ \Delta(2T_0)\in \mathcal{C}_i\setminus \mathcal{C}_{i-1}
\]
for some $i\in\{1,\dots,N_0\}$, then
\[
\Delta(T_0)\in X^+\cup X^-.
\]
\end{lemma}

\begin{proof}
Proposition~\ref{pertur-pro}(3) gives $\Delta(T_0)\notin\mathcal H$.
The conclusion follows from $X\setminus\mathcal H=X^+\cup X^-$.
\end{proof}

\begin{lemma}\label{liein-halfspace}
Let $(x,\theta),(y,\theta)\in\mathcal{V}_0$ be two distinct points such that
$\tilde\phi(t,x,\theta),\tilde\phi(t,y,\theta)\in\mathcal{V}_0(\theta\cdot t)$ for all $t>0$.
Assume that the two trajectories do not meet at a finite time and that
\eqref{exponential-decay} does not occur. Then there exists $T>0$ such that
\[
\tilde\phi(t,x,\theta)-\tilde\phi(t,y,\theta)\in X^+
\quad\text{for all }t\ge T,
\]
or
\[
\tilde\phi(t,x,\theta)-\tilde\phi(t,y,\theta)\in X^-
\quad\text{for all }t\ge T.
\]
\end{lemma}

\begin{proof}
Theorem~\ref{prop-orbit}(1) gives an index $i\in\{1,\ldots,N_0\}$
and a time $T^*\ge0$ such that
\begin{equation}\label{stay-in-cone}
\tilde\phi(t,x,\theta)-\tilde\phi(t,y,\theta)
\in\operatorname{int}(\mathcal C_i)\setminus\mathcal C_{i-1},
\qquad t\ge T^*.
\end{equation}
For each $t\ge T^*+T_0$, apply Lemma~\ref{halfspace-entry} to the
orbit segment from $t-T_0$ to $t+T_0$. The difference at time $t$
therefore lies outside $\mathcal H$. By continuity, it stays in one
fixed connected component of $X\setminus\mathcal H$ for all
$t\ge T^*+T_0$.
\end{proof}

\begin{definition}
\rm
For any $\theta\in\Theta$, define relations $\ge$ and $>$ on
$\mathcal U_0\cap P^{-1}(\theta)$ as follows: for
$(x,\theta),(y,\theta)\in\mathcal U_0\cap P^{-1}(\theta)$ we write $(x,\theta)\ge(y,\theta)$
(resp. $(x,\theta)>(y,\theta)$) if there exists $T>0$ such that
\[
\tilde\phi(t,x,\theta)-\tilde\phi(t,y,\theta)\in X^+\cup\{0\}
\quad\text{(resp. }\tilde\phi(t,x,\theta)-\tilde\phi(t,y,\theta)\in X^+\text{)}
\]
for all $t\ge T$.
\end{definition}

\begin{lemma}\label{ordered}
Assume that $E\subset\mathcal{V}_0$ is a compact invariant set of $\tilde\Lambda_t$.
For any $\theta\in\Theta$, the relation $\ge$ defines a total order on $P^{-1}(\theta)\cap E$.
\end{lemma}

\begin{proof}
By Corollary~\ref{constant-invar-set}, the restriction of
$\tilde\Lambda_t$ to $E$ extends to a flow. In particular, its
time maps are injective.

For any two distinct points in the same fibre
$P^{-1}(\theta)\cap E$, the same corollary gives an eventual
fixed cone layer for their orbit difference.
Proposition~\ref{pertur-pro}(3), applied to sufficiently late
orbit segments, then shows that this difference avoids
$\mathcal H$ for all sufficiently large times.
By continuity, it remains in one of the two components
$X^+$ or $X^-$ of $X\setminus\mathcal H$.
Thus any two points in the fibre are comparable.

Reflexivity is immediate. To verify antisymmetry, suppose that
\[
(x,\theta)\ge(y,\theta),
\qquad
(y,\theta)\ge(x,\theta).
\]
Since
\[
(X^+\cup\{0\})\cap(X^-\cup\{0\})=\{0\},
\]
the defining relations imply
\[
\tilde\phi(t,x,\theta)=\tilde\phi(t,y,\theta)
\]
for all sufficiently large $t$.
Injectivity of the time maps on $E$ therefore gives $x=y$.

Finally, suppose that
\[
(x_1,\theta)\ge(x_2,\theta),
\qquad
(x_2,\theta)\ge(x_3,\theta).
\]
For all sufficiently large $t$, both successive orbit
differences belong to $X^+\cup\{0\}$.
Since this set is closed under addition, we obtain
\[
\begin{aligned}
&\tilde\phi(t,x_1,\theta)-\tilde\phi(t,x_3,\theta)\\
&\quad=
\bigl(\tilde\phi(t,x_1,\theta)-\tilde\phi(t,x_2,\theta)\bigr)
+\bigl(\tilde\phi(t,x_2,\theta)-\tilde\phi(t,x_3,\theta)\bigr)
\in X^+\cup\{0\}.
\end{aligned}
\]
Hence $(x_1,\theta)\ge(x_3,\theta)$, proving transitivity.
Therefore $\ge$ defines a total order on
$P^{-1}(\theta)\cap E$.
\end{proof}

We now introduce some fibrewise notation that will be used in the proof of
Theorem~\ref{a-a-one}. For any $\theta\in\Theta$, set
\[
A(\theta):=P^{-1}(\theta)\cap \tilde{\mathcal M}.
\]
If $A(\theta)$ is not a singleton, then for any two distinct points
$(x_1,\theta),(x_2,\theta)\in A(\theta)$, Corollary~\ref{constant-invar-set} yields an index
\[
I(x_1,x_2,\theta)=i
\]
such that
\[
\tilde\phi(t,x_1,\theta)-\tilde\phi(t,x_2,\theta)\in
\operatorname{int}(\mathcal C_i)\setminus \mathcal C_{i-1}
\qquad\text{for all sufficiently large }t.
\]
Define
\[
I_m(A(\theta))
:=\min\Bigl\{ I(x_1,x_2,\theta):(x_1,\theta),(x_2,\theta)\in A(\theta),\ x_1\neq x_2\Bigr\}.
\]
For any $(x_0,\theta)\in A(\theta)$, set
\begin{equation}\label{neighborhood}
A(x_0,\theta)
:=\Bigl\{(x,\theta)\in A(\theta):x\neq x_0,\ I(x,x_0,\theta)=I_m(A(\theta))\Bigr\}.
\end{equation}

The next two lemmas describe the structure of these sets.

\begin{lemma}\label{non-empty}
Let $\Theta'\subset\Theta$ be the residual subset in Lemma~\ref{epimorphism-thm}. Assume
that $\theta\in\Theta'$ and that $A(\theta)$ contains at least two points. Then
$A(x_0,\theta)\neq\varnothing$ for every $(x_0,\theta)\in A(\theta)$.
\end{lemma}

\begin{proof}
Choose $(x_1^*,\theta),(x_2^*,\theta)\in A(\theta)$ such that
\[
I(x_1^*,x_2^*,\theta)=I_m(A(\theta)).
\]
Let $T>0$ be such that
\[
\tilde\phi(t,x_1^*,\theta)-\tilde\phi(t,x_2^*,\theta)
\in \operatorname{int}\bigl(\mathcal{C}_{I(x_1^*,x_2^*,\theta)}\bigr)
\setminus \mathcal{C}_{I(x_1^*,x_2^*,\theta)-1}
\quad\text{for all }t\ge T.
\]
By Lemma~\ref{epimorphism-thm} and the minimality of $\tilde{\mathcal M}$, for any
$(x_0,\theta)\in A(\theta)$ there exist a sequence $t_n\to+\infty$ and points
$(x_n,\theta)\in A(\theta)$ such that
\[
\theta\cdot t_n\to\theta,\qquad
\tilde\Lambda_{t_n}(x_0,\theta)\to(x_1^*,\theta),\qquad
\tilde\Lambda_{t_n}(x_n,\theta)\to(x_2^*,\theta).
\]
Hence there exists $n_0\in\NN$ such that, for all $n\ge n_0$,
\[
\tilde\phi(t_n+T,x_0,\theta)-\tilde\phi(t_n+T,x_n,\theta)
\in \operatorname{int}\bigl(\mathcal{C}_{I(x_1^*,x_2^*,\theta)}\bigr)
\setminus \mathcal{C}_{I(x_1^*,x_2^*,\theta)-1}.
\]
By minimality of the index $I_m(A(\theta))$, this implies
\[
I(x_0,x_n,\theta)=I_m(A(\theta))\qquad\text{for all }n\ge n_0.
\]
Therefore $(x_n,\theta)\in A(x_0,\theta)$ for all sufficiently large $n$, and
$A(x_0,\theta)\neq\varnothing$.
\end{proof}

\begin{lemma}\label{maximum-value}
For any $\theta\in\Theta'$ such that $A(\theta)$ contains at least two points, there exists
$(x_0,\theta)\in A(\theta)$ with the property that
\[
(x_0,\theta)\ge (x,\theta)\qquad\text{for all }(x,\theta)\in A(x_0,\theta).
\]
\end{lemma}

\begin{proof}
Assume the contrary. Then there exists $\theta^*\in\Theta'$ such that for every
$(x^*,\theta^*)\in A(\theta^*)$ there is a point $(x_0^*,\theta^*)\in A(x^*,\theta^*)$ with
\[
(x_0^*,\theta^*)>(x^*,\theta^*).
\]
Since $I(x_0^*,x^*,\theta^*)=I_m(A(\theta^*))$, Corollary~\ref{constant-invar-set} yields $T>0$
such that
\[
\tilde\phi(t,x_0^*,\theta^*)-\tilde\phi(t,x^*,\theta^*)
\in \Bigl(\operatorname{int}(\mathcal{C}_{I_m(A(\theta^*))})
\setminus \mathcal{C}_{I_m(A(\theta^*))-1}\Bigr)\cap X^+
\quad\text{for all }t\ge T.
\]
By continuity at the two times $T$ and $T+2T_0$, there exists an open neighborhood
$B(x^*,\theta^*)$ of $(x^*,\theta^*)$ in $A(\theta^*)$ such that
\[
\tilde\phi(T,x_0^*,\theta^*)-\tilde\phi(T,x,\theta^*)
\in \Bigl(\operatorname{int}(\mathcal{C}_{I_m(A(\theta^*))})
\setminus \mathcal{C}_{I_m(A(\theta^*))-1}\Bigr)\cap X^+
\]
for all $(x,\theta^*)\in B(x^*,\theta^*)$, and the corresponding
difference at time $T+2T_0$ belongs to $X^+$ as well.
No difference between two distinct points of $A(\theta^*)$ can enter
$\mathcal C_{I_m(A(\theta^*))-1}$: by Proposition~\ref{pertur-pro}(1),
such entry would force a smaller eventual index, contradicting its
minimality (for index one, use injectivity on the minimal set).
Proposition~\ref{pertur-pro}(1) therefore places each of the above
differences in the same cone layer for every $t\ge T+T_0$.
Part (3) of that proposition then excludes $\mathcal H$ for every
$t\ge T+2T_0$. Continuity and the sign at $T+2T_0$ imply that these
differences remain in $X^+$ thereafter.
Thus
\[
I(x_0^*,x,\theta^*)=I_m(A(\theta^*))\qquad\text{for all }(x,\theta^*)\in B(x^*,\theta^*),
\]
and consequently
\[
(x_0^*,\theta^*)>(x,\theta^*)\qquad\text{for all }(x,\theta^*)\in B(x^*,\theta^*).
\]
So $(x_0^*,\theta^*)$ is an upper bound of $B(x^*,\theta^*)$ in $A(\theta^*)$.

Since $\{B(x^*,\theta^*):(x^*,\theta^*)\in A(\theta^*)\}$ is an open cover of the compact set
$A(\theta^*)$, we may choose finitely many points $(x_i^*,\theta^*)\in A(\theta^*)$,
$i=1,\dots,n_0$, such that
\[
\bigcup_{i=1}^{n_0}B(x_i^*,\theta^*)=A(\theta^*).
\]
Let $(x_0,\theta^*)$ be a maximum element of the finite set
$\{(x_{0i}^*,\theta^*)\}_{i=1}^{n_0}$, where each $(x_{0i}^*,\theta^*)$ is the upper bound
corresponding to $B(x_i^*,\theta^*)$. Then
\[
(x_0,\theta^*)\ge (x_{0i}^*,\theta^*)\qquad\text{for }i=1,\dots,n_0.
\]
It follows that $(x_0,\theta^*)$ is an upper bound of all of $A(\theta^*)$. By
Lemma~\ref{non-empty}, the set $A(x_0,\theta^*)$ is nonempty. But then for any
$(x,\theta^*)\in A(x_0,\theta^*)$ we have
\[
(x_0,\theta^*)>(x,\theta^*),
\]
contradicting the assumption that no such $(x_0,\theta^*)$ exists. The lemma follows.
\end{proof}

\begin{proof}[Proof of Theorem~\ref{a-a-one}]
Suppose, by contradiction, that $A(\theta)$ is not a singleton for some $\theta\in\Theta'$.
By Lemma~\ref{maximum-value}, there exists $(x_0,\theta)\in A(\theta)$ such that
\[
(x_0,\theta)\ge(x,\theta)\qquad\text{for all }(x,\theta)\in A(x_0,\theta).
\]

Choose distinct points $(x_1^*,\theta),(x_2^*,\theta)\in A(\theta)$ such that
\[
I(x_1^*,x_2^*,\theta)=I_m(A(\theta)).
\]
Without loss of generality, assume
\[
(x_1^*,\theta)>(x_2^*,\theta),
\]
and choose $T>0$ such that
\[
\tilde\phi(t,x_1^*,\theta)-\tilde\phi(t,x_2^*,\theta)
\in \Bigl(\operatorname{int}(\mathcal{C}_{I_m(A(\theta))})
\setminus \mathcal{C}_{I_m(A(\theta))-1}\Bigr)\cap X^+
\quad\text{for all }t\ge T.
\]

By the minimality of $\tilde{\mathcal M}$ and Lemma~\ref{epimorphism-thm}, there exist a
sequence $t_n\to+\infty$ and points $(x_n,\theta)\in A(\theta)$ such that
\[
\tilde\Lambda_{t_n}(x_0,\theta)\to(x_2^*,\theta),\qquad
\tilde\Lambda_{t_n}(x_n,\theta)\to(x_1^*,\theta).
\]
Arguing exactly as in the proofs of Lemmas~\ref{non-empty} and \ref{maximum-value}, we obtain
$n_0\in\NN$ such that for all $n\ge n_0$,
\[
\tilde\phi(t_n+T,x_n,\theta)-\tilde\phi(t_n+T,x_0,\theta)
\in \Bigl(\operatorname{int}(\mathcal{C}_{I_m(A(\theta))})
\setminus \mathcal{C}_{I_m(A(\theta))-1}\Bigr)\cap X^+,
\]
and
\[
I(x_n,x_0,\theta)=I_m(A(\theta)).
\]
By continuity, the difference at time $t_n+T+2T_0$ also belongs to $X^+$
for all sufficiently large $n$. The minimal-index and transversality
argument in Lemma~\ref{maximum-value} shows that it remains in $X^+$
thereafter.
Hence $(x_n,\theta)\in A(x_0,\theta)$ and $(x_n,\theta)>(x_0,\theta)$ for all $n\ge n_0$,
contradicting the maximality property of $(x_0,\theta)$. Therefore $A(\theta)$ must be a
singleton for every $\theta\in\Theta'$. This shows that $\tilde{\mathcal M}$ is an almost
$1$-cover of the base flow, and the proof is complete.
\end{proof}

\subsubsection{Proof of Theorem~\ref{perturbation-thm}}

We next study the codimension-one structure of compact $\omega$-limit sets. The first lemma
shows that any two minimal sets are separated by a single cone layer.

\begin{lemma}\label{constant-on-twosets}
Let $\tilde{\mathcal M}_1,\tilde{\mathcal M}_2\subset\mathcal{V}_0$ be two
distinct minimal sets of $\tilde\Lambda_t$. Then there exists an integer
$j_*\in\{1,\dots,N_0\}$ such that, for any
$\theta\in\Theta$ and any $(x_i,\theta)\in\tilde{\mathcal M}_i\cap P^{-1}(\theta)$ ($i=1,2$),
\begin{equation}\label{constant-on-twosets1}
\tilde\phi(t,x_1,\theta)-\tilde\phi(t,x_2,\theta)
\in \operatorname{int}(\mathcal{C}_{j_*})\setminus\mathcal{C}_{j_*-1}
\quad\text{for all }t\in\RR.
\end{equation}
\end{lemma}

\begin{proof}
Fix $\theta\in\Theta$ and $(x_i,\theta)\in\tilde{\mathcal M}_i\cap P^{-1}(\theta)$ ($i=1,2$).
Since distinct minimal sets are disjoint compact sets, their distance is
positive. Thus the two complete orbits remain uniformly separated in both
time directions, and the backward nonconvergence hypothesis in
Corollary~\ref{constant-invar-set} is satisfied.
We first claim that there exists some $j_*\in\{1,\dots,N_0\}$ such that
\begin{equation}\label{constant-onfibre}
\tilde\phi(t,x_1,\theta)-\tilde\phi(t,x_2,\theta)
\in \operatorname{int}(\mathcal{C}_{j_*})\setminus\mathcal{C}_{j_*-1}
\quad\text{for all }t\in\RR.
\end{equation}

By Corollary~\ref{constant-invar-set}, there exist $T>0$ and integers
$1\le j_+\le j_-\le N_0$ such that
\begin{equation}\label{positive-constant}
\tilde\phi(t,x_1,\theta)-\tilde\phi(t,x_2,\theta)
\in \operatorname{int}(\mathcal{C}_{j_+})\setminus\mathcal{C}_{j_+-1}
\quad\text{for all }t\ge T,
\end{equation}
and
\begin{equation}\label{negative-constant}
\tilde\phi(t,x_1,\theta)-\tilde\phi(t,x_2,\theta)
\in \operatorname{int}(\mathcal{C}_{j_-})\setminus\mathcal{C}_{j_- -1}
\quad\text{for all }t\le -T.
\end{equation}
It suffices to prove $j_+=j_-$.

Choose $\theta^*$ in the intersection of the residual singleton-fibre
sets of $\tilde{\mathcal M}_1$ and $\tilde{\mathcal M}_2$, and denote
the unique points in those fibres by
\[
\tilde{\mathcal M}_i\cap P^{-1}(\theta^*)=\{(x_i^*,\theta^*)\},\qquad i=1,2,
\]
respectively.  Minimality of the base flow supplies sequences
$t_n^+\to+\infty$ and $t_n^-\to-\infty$ such that
$\theta\cdot t_n^\pm\to\theta^*$.  Compactness, followed by passage to
subsequences, and uniqueness of the two fibres give
\[
\tilde\Lambda_{t_n^+}(x_i,\theta)\to(x_i^*,\theta^*),
\qquad
\tilde\Lambda_{t_n^-}(x_i,\theta)\to(x_i^*,\theta^*),
\qquad i=1,2.
\]

Lemma~\ref{sequence-limit}, applied to both sequences, puts the same
limiting pair in the layers indexed by $j_+$ and $j_-$ for every
$t\in\RR$. Since distinct cone layers are disjoint, $j_+=j_-$.

Write their common value as $j_*$. Fix $s\in\RR$ and choose
$a\le -T$ with $s-a\ge T_0$. By \eqref{negative-constant} and
Proposition~\ref{pertur-pro}(1), the difference at time $s$ belongs to
$\operatorname{int}\mathcal C_{j_*}$.
If $j_*>1$ and this difference belonged to $\mathcal C_{j_*-1}$,
the same proposition would place all sufficiently late differences in
$\operatorname{int}\mathcal C_{j_*-1}$, contradicting
\eqref{positive-constant}. For $j_*=1$, exclusion of
$\mathcal C_0=\{0\}$ follows from injectivity of the flow.
Since $s$ is arbitrary, \eqref{constant-onfibre} follows.

Finally, fix the singleton-fibre base point $\theta^*$ chosen above.
Every pair drawn from the two minimal sets in any common fibre has an
all-time cone-layer index by \eqref{constant-onfibre}. Minimality of the
base supplies positive times at which its base coordinate tends to
$\theta^*$. Compactness and uniqueness of the two fibres give convergence
to $(x_1^*,\theta^*)$ and $(x_2^*,\theta^*)$. By
Lemma~\ref{sequence-limit}, the original pair and this fixed limiting
pair have the same index. Thus the index is independent of the base point
and of the points chosen in each fibre.
\end{proof}

Recall that
\[
X=\mathcal H\oplus \operatorname{span}\{u^+\}.
\]
The next lemma shows that different minimal sets are separated along the scalar coordinate in the
$u^+$-direction.

\begin{lemma}\label{minimal-separated}
Any two distinct minimal invariant sets
$\tilde{\mathcal M}_1,\tilde{\mathcal M}_2\subset\mathcal{V}_0$
of $\tilde\Lambda_t$ are separated in the following sense:

\begin{itemize}
\item[{\rm(i)}] For all $\theta\in\Theta$,
\[
[m_1(\theta),M_1(\theta)]\cap [m_2(\theta),M_2(\theta)]=\varnothing,
\]
where
\[
\begin{split}
m_i(\theta)
&:=\min\{s_x:x=x_u+s_xu^+,\ x_u\in\mathcal H,\ (x,\theta)\in\tilde{\mathcal M}_i\cap P^{-1}(\theta)\},\\
M_i(\theta)
&:=\max\{s_x:x=x_u+s_xu^+,\ x_u\in\mathcal H,\ (x,\theta)\in\tilde{\mathcal M}_i\cap P^{-1}(\theta)\},
\end{split}
\]
for $i=1,2$.

\item[{\rm(ii)}] If $m_2(\theta^*)-M_1(\theta^*)>0$ for some $\theta^*\in\Theta$, then there
exists $\delta>0$ such that
\[
m_2(\theta)-M_1(\theta)\ge \delta\qquad\text{for all }\theta\in\Theta.
\]
\end{itemize}
\end{lemma}

\begin{proof}
Let $q:X\to\RR$ be the continuous functional determined by
$x=x_u+q(x)u^+$, $x_u\in\mathcal H$.  Choose
$\bar\theta$ in the intersection of the two residual sets supplied by
Theorem~\ref{a-a-one}.  Then
\[
 \tilde{\mathcal M}_i\cap P^{-1}(\bar\theta)
 =\{(\bar x_i,\bar\theta)\},\qquad i=1,2.
\]
Lemma~\ref{constant-on-twosets} and Proposition~\ref{pertur-pro}(3), applied
to time translates with endpoints in the common cone layer, imply that every
difference between the two minimal sets avoids $\mathcal H$ at every time.
In particular $q(\bar x_2)-q(\bar x_1)\ne0$.  Relabeling the two minimal
sets is not needed; retain their given labels and set
\[
\varsigma=\operatorname{sgn}\bigl(q(\bar x_2)-q(\bar x_1)\bigr)
\in\{-1,1\}.
\]

Fix $\theta\in\Theta$ and arbitrary
$(x_i,\theta)\in\tilde{\mathcal M}_i\cap P^{-1}(\theta)$, $i=1,2$.
Choose $t_n\to+\infty$ such that
$\theta\cdot t_n\to\bar\theta$.  Compactness and invariance allow us to pass
to a subsequence for which
\[
 \tilde\Lambda_{t_n}(x_i,\theta)\longrightarrow
 (\bar x_i,\bar\theta),\qquad i=1,2;
\]
the limits are forced by the singleton fibres at $\bar\theta$.  The scalar
function
\[
 t\longmapsto q\bigl(\tilde\phi(t,x_2,\theta)
                  -\tilde\phi(t,x_1,\theta)\bigr)
\]
never vanishes, again by the preceding avoidance of $\mathcal H$, and hence
has constant sign. Its limit along $t_n$ has sign $\varsigma$.
Therefore every cross-pair in every fibre has the same orientation:
\[
\varsigma\bigl(q(x_2)-q(x_1)\bigr)>0.
\]
If $\varsigma=1$, taking extrema gives
$M_1(\theta)<m_2(\theta)$; if $\varsigma=-1$, it gives
$M_2(\theta)<m_1(\theta)$. This proves part {\rm(i)} without changing the
labels of the two minimal sets.

Now assume the hypothesis of part {\rm(ii)}. It forces
$\varsigma=1$.  The fibre product
\[
 \mathcal P_{12}:=
 \{(x_1,x_2,\theta):(x_i,\theta)\in\tilde{\mathcal M}_i,\ i=1,2\}
\]
is compact.  The continuous function
$(x_1,x_2,\theta)\mapsto q(x_2)-q(x_1)$ is strictly positive on
$\mathcal P_{12}$, and therefore has a positive minimum $\delta$.
Consequently
\[
 m_2(\theta)-M_1(\theta)\ge\delta
 \qquad\text{for every }\theta\in\Theta,
\]
which proves part {\rm(ii)}.
\end{proof}
\begin{proof}[Proof of Theorem~\ref{perturbation-thm}]
We divide the proof into two steps.

\medskip
\noindent\textbf{Step 1. $\tilde\Omega$ contains at most two minimal sets.}
Suppose, to the contrary, that
$\tilde\Omega=\tilde\omega(x_0,\theta_0)$ contains three distinct
minimal sets $\tilde{\mathcal M}_i$, $i=1,2,3$.
For each $\theta\in\Theta$, define
\[
A_i(\theta)
:=\{s_x:(x,\theta)\in\tilde{\mathcal M}_i\},
\qquad
m_i(\theta):=\min A_i(\theta),
\qquad
M_i(\theta):=\max A_i(\theta),
\]
where $x=x_u+s_xu^+$ with $x_u\in\mathcal H$.
These extrema exist because each fibre of
$\tilde{\mathcal M}_i$ is nonempty and compact.
By Lemma~\ref{minimal-separated}(ii), after relabeling the minimal
sets if necessary, there exists $\delta>0$ such that
\begin{equation}\label{inequality}
M_1(\theta)+\delta\le m_2(\theta)\le M_2(\theta),
\qquad
M_2(\theta)+\delta\le m_3(\theta),
\qquad \theta\in\Theta.
\end{equation}

Precompactness and $\tilde\Omega\subset\mathcal V_0$ imply that
the positive orbit is eventually contained in $\mathcal V_0$: otherwise,
a sequence of points outside this open set at times tending to infinity
would have a limit in $\tilde\Omega$. After replacing
$(x_0,\theta_0)$ by a sufficiently large positive-time translate, we may assume that its entire positive
orbit is contained in $\mathcal V_0$; this does not change
$\tilde\Omega$.
Choose
$(x_2,\theta_0)\in\tilde{\mathcal M}_2\cap P^{-1}(\theta_0)$
and apply Theorem~\ref{prop-orbit} to the pair
$(x_0,\theta_0)$ and $(x_2,\theta_0)$.

If the two trajectories coincide at a finite time, then their
$\omega$-limit sets are equal. The same conclusion holds if
\begin{equation}\label{exponential-decay3}
\|\tilde\phi(t,x_2,\theta_0)
-\tilde\phi(t,x_0,\theta_0)\|
\longrightarrow0
\qquad(t\to+\infty).
\end{equation}
Since $\tilde{\mathcal M}_2$ is minimal,
$\tilde\omega(x_2,\theta_0)=\tilde{\mathcal M}_2$.
Thus either possibility would imply
$\tilde\Omega=\tilde{\mathcal M}_2$, contradicting the presence
of three distinct minimal sets.

The remaining alternative in Theorem~\ref{prop-orbit}, together
with Lemma~\ref{liein-halfspace}, therefore yields
$j\in\{1,\ldots,N_0\}$ and $T>0$ such that
\begin{equation}\label{stayin-cone2}
\tilde\phi(t,x_2,\theta_0)-\tilde\phi(t,x_0,\theta_0)
\in
\bigl(\operatorname{int}(\mathcal C_j)
\setminus\mathcal C_{j-1}\bigr)\cap X^\pm,
\qquad t\ge T,
\end{equation}
where the choice of sign is fixed for all $t\ge T$.

Suppose first that the $X^+$ alternative holds. Then
\[
s_{\tilde\phi(t,x_2,\theta_0)}
>
s_{\tilde\phi(t,x_0,\theta_0)},
\qquad t\ge T.
\]
Choose a point
$(x_3^*,\theta^*)\in\tilde{\mathcal M}_3$.
Since $\tilde{\mathcal M}_3\subset\tilde\Omega$, there exists
a sequence $t_n\to+\infty$ such that
\[
\tilde\Lambda_{t_n}(x_0,\theta_0)
\longrightarrow(x_3^*,\theta^*).
\]
By compactness of $\tilde{\mathcal M}_2$, we may pass to a
subsequence such that
\[
\tilde\Lambda_{t_n}(x_2,\theta_0)
\longrightarrow(x_2^*,\theta^*)
\in\tilde{\mathcal M}_2.
\]
The limiting base point is the same because the two trajectories
have the common base coordinate $\theta_0\cdot t_n$.
Continuity of the scalar projection gives
\[
s_{x_3^*}\le s_{x_2^*}.
\]
On the other hand, \eqref{inequality} implies
\[
s_{x_2^*}\le M_2(\theta^*)
<m_3(\theta^*)\le s_{x_3^*},
\]
a contradiction.

If the $X^-$ alternative holds, then
\[
s_{\tilde\phi(t,x_2,\theta_0)}
<
s_{\tilde\phi(t,x_0,\theta_0)},
\qquad t\ge T.
\]
Choose $(x_1^*,\theta^*)\in\tilde{\mathcal M}_1$.
By the same argument, there exist a sequence $t_n\to+\infty$
and a point $(x_2^*,\theta^*)\in\tilde{\mathcal M}_2$ such that
\[
\tilde\Lambda_{t_n}(x_0,\theta_0)
\longrightarrow(x_1^*,\theta^*),
\qquad
\tilde\Lambda_{t_n}(x_2,\theta_0)
\longrightarrow(x_2^*,\theta^*).
\]
Passing to the limit yields $s_{x_2^*}\le s_{x_1^*}$, whereas
\eqref{inequality} gives
\[
s_{x_1^*}\le M_1(\theta^*)
<m_2(\theta^*)\le s_{x_2^*}.
\]
This is again a contradiction. Hence $\tilde\Omega$ contains
at most two minimal sets.

\medskip
\noindent\textbf{Step 2. Classification of the remaining part.}
The $\omega$-limit set $\tilde\Omega$ is nonempty, compact,
connected, and invariant. By
Corollary~\ref{constant-invar-set}, the restriction of
$\tilde\Lambda_t$ to $\tilde\Omega$ extends uniquely to a flow.
Consequently, the $\omega$- and $\alpha$-limit sets of every
point of $\tilde\Omega$ are nonempty compact invariant subsets
of $\tilde\Omega$, and each contains a minimal set.

Any compact invariant subset of $\tilde\Omega$ that intersects a
minimal set must contain it. Indeed, their intersection is a
nonempty compact invariant subset of the minimal set and therefore
coincides with it.

Suppose first that $\tilde\Omega$ contains exactly one minimal
set $\tilde{\mathcal M}_1$.
If $\tilde\Omega=\tilde{\mathcal M}_1$, then alternative
{\rm(1)} holds.
Otherwise, define
\[
\tilde{\mathcal M}_{11}
:=\tilde\Omega\setminus\tilde{\mathcal M}_1
\ne\emptyset.
\]
For every $(x_{11},\theta)\in\tilde{\mathcal M}_{11}$, both
$\tilde\omega(x_{11},\theta)$ and
$\tilde\alpha(x_{11},\theta)$ contain a minimal set.
Uniqueness of the minimal set in $\tilde\Omega$ therefore gives
\[
\tilde{\mathcal M}_1
\subset\tilde\omega(x_{11},\theta),
\qquad
\tilde{\mathcal M}_1
\subset\tilde\alpha(x_{11},\theta).
\]
Thus alternative {\rm(2)} holds.

Suppose now that $\tilde\Omega$ contains exactly two distinct
minimal sets $\tilde{\mathcal M}_1$ and
$\tilde{\mathcal M}_2$, and define
\[
\tilde{\mathcal M}_{12}
:=\tilde\Omega\setminus
\bigl(\tilde{\mathcal M}_1\cup\tilde{\mathcal M}_2\bigr).
\]
Distinct minimal sets are disjoint and compact. Hence
connectedness of $\tilde\Omega$ implies
$\tilde{\mathcal M}_{12}\ne\emptyset$.

For every $(x_{12},\theta)\in\tilde{\mathcal M}_{12}$, each of
$\tilde\omega(x_{12},\theta)$ and
$\tilde\alpha(x_{12},\theta)$ contains at least one of
$\tilde{\mathcal M}_1$ and $\tilde{\mathcal M}_2$.
Moreover, either limit set is disjoint from any minimal set
that it does not contain, by the preceding observation.
Thus each limit set contains $\tilde{\mathcal M}_1$ and is
disjoint from $\tilde{\mathcal M}_2$, contains
$\tilde{\mathcal M}_2$ and is disjoint from
$\tilde{\mathcal M}_1$, or contains both.
These are precisely the alternatives in {\rm(3)}.

The three cases are mutually exclusive and exhaustive.
This completes the proof.
\end{proof}

\section{Applications}\label{applications-sec}

In this section, we apply Theorems~\ref{prop-orbit}, \ref{a-a-one}, and
\ref{perturbation-thm} to three representative classes of systems:
\begin{itemize}
\item[(i)] scalar parabolic equations with almost periodic time-dependence;
\item[(ii)] nonlocal and chemotaxis-type perturbations of parabolic equations;
\item[(iii)] competitive--cooperative tridiagonal systems.
\end{itemize}
Each application follows the same scheme: we first establish dissipativity and
precompactness of the corresponding perturbed semiflow, then verify that the
unperturbed system admits a nested invariant cone structure, and finally invoke
the abstract perturbation results from Section~2 to describe the minimal sets
and $\omega$-limit sets of the perturbed dynamics.

\subsection{Perturbations of parabolic equations}\label{parabolic-sec}

In this subsection, we consider the almost-periodic scalar parabolic equation
\eqref{parabolic-eq-01} and two of its perturbations. We first verify that the
unperturbed equation fits into the abstract NIC framework, and then treat
a nonlocal perturbation and a chemotaxis-type system.

\subsubsection{The unperturbed almost-periodic parabolic equation}

Assume that
\[
f(t,x,u,p)\in C^2(\R\times[0,1]\times\R\times\R,\R)
\]
and that $f$ and all its partial derivatives up to order $2$ are uniformly
almost periodic in $t$. Let $\Theta=H(f)$ be the hull of $f$. By
Remark~\ref{a-p-to-minial}(3), there exists a unique continuous function
$F$ on $\Theta\times[0,1]\times\R\times\R$ that is twice continuously
differentiable in $(x,u,p)$, with derivatives depending continuously on
$\theta$, and such that, for any $\theta\in\Theta$ and $t\in\R$,
\[
F(\theta\cdot t,x,u,p)=\theta(t,x,u,p),
\]
where $\theta\cdot t$ denotes the translation of $\theta$ by $t$
(cf.~\cite[Theorem~I.3.2]{Shen1998}). For simplicity, we continue to write
$F(\theta\cdot t,x,u,p)$ as $f(\theta\cdot t,x,u,p)$.

For each $\theta\in\Theta$, consider the induced equation
\begin{equation}\label{induced-para-eq}
\begin{cases}
u_t=u_{xx}+f(\theta\cdot t,x,u,u_x),\quad t>0,\ 0<x<1,\\
u_x(t,0)=u_x(t,1)=0,\quad t>0.
\end{cases}
\end{equation}
Let $A=I-\partial_{xx}$ be the Neumann realization on $L^2(0,1)$,
with domain $D(A)=\{u\in H^2(0,1):u_x(0)=u_x(1)=0\}$, and set
$X=X^\alpha=D(A^\alpha)$, where $3/4<\alpha<1$.
Then $X\hookrightarrow C^1[0,1]$, and every $u\in X$ satisfies the
Neumann boundary conditions. For any $u\in X$, \eqref{induced-para-eq}
defines a local solution $\phi(t,\cdot;u,\theta)$ in $X$ with
$\phi(0,\cdot;u,\theta)=u(\cdot)$, depending continuously on $(u,\theta)\in
X\times\Theta$. Consequently, \eqref{induced-para-eq} generates a local
skew-product semiflow
\begin{equation}\label{para-skew-produ}
\Lambda_t(u,\theta)=(\phi(t,\cdot;u,\theta),\theta\cdot t),
\qquad 0\le t<T_{\max}(u,\theta),
\end{equation}
on $X\times\Theta$, where $T_{\max}(u,\theta)$ is the maximal
existence time.

We impose the following dissipative condition {\bf(D)} (cf.\
\cite{CCH,Chen-P,JR2,Pol}): there exist constants $\delta,\zeta>0$, an exponent $\eta\in(0,2)$,
and a continuous function $C:[0,\infty)\to[0,\infty)$ such that
\begin{equation}\label{dissipative-condition}
\begin{split}
&\forall\,l>0,\ p\in\R,\quad
\sup_{(t,x,u)\in\R\times[0,1]\times[-l,l]}
|f(t,x,u,p)| \le C(l)\bigl(1+|p|^{2-\eta}\bigr),\\
&\forall\,|u|\ge\delta,\ (t,x)\in\R\times[0,1],\quad
u\,f(t,x,u,0)\le -\zeta.
\end{split}
\end{equation}
Then for every $\theta\in\Theta$,
\begin{equation}\label{dissipative-condition-2}
\begin{split}
&\forall\,l>0,\ p\in\R,\quad
\sup_{(t,x,u)\in\R\times[0,1]\times[-l,l]}
|f(\theta\cdot t,x,u,p)| \le C(l)\bigl(1+|p|^{2-\eta}\bigr),\\
&\forall\,|u|\ge\delta,\ (t,x)\in\R\times[0,1],\quad
u\,f(\theta\cdot t,x,u,0)\le -\zeta.
\end{split}
\end{equation}

\medskip
\noindent{\it A common regularity estimate.}
We record the regularity argument used for the unperturbed equation
and both perturbations below. Let $u\in C([0,T_{\max}),X)$ solve
$u_t=u_{xx}+F(t,x)$ with homogeneous Neumann boundary conditions.
Suppose that, throughout its existence interval,
\begin{equation}\label{parabolic-common-growth}
\|u(t)\|_\infty\le M,
\qquad
\|F(t,\cdot)\|_\infty
\le C_M\bigl(1+\|u_x(t)\|_\infty^q\bigr),
\qquad 1\le q<2.
\end{equation}
Here $F$ denotes the full nonlinear term evaluated along the solution;
it may include nonlocal terms. The estimates below are uniform over
families for which $M$, $C_M$, and $q$ can be chosen uniformly.

Let $S(t)=e^{t\partial_{xx}}$ denote the Neumann heat semigroup. Its
gradient estimate is
\[
\|\partial_xS(h)v\|_\infty\le Ch^{-1/2}\|v\|_\infty,
\qquad 0<h\le1.
\]
Define
\[
K(h):=\sup_{0<s\le h}s^{1/2}\|u_x(s)\|_\infty,
\qquad 0<h<T_{\max},\quad h\le1.
\]
For $0<s\le h$, the variation-of-constants formula gives
\[
\begin{aligned}
s^{1/2}\|u_x(s)\|_\infty
&\le CM+C_Ms^{1/2}\int_0^s
(s-r)^{-1/2}\bigl(1+\|u_x(r)\|_\infty^q\bigr)\,dr\\
&\le CM+C_Ms+C_MK(h)^q s^{1/2}
\int_0^s(s-r)^{-1/2}r^{-q/2}\,dr\\
&\le CM+C_Ms+C_MK(h)^q s^{1-q/2}.
\end{aligned}
\]
The last integral is finite because $q<2$. Here $C_M$ may be
enlarged, independently of the particular solution. Consequently,
\[
K(h)\le CM+C_Mh+C_Mh^{1-q/2}K(h)^q.
\]
Since $X\hookrightarrow C^1([0,1])$, the function $K$ is continuous
and $K(h)\to0$ as $h\downarrow0$. Choose $R_1>CM+C_M$ and then
$h_0\in(0,1]$ sufficiently small that
\[
CM+C_Mh_0+C_Mh_0^{1-q/2}R_1^q<R_1.
\]
A first-crossing argument yields $K(h)\le R_1$ whenever
$0<h<T_{\max}$ and $h\le h_0$.
The same argument applies after any time translation within the
existence interval. Thus, for $0<\tau<T_{\max}$, applying it on
$[t-h,t]$ with $h=\min\{\tau,h_0\}$ gives
\begin{equation}\label{nonlocal-gradient-estimate}
\sup_{\tau\le t<T_{\max}}\|u_x(t)\|_\infty
\le R_1h^{-1/2}=:C_{\tau,M}.
\end{equation}
In particular, the constant does not depend on the initial gradient.

Fix $\alpha<\beta<1$ and put $h_1=\min\{\tau/2,1\}$.
For $\tau\le t<T_{\max}$, the formula
\[
u(t)=e^{-h_1A}u(t-h_1)
+\int_{t-h_1}^t e^{-(t-s)A}\bigl(u(s)+F(s,\cdot)\bigr)\,ds
\]
and $t-h_1\ge\tau/2$ imply
\begin{equation}\label{parabolic-common-smoothing}
\|u(t)\|_{X^\beta}
\le Ch_1^{-\beta}M+C_{\tau,M}\int_0^{h_1}r^{-\beta}\,dr
\le C_{\tau,M,\beta}.
\end{equation}
We used \eqref{nonlocal-gradient-estimate} to bound $F$ in $L^2$
and the analytic-semigroup estimate
$\|A^\beta e^{-rA}\|_{L(L^2,L^2)}\le Cr^{-\beta}$.
The integral is finite because $\beta<1$.
If $T_{\max}<\infty$, this bounds the $X$ norm up to
$T_{\max}$, contradicting the continuation criterion for the
locally well-posed equations considered here. Hence a bound of the
form \eqref{parabolic-common-growth} rules out finite-time blow-up.
For global solutions, \eqref{parabolic-common-smoothing} is uniform
for all $t\ge\tau$. Since $X^\beta\hookrightarrow X$ is compact,
a common invariant and absorbing sup-norm bound therefore supplies
a common relatively compact absorbing set in $X\times\Theta$.

\medskip
For the unperturbed equation, write
$U(t)=\|\phi(t,\cdot;u_0,\theta)\|_\infty$ on its maximal
existence interval. The maximum principle and
\eqref{dissipative-condition-2} give
\[
D^+U(t)\le-\frac{\zeta}{U(t)}
\qquad\text{whenever }U(t)\ge\delta.
\]
Thus $U(t)\le\max\{\|u_0\|_\infty,\delta\}$, and {\bf(D)}
verifies \eqref{parabolic-common-growth}, with
$q=\max\{1,2-\eta\}<2$. The preceding argument gives global
existence. Fix $\delta_0>\delta$ and set
\[
\mathcal B:=\{u\in X:\ \|u\|_\infty\le\delta_0\}.
\]
The differential inequality yields absorption into $\mathcal B$
in finite time and strict inward motion at its boundary, uniformly
in $\theta$. Together with \eqref{parabolic-common-smoothing},
this gives the following properties:

\begin{itemize}
\item[{\bf(D1)}] For any $(u_0,\theta)\in X\times\Theta$, there exists
$t_{(u_0,\theta)}>0$ such that
$\phi(t,\cdot;u_0,\theta)\in\mathcal{B}$ for all $t>t_{(u_0,\theta)}$.

\item[{\bf(D2)}] $\mathcal{B}$ is positively invariant:
$\phi(t,\cdot;u_0,\theta)\in\mathcal{B}$ for all $u_0\in\mathcal B$ and
$t\ge0$. In particular,
\[
\phi(t,\cdot;u_0,\theta)\in\mathring{\mathcal{B}}
:=\{u\in X:\ \|u\|_\infty<\delta_0\}
\]
for all $u_0\in\mathcal B$ and $t>0$.

\item[{\bf(D3)}] Equation \eqref{induced-para-eq} (equivalently,
\eqref{para-skew-produ}) admits a compact global attractor $\mathcal{A}$:
\[
\mathcal A
=\bigcap_{T>0}\overline{\bigcup_{t\ge T}
\Lambda_t(\mathcal B\times\Theta)},
\]
\end{itemize}
In {\bf(D3)}, $\mathcal A\subset X\times\Theta$ and the closure is
taken in $X\times\Theta$. The entrance time in {\bf(D1)} is uniform
for $u_0$ in a bounded subset of $X$ and $\theta\in\Theta$.

Given a $C^1$-smooth function $u:[0,1]\to\R$, define its zero number by
\[
Z(u(\cdot)):=\#\{x\in[0,1]: u(x)=0\}.
\]
For each $i\in\mathbb{N}$, define
\begin{equation}\label{parabolic-cone-eq}
C_i:=\overline{\left\{u\in X:\
 (u(x),u_x(x))\ne(0,0)\ \text{for every }x\in[0,1],\quad
 Z(u)\le i-1\right\}}^{\,X}.
\end{equation}
The closure is essential: it includes the zero function and functions with
multiple zeros and makes $C_i$ a closed cone.
Define the linear functional $\pi:X\to\R$ by
\[
\pi(u)=u(0),
\]
and let
\begin{equation}\label{parabolic-hyper-plane-eq}
\mathcal H:=\pi^{-1}(0).
\end{equation}
Then $\mathcal H$ consists of functions vanishing at $x=0$. Since
$\pi\in X^*\setminus\{0\}$, the subspace $\mathcal H$ has codimension one in
$X$.

\begin{lemma}\label{cone-lm}
For each $i\ge1$, $C_i$ is an $i$-cone in $X$.
\end{lemma}

\begin{proof}
Define $W_0=\mathrm{span}\{1\}$, and for any $n\in\mathbb{N}$ let
$W_n=\mathrm{span}\{\cos(n\pi x)\}$. Set
\[
V_i=\bigoplus_{j=0}^{i-1}W_j,\qquad
L_i=\overline{\operatorname{span}}^{\,X}
\{\cos(n\pi x): n\ge i\}.
\]
Then $\dim V_i=i$ and $\operatorname{codim}L_i=i$.
Every element of $V_i$ has the form $p(\cos(\pi x))$, where $p$ is a
polynomial of degree at most $i-1$. Polynomials of this degree with simple
roots and with $p(1)p(-1)\ne0$ are dense in the coefficient topology.
Their compositions with $\cos(\pi x)$ have at most $i-1$ zeros in $(0,1)$,
all simple, and do not vanish at the endpoints. The resulting
approximations converge in the finite-dimensional space $V_i$, hence in
$X$. This proves $V_i\subset C_i$.

For $0\ne w\in L_i$, let $m\ge i$ be the first nonzero Fourier mode of $w$. Then
\[
 e^{m^2\pi^2t}S(t)w\longrightarrow a_m\cos(m\pi x)
 \quad\text{in }X\quad(t\to\infty).
\]
If $w\in C_i$, the zero-number semicontinuity and the strong
variation-diminishing property of $S(t)$ give $S(t)w\in C_i$ for every
$t>0$.  The displayed convergence would then imply that
$\cos(m\pi x)\in C_i$, which is impossible because it has $m\ge i$ simple
zeros. Thus $L_i\cap C_i=\{0\}$.

The set in \eqref{parabolic-cone-eq} is invariant under multiplication by
nonzero real numbers, and its closure is therefore a closed homogeneous cone.
The simple-zero condition is open in $C^1[0,1]$, so $C_i$ has nonempty
interior. Together with $V_i\subset C_i$ and $L_i\cap C_i=\{0\}$, this proves
that $C_i$ is a solid $i$-cone.  This is the standard Sturm cone construction;
see \cite{H.MATANO:1982,Te}.
\end{proof}

Let
\[
\mathcal U_0=X\times\Theta,\qquad
K_0=\mathcal A.
\]
For completeness, we describe the associated difference cocycle.  If
$u_j(t)=\phi(t,\cdot;u_j^0,\theta)$, $j=1,2$, and $w=u_1-u_2$, then
\[
 w_t=w_{xx}+a_1(t,x)w_x+a_0(t,x)w,
 \qquad w_x(t,0)=w_x(t,1)=0,
\]
where, writing $u_s=u_2+s(u_1-u_2)$,
\[
 a_0(t,x)=\int_0^1 f_u(\theta\cdot t,x,u_s,(u_s)_x)\,ds,
 \qquad
 a_1(t,x)=\int_0^1 f_p(\theta\cdot t,x,u_s,(u_s)_x)\,ds.
\]
Let $T(t,u_1^0,u_2^0,\theta)$ be the evolution operator of this linear
equation acting on an arbitrary initial datum.  Uniqueness and time
translation of the coefficient equation give the cocycle identity, while the
choice $w(0)=u_1^0-u_2^0$ gives \eqref{operator-T-eq1}.
Continuous dependence of the linear parabolic equation on its coefficients
and initial data gives the operator-norm continuity for positive times
and the local bounds required in {\bf(H1)}. Parabolic smoothing
gives {\bf(H2)}, and backward uniqueness gives {\bf(H3)}.

Write $z=(u_1^0,u_2^0,\theta)\in\mathcal U_0^e$ and abbreviate
the evolution operator as $T(t,z)$.
To verify {\bf(H4)}, fix $t>0$ and $v\in C_i\setminus\{0\}$.
The Sturm zero-number theorem, applied first to the defining dense subset
of $C_i$ and then by continuity, implies $T(s,z)C_i\subset C_i$ for
$s>0$. Choose $s\in(0,t)$ at which $T(s,z)v$ has only simple interior
zeros and does not vanish at the endpoints. Such regular times exist by
the strict-drop property. This function has at most $i-1$ zeros, and the
same holds in a sufficiently small $X$-neighborhood of it, since
$X\hookrightarrow C^1[0,1]$. By continuity of $T(s,z)$, there is a
neighborhood $\mathcal N$ of $v$, not containing zero, whose image at
time $s$ lies in $C_i$. Cone preservation for the remaining time and
{\bf(H3)} yield $T(t,z)\mathcal N\subset C_i\setminus\{0\}$.
For an arbitrary nonzero $v\in X$, the zero number is finite at positive
times. Approximating $T(t,z)v$ by its images at nearby later regular
times shows that it belongs to some $C_j$. Thus
\[
T(t,z)(X\setminus\{0\})
\subset\bigcup_{j\ge1}(C_j\setminus\{0\}),
\]
which verifies \eqref{regularizing-image}.

For {\bf(H5)}, suppose that $T(t,z)v\in
(C_i\setminus C_{i-1})\cap\mathcal H$ is nonzero and, contrary to the
claim, that $v\in C_i$. The evolved function has a multiple zero at
$x=0$, because both its value and its Neumann derivative vanish there.
At regular times before $t$, its zero number is at most $i-1$.
The strict drop at the boundary therefore gives at most $i-2$ zeros at
regular times immediately after $t$. For $i\ge2$, taking these times
down to $t$ places $T(t,z)v$ in the closed cone $C_{i-1}$, a
contradiction. For $i=1$, the same strict drop is impossible for a
nonzero solution whose zero number is already zero. This proves
{\bf(H5)} with $d=1$.
Hence $\Lambda_t$ satisfies {\bf(H1)--(H5)}.
Therefore, by
Theorems~\ref{prop-orbit}, \ref{a-a-one}, and \ref{perturbation-thm}, we obtain
the following result for perturbations of
\eqref{induced-para-eq}.

\begin{proposition}\label{parobolic-thm}
There exist open bounded neighborhoods $\mathcal V_0$ and $\mathcal W_0$ of
$K_0$, with
$\overline{\mathcal V_0}\subset\mathcal W_0\subset\mathcal U_0$, and a
number $\epsilon_0>0$ with the following property.

For any skew-product semiflow
\[
\tilde\Lambda_t=(\tilde\phi(t,\cdot,\theta),\theta\cdot t)
\]
on $\mathcal U_0$, with
$\widetilde\phi(t,\cdot,\theta)\in C^1(\mathcal W_0(\theta),X)$ for
$t\in[1/2,1]$, such that
\begin{equation}\label{parabolic-c1}
\sup_{\substack{\frac12\le t\le1\\ \theta\in\Theta}}
\bigl\|\widetilde\phi(t,\cdot,\theta)
-\phi(t,\cdot,\theta)\bigr\|_{C^1(\mathcal W_0(\theta),X)}
<\epsilon_0,
\end{equation}
any minimal set of $\tilde\Lambda_t$ contained in $\mathcal V_0$ is an almost
$1$-cover of the base flow $\Theta$. Moreover, every precompact positive orbit
that is eventually contained in $\mathcal V_0$ and whose $\omega$-limit set
is contained in $\mathcal V_0$ satisfies the three alternatives of
Theorem~\ref{perturbation-thm}; in particular, its $\omega$-limit set contains
at most two minimal sets.
\end{proposition}

For semilinear parabolic equations whose nonlinearities are $C^1$-close on a
common bounded set in $X$, standard parameter-dependence estimates for the
semiflow and its fibre derivative give \eqref{parabolic-c1}. We verify
this condition for the two perturbations below.

\subsubsection{A nonlocal perturbation}

We now consider the perturbed system \eqref{parabolic-eq-induced0}.
Assume that $c$ is uniformly almost periodic in $t$ and that
$\nu\in C([0,1])$. For this perturbation, replace $H(f)$ by the
joint hull $H(f,c)$, and continue to denote it by $\Theta$.
The two components of a hull element define $f(\theta,x,u,p)$
with the regularity described above and
$c\in C(\Theta\times[0,1],\R)$.
The unperturbed equation is considered over this same base through
the natural projection $H(f,c)\to H(f)$. Its difference cocycle and
NIC properties carry over to the joint hull. We retain the notation
$\mathcal A$, $\mathcal V_0$, and $\mathcal W_0$ for its attractor
and the neighborhoods in Proposition~\ref{parobolic-thm} over this base.
Thus we consider the family
\begin{equation}\label{parabolic-eq-induced1}
\begin{cases}
u_t=u_{xx}+f(\theta\cdot t,x,u,u_x)
+\epsilon c(\theta\cdot t,x)
 \displaystyle\int_0^1\nu(y)u(t,y)\,dy,
& t>0,\ 0<x<1,\\[4pt]
u_x(t,0)=u_x(t,1)=0,
\end{cases}
\qquad \theta\in\Theta.
\end{equation}
For every $(u_0,\theta)\in X\times\Theta$, this equation admits a
local solution $\phi^N(t,\cdot;u_0,\theta)$ in $X$ on its maximal
existence interval $[0,T_{\max}(u_0,\theta))$, with
$\phi^N(0,\cdot;u_0,\theta)=u_0$. The associated local skew-product
semiflow is
\begin{equation}\label{perturbed-para-skew-produ}
\Lambda_t^N(u_0,\theta)
=
\bigl(\phi^N(t,\cdot;u_0,\theta),\theta\cdot t\bigr),
\qquad 0\le t<T_{\max}(u_0,\theta).
\end{equation}
Set
\[
\|c\|_\infty
:=\max_{(\theta,x)\in\Theta\times[0,1]}|c(\theta,x)|,
\qquad
\|\nu\|_\infty
:=\max_{x\in[0,1]}|\nu(x)|.
\]

To control the nonlocal term, we strengthen the dissipativity part
of {\bf(D)} as follows.

\medskip
\noindent{\bf(N)}
There exist $\xi>0$, $R>0$, and $\eta_0>0$ such that, uniformly in
$(\theta,x)\in\Theta\times[0,1]$,
\begin{equation}\label{condition-eq1}
\begin{cases}
f(\theta,x,u,p)\le-\xi u,
&u\ge R,\quad |p|\le\eta_0,\\
f(\theta,x,u,p)\ge-\xi u,
&u\le-R,\quad |p|\le\eta_0.
\end{cases}
\end{equation}

The following theorem describes the asymptotic dynamics of the
nonlocal perturbation for sufficiently small $|\epsilon|$.

\begin{theorem}\label{non-local-thm}
Assume {\bf(D)} and {\bf(N)}. There exists $\epsilon_1>0$ such that,
whenever $|\epsilon|<\epsilon_1$, every positive orbit of
\eqref{parabolic-eq-induced1} is globally defined and precompact in
$X\times\Theta$. For each $(u_0,\theta)\in X\times\Theta$, its
$\omega$-limit set $\omega^N(u_0,\theta)$ satisfies the three
alternatives in Theorem~\ref{perturbation-thm}. In particular, it
contains at most two minimal sets, each of which is an almost
$1$-cover of the base flow $\Theta$.
\end{theorem}

We first establish global existence and the dissipative estimates
needed in the proof. By \eqref{condition-eq1} and continuity on
bounded sets, there exists $M_0>0$ such that
\begin{equation}\label{condition-eq2}
\begin{cases}
f(\theta,x,u,p)\le-\xi u+M_0,
&u\ge0,\quad |p|\le\eta_0,\\
f(\theta,x,u,p)\ge-\xi u-M_0,
&u\le0,\quad |p|\le\eta_0,
\end{cases}
\end{equation}
uniformly in $(\theta,x)\in\Theta\times[0,1]$.
Assume
\begin{equation}\label{condition-eq0}
|\epsilon|<\frac{\xi}{\|c\|_\infty\|\nu\|_\infty},
\end{equation}
with the right-hand side understood as $+\infty$ if
$\|c\|_\infty\|\nu\|_\infty=0$, and define
\begin{equation}\label{m-star-eq}
M^*:=\frac{M_0}{\xi-|\epsilon|\|c\|_\infty\|\nu\|_\infty}.
\end{equation}

\begin{lemma}\label{c0-estimate-lem}
Given $u_0\in X$, let $\phi^N(t,x;u_0,\theta)$ be the solution of
\eqref{parabolic-eq-induced1} with
$\phi^N(0,x;u_0,\theta)=u_0(x)$ on its maximal existence interval
$[0,T_{\max}(u_0,\theta))$. Then
$T_{\max}(u_0,\theta)=\infty$ and
\[
\limsup_{t\to\infty}
\|\phi^N(t,\cdot;u_0,\theta)\|_\infty\le M^*.
\]
\end{lemma}

\begin{proof}
Write
\[
\phi^N(t,x)=\phi^N(t,x;u_0,\theta),\qquad
U(t)=\|\phi^N(t,\cdot)\|_\infty,\qquad
T_{\max}=T_{\max}(u_0,\theta).
\]
The estimate on the nonlocal term gives
\begin{equation}\label{parabolic-eq-induced-eq1}
\begin{cases}
\phi^N_t\le\phi^N_{xx}
+f(\theta\cdot t,x,\phi^N,\phi^N_x)
+|\epsilon|\|c\|_\infty\|\nu\|_\infty U(t),
&0<x<1,\\
\phi^N_x(t,0)=\phi^N_x(t,1)=0,
\end{cases}
\end{equation}
and
\begin{equation}\label{parabolic-eq-induced-eq2}
\begin{cases}
\phi^N_t\ge\phi^N_{xx}
+f(\theta\cdot t,x,\phi^N,\phi^N_x)
-|\epsilon|\|c\|_\infty\|\nu\|_\infty U(t),
&0<x<1,\\
\phi^N_x(t,0)=\phi^N_x(t,1)=0,
\end{cases}
\end{equation}
for $0<t<T_{\max}$.

Fix $t_0\in[0,T_{\max})$. Let $\phi^+(t,x)$ and $\phi^-(t,x)$
be the local solutions of
\begin{equation}\label{parabolic-eq-induced-eq1-aux}
\begin{cases}
u_t=u_{xx}+f((\theta\cdot t_0)\cdot t,x,u,u_x)
+|\epsilon|\|c\|_\infty\|\nu\|_\infty U(t_0+t),\\
u_x(t,0)=u_x(t,1)=0,\\
u(0,x)=u_0^+(x):=U(t_0),
\end{cases}
\end{equation}
and
\begin{equation}\label{parabolic-eq-induced-eq2-aux}
\begin{cases}
u_t=u_{xx}+f((\theta\cdot t_0)\cdot t,x,u,u_x)
-|\epsilon|\|c\|_\infty\|\nu\|_\infty U(t_0+t),\\
u_x(t,0)=u_x(t,1)=0,\\
u(0,x)=u_0^-(x):=-U(t_0),
\end{cases}
\end{equation}
respectively. On a sufficiently small common existence interval
$[0,\tau_0)$, with $\tau_0<T_{\max}-t_0$, the comparison principle gives
\[
\phi^-(t,x)\le\phi^N(t_0+t,x)\le\phi^+(t,x),
\qquad 0\le t<\tau_0,\quad x\in[0,1].
\]

Let $v^+(t)$ and $v^-(t)$ solve
\begin{equation}\label{ode-eq1}
\begin{cases}
v_t=-\xi v+M_0
+|\epsilon|\|c\|_\infty\|\nu\|_\infty U(t_0+t),\\
v(0)=U(t_0),
\end{cases}
\end{equation}
and
\begin{equation}\label{ode-eq2}
\begin{cases}
v_t=-\xi v-M_0
-|\epsilon|\|c\|_\infty\|\nu\|_\infty U(t_0+t),\\
v(0)=-U(t_0),
\end{cases}
\end{equation}
respectively. Uniqueness and the variation-of-constants formula give
$v^-(t)=-v^+(t)$ and $v^+(t)\ge0$.
By \eqref{condition-eq2}, the spatially constant functions $v^+$
and $v^-$ are, respectively, a supersolution of
\eqref{parabolic-eq-induced-eq1-aux} and a subsolution of
\eqref{parabolic-eq-induced-eq2-aux}. Hence
\[
v^-(t)\le\phi^-(t,x)
\le\phi^N(t_0+t,x)
\le\phi^+(t,x)\le v^+(t)
\]
for $0\le t<\tau_0$ and $x\in[0,1]$. In particular,
\[
U(t_0+t)\le v^+(t),\qquad U(t_0)=v^+(0).
\]
Taking the upper right Dini derivative at $t_0$, we obtain
\[
D^+U(t_0)\le(v^+)'(0)=-aU(t_0)+M_0,
\qquad
a:=\xi-|\epsilon|\|c\|_\infty\|\nu\|_\infty>0.
\]
Since $t_0$ is arbitrary, scalar comparison yields
\begin{equation}\label{sup-norm-comparison}
U(t)\le U(0)e^{-at}+\frac{M_0}{a}(1-e^{-at})
=M^*+\bigl(U(0)-M^*\bigr)e^{-at},
\qquad 0\le t<T_{\max}.
\end{equation}
Consequently,
\[
\sup_{0\le t<T_{\max}}U(t)\le M,
\qquad M:=\max\{\|u_0\|_\infty,M^*\}.
\]

The full nonlinear term satisfies
\[
\biggl\|f(\theta\cdot t,\cdot,u,u_x)
+\epsilon c(\theta\cdot t,\cdot)\int_0^1\nu(y)u(t,y)\,dy\biggr\|_\infty
\le C_M\bigl(1+\|u_x\|_\infty^q\bigr),
\qquad q=\max\{1,2-\eta\}<2,
\]
where $u(t)=\phi^N(t,\cdot)$. Thus
\eqref{parabolic-common-growth} holds on the maximal existence
interval. Estimates \eqref{nonlocal-gradient-estimate} and
\eqref{parabolic-common-smoothing}, followed by the continuation
criterion, give $T_{\max}(u_0,\theta)=\infty$.

Finally, since $a>0$, letting $t\to\infty$ in
\eqref{sup-norm-comparison} yields
\[
\limsup_{t\to\infty}
\|\phi^N(t,\cdot;u_0,\theta)\|_\infty\le M^*.
\]
\end{proof}

We now obtain estimates uniform with respect to the perturbation
parameter. Fix $\bar\epsilon>0$ such that
\[
\bar\epsilon\|c\|_\infty\|\nu\|_\infty<\xi,
\]
choose
\[
M_b>\frac{M_0}{\xi-\bar\epsilon\|c\|_\infty\|\nu\|_\infty},
\]
and set
\[
\mathcal B^N:=\{u\in X:\ \|u\|_\infty\le M_b\}.
\]
We write $\phi^N_\epsilon$ and $\Lambda^N_{\epsilon,t}$ when the
dependence on $\epsilon$ needs to be explicit.

By \eqref{sup-norm-comparison}, solutions starting in $\mathcal B^N$
remain there for all positive times. The constants in
\eqref{parabolic-common-growth} are uniform for
$|\epsilon|\le\bar\epsilon$, $\theta\in\Theta$, and
$\|u\|_\infty\le M_b$. Hence \eqref{parabolic-common-smoothing}
gives, for every $\alpha<\beta<1$ and $\tau>0$,
\[
\sup_{\substack{|\epsilon|\le\bar\epsilon,\ \theta\in\Theta\\
u_0\in\mathcal B^N,\ t\ge\tau}}
\|\phi^N_\epsilon(t,\cdot;u_0,\theta)\|_{X^\beta}
\le C_{\tau,M_b,\beta}.
\]
Together with the comparison estimate, this yields the following
properties for every $|\epsilon|\le\bar\epsilon$:

\begin{itemize}
\item[{\bf(N1)}]
For every $(u_0,\theta)\in X\times\Theta$, there exists
$t_{(u_0,\theta)}>0$ such that
\[
\phi^N_\epsilon(t,\cdot;u_0,\theta)\in\mathcal B^N
\qquad\text{for all }t\ge t_{(u_0,\theta)}.
\]
The entrance time can be chosen uniformly for $\theta\in\Theta$,
$|\epsilon|\le\bar\epsilon$, and $u_0$ in a bounded subset of $X$.

\item[{\bf(N2)}]
The set $\mathcal B^N$ is positively invariant. Moreover,
\[
\phi^N_\epsilon(t,\cdot;u_0,\theta)\in\mathring{\mathcal B}^N
:=\{u\in X:\ \|u\|_\infty<M_b\}
\]
for every $u_0\in\mathcal B^N$, $\theta\in\Theta$, and $t>0$.

\item[{\bf(N3)}]
The skew-product semiflow $\Lambda^N_\epsilon$ admits a compact
global attractor $\mathcal A^N_\epsilon\subset X\times\Theta$, given by
\[
\mathcal A^N_\epsilon
=\bigcap_{T>0}
\overline{\bigcup_{t\ge T}
\Lambda^N_{\epsilon,t}(\mathcal B^N\times\Theta)},
\]
where the closure is taken in $X\times\Theta$.
\end{itemize}

\begin{proof}[Proof of Theorem~\ref{non-local-thm}]
The preceding estimates give global existence and precompactness
of every positive orbit for $|\epsilon|\le\bar\epsilon$.
It remains to verify \eqref{parabolic-c1} and to show that every
positive orbit eventually enters $\mathcal V_0$.

The nonlocal term is a bounded linear map from $X$ to $L^2(0,1)$,
with norm bounded by a constant times
$|\epsilon|\|c\|_\infty\|\nu\|_\infty$, uniformly in $\theta$.
The neighborhoods in Proposition~\ref{parobolic-thm} may be
chosen sufficiently small that the solution segments starting in
$\mathcal W_0$ for $0\le t\le1$ remain in a common bounded subset
of $X$ for sufficiently small $|\epsilon|$. On this set, the nonlinearity and its fibre derivative
satisfy uniform Lipschitz estimates by the regularity assumptions
on $f$ and the embedding $X\hookrightarrow C^1([0,1])$.

The variation-of-constants formulas for the solutions and their
variational equations, together with the integrable smoothing kernel
$t^{-\alpha}$, therefore give
\[
\sup_{\substack{\frac12\le t\le1\\ \theta\in\Theta}}
\|\phi^N_\epsilon(t,\cdot,\theta)-\phi(t,\cdot,\theta)\|_
{C^1(\mathcal W_0(\theta),X)}
\longrightarrow0
\qquad(\epsilon\to0).
\]
This verifies \eqref{parabolic-c1} for sufficiently small
$|\epsilon|$.

To establish eventual entry into $\mathcal V_0$, we prove upper
semicontinuity of the global attractors. The uniform absorbing and
smoothing estimates give a compact set $B\subset X\times\Theta$
containing $\mathcal A^N_\epsilon$ for
$|\epsilon|\le\bar\epsilon$, with $\mathcal A^N_0=\mathcal A$.
The same parameter-dependence argument gives uniform convergence
of the time-$t$ maps on $B$ for each fixed $t>0$; see also
\cite[Chapter~3]{Hen}.

For a compatible product metric $d$ on $X\times\Theta$, write
\[
\operatorname{dist}(A_1,A_2)
:=\sup_{a\in A_1}\inf_{b\in A_2}d(a,b)
\]
for the one-sided Hausdorff semidistance. Invariance of
$\mathcal A^N_\epsilon$ gives, for every $t>0$,
\[
\operatorname{dist}(\mathcal A^N_\epsilon,\mathcal A)
\le
\sup_{z\in B}d(\Lambda^N_{\epsilon,t}z,\Lambda_tz)
+\operatorname{dist}(\Lambda_tB,\mathcal A).
\]
First choose $t$ large, using attraction of $B$ by $\mathcal A$,
and then let $\epsilon\to0$. We obtain
\begin{equation}\label{asym-attractor2}
\operatorname{dist}_{X\times\Theta}
\bigl(\mathcal A^N_\epsilon,\mathcal A\bigr)
\longrightarrow0
\qquad(\epsilon\to0).
\end{equation}
This is the standard upper-semicontinuity argument; see
\cite{HLR,CCDL} for the general semigroup and skew-product settings.

Since $\mathcal V_0$ is an open neighborhood of the compact set
$\mathcal A$, \eqref{asym-attractor2} implies
$\mathcal A^N_\epsilon\subset\mathcal V_0$ for sufficiently small
$|\epsilon|$. For each such $\epsilon$, compactness of
$\mathcal A^N_\epsilon$ gives a neighborhood of
$\mathcal A^N_\epsilon$ contained in $\mathcal V_0$.
Attraction by $\mathcal A^N_\epsilon$ then shows that every positive
orbit eventually enters $\mathcal V_0$ and remains there.
Its $\omega$-limit set is contained in
$\mathcal A^N_\epsilon\subset\mathcal V_0$.

Choosing $\epsilon_1>0$ sufficiently small to satisfy all the
preceding requirements, Proposition~\ref{parobolic-thm} gives the
almost-covering property of minimal sets and the three alternatives
for $\omega$-limit sets.
\end{proof}

\begin{example}[A quasiperiodic nonlocal equation]\label{example-quasiperiodic-nonlocal}
{\rm
Consider
\[
\begin{cases}
u_t=u_{xx}-u^3+(2+\sin t\cos\pi x)u+\cos(\sqrt2\,t)
+\epsilon\cos\pi x\displaystyle\int_0^1u(t,y)\,dy,\\
u_x(t,0)=u_x(t,1)=0.
\end{cases}
\]
Here $f(t,x,u,p)=-u^3+(2+\sin t\cos\pi x)u+\cos(\sqrt2\,t)$,
$c(x)=\cos\pi x$, and $\nu\equiv1$. The joint hull is the
two-dimensional torus with frequencies $1$ and $\sqrt2$.
The nonlinearity is independent of $p$ and has the required admissibility
and gradient-growth properties. Uniformly in $t$ and $x$, for $|u|\ge3$,
\[
uf(t,x,u,p)\le-u^4+3u^2+|u|\le-1,
\qquad
\operatorname{sgn}(u)f(t,x,u,p)
\le-|u|^3+3|u|+1\le-|u|.
\]
Thus {\bf(D)} and {\bf(N)} hold, with $\delta=R=3$ and
$\zeta=\xi=1$. Since $\|c\|_\infty=\|\nu\|_\infty=1$,
the dissipative estimate gives global existence for $|\epsilon|<1$.
For all sufficiently small $|\epsilon|$, Theorem~\ref{non-local-thm}
also gives the almost $1$-covering property and the stated trichotomy
for every $\omega$-limit set.

For every $\epsilon\ne0$, this equation fails to preserve the usual
pointwise order. To see this, choose $x_*\in(0,1)$ with
$\epsilon\cos\pi x_*<0$ and a nonzero nonnegative
$h\in C_c^\infty(0,1)$ that vanishes near $x_*$. Let $u^h$ and $u^0$
be the solutions with initial values $h$ and $0$, at the same base
point. Their initial difference and all its spatial derivatives vanish
near $x_*$, so the equations give
\[
\partial_t(u^h-u^0)(0,x_*)
=\epsilon\cos\pi x_*\int_0^1h(y)\,dy<0.
\]
Consequently, $u^h(t,x_*)<u^0(t,x_*)$ for sufficiently small $t>0$,
although $h\ge0$. This illustrates the application of the theorem to
small perturbations that destroy monotonicity in the pointwise order.
}
\end{example}

\subsubsection{A chemotaxis-type perturbation}

We next consider the perturbation
\begin{equation}\label{parabolic-eq-03-1}
\begin{cases}
u_t=u_{xx}-\epsilon v_xu_x-\epsilon(v-u)u+f(t,x,u,u_x),
&t>0,\ 0<x<1,\\
u_x(t,0)=u_x(t,1)=0,
\end{cases}
\end{equation}
where the elliptic component satisfies
\begin{equation}\label{parabolic-eq-03-2}
\begin{cases}
0=v_{xx}-v+u,&0<x<1,\\
v_x(t,0)=v_x(t,1)=0.
\end{cases}
\end{equation}
Since $v_{xx}=v-u$, the pair $(u,v)$ equivalently solves the
parabolic--elliptic system
\begin{equation}\label{parabolic-eq-03}
\begin{cases}
u_t=u_{xx}-\epsilon(uv_x)_x+f(t,x,u,u_x),&t>0,\ 0<x<1,\\
0=v_{xx}-v+u,&0<x<1,\\
u_x(t,0)=u_x(t,1)=v_x(t,0)=v_x(t,1)=0.
\end{cases}
\end{equation}
This is a one-dimensional Keller--Segel type chemotaxis system;
see \cite{KS1970,KS1971,JL,DN,IS,Nagai1,Nagai2,NT1}.

For this application, we again take $\Theta=H(f)$ and use
$\mathcal A$, $\mathcal V_0$, and $\mathcal W_0$ for the
unperturbed attractor and neighborhoods over this base.
The induced family is
\begin{equation}\label{para-ellip-induced}
\begin{cases}
u_t=u_{xx}-\epsilon(uv_x)_x+f(\theta\cdot t,x,u,u_x),
&t>0,\ 0<x<1,\\
0=v_{xx}-v+u,&0<x<1,\\
u_x(t,0)=u_x(t,1)=v_x(t,0)=v_x(t,1)=0,
\end{cases}
\qquad\theta\in\Theta.
\end{equation}
Its parabolic equation can be written as
\begin{equation}\label{para-ellip-induced-1}
u_t=u_{xx}-\epsilon u_xv_x-\epsilon u(v-u)
+f(\theta\cdot t,x,u,u_x).
\end{equation}

Let $\mathcal R=A^{-1}=(I-\partial_{xx})^{-1}$ denote the Neumann
elliptic resolvent. For each fixed $t$, the elliptic problem has
the unique solution $v(t)=\mathcal Ru(t)$, explicitly given by
\begin{equation}\label{solution-expresion}
\begin{aligned}
v(t,x)
&=c_0(t)(e^x+e^{-x})
-\int_0^x\frac{e^{x-y}-e^{y-x}}{2}u(t,y)\,dy,\\
c_0(t)
&=\int_0^1\frac{e^{2-y}+e^y}{2(e^2-1)}u(t,y)\,dy.
\end{aligned}
\end{equation}
The elliptic maximum principle gives the sharp bound
\begin{equation}\label{estimate-v}
\|\mathcal Ru\|_\infty\le\|u\|_\infty.
\end{equation}
Moreover, $(\mathcal Ru)_{xx}=\mathcal Ru-u$ and
$(\mathcal Ru)_x(0)=0$ imply
\begin{equation}\label{chemotaxis-elliptic-bound}
\|\mathcal Ru\|_{C^2([0,1])}\le C_R\|u\|_\infty
\end{equation}
for a constant $C_R$ independent of $u\in C([0,1])$.

The resulting equation for $u$ is locally well posed in $X$.
We write
\[
u(t,x)=\phi^u(t,x;u_0,\theta),\qquad
v(t,x)=\phi^v(t,x;u_0,\theta)
=\bigl(\mathcal R\phi^u(t,\cdot;u_0,\theta)\bigr)(x).
\]
Thus the elliptic component requires no independent initial datum,
and the local skew-product semiflow is
\[
\Lambda_t^u(u_0,\theta)
=\bigl(\phi^u(t,\cdot;u_0,\theta),\theta\cdot t\bigr),
\qquad 0\le t<T_{\max}(u_0,\theta).
\]

To control the quadratic terms in the perturbation, we impose the
following strengthening of the dissipativity part of {\bf(D)}.

\medskip
\noindent{\bf(P)}
There exist $\xi>0$, $R>0$, and $\eta_1>0$ such that, uniformly in
$(\theta,x)\in\Theta\times[0,1]$,
\begin{equation}\label{condition-eq1-1}
\begin{cases}
f(\theta,x,u,p)\le-\xi u^2,&u\ge R,\quad |p|\le\eta_1,\\
f(\theta,x,u,p)\ge\xi u^2,&u\le-R,\quad |p|\le\eta_1.
\end{cases}
\end{equation}

\begin{theorem}\label{chemotaxis-thm}
Assume {\bf(D)} and {\bf(P)}. There exists $\epsilon_1>0$ such that,
whenever $|\epsilon|<\epsilon_1$, every positive orbit of the
skew-product semiflow associated with \eqref{para-ellip-induced}
is globally defined and precompact in $X\times\Theta$.
For each $(u_0,\theta)\in X\times\Theta$, its $\omega$-limit set
$\omega^u(u_0,\theta)$ satisfies the three alternatives in
Theorem~\ref{perturbation-thm}. In particular, it contains at most
two minimal sets, each of which is an almost $1$-cover of the base
flow $\Theta$.
\end{theorem}

By {\bf(P)}, there exists $M_1>0$ such that
\begin{equation}\label{condition-eq1-2}
\begin{cases}
f(\theta,x,u,p)\le-\xi u^2+M_1,&u\ge0,\quad |p|\le\eta_1,\\
f(\theta,x,u,p)\ge\xi u^2-M_1,&u\le0,\quad |p|\le\eta_1,
\end{cases}
\end{equation}
uniformly in $(\theta,x)\in\Theta\times[0,1]$.
Assume
\begin{equation}\label{condition-eq0-0}
|\epsilon|<\frac{\xi}{2},
\end{equation}
and set
\[
M_1^*:=\sqrt{\frac{M_1}{\xi-2|\epsilon|}}.
\]

\begin{lemma}
For $u_0\in X$, let $(\phi^u,\phi^v)$ be the solution of
\eqref{para-ellip-induced}, with
$\phi^u(0,x;u_0,\theta)=u_0(x)$, on its maximal existence interval
$[0,T_{\max}(u_0,\theta))$. Then
$T_{\max}(u_0,\theta)=\infty$ and
\[
\limsup_{t\to\infty}
\|\phi^u(t,\cdot;u_0,\theta)\|_\infty\le M_1^*.
\]
\end{lemma}

\begin{proof}
Write $u(t,x)=\phi^u(t,x;u_0,\theta)$,
$U(t)=\|u(t,\cdot)\|_\infty$, and
$T_{\max}=T_{\max}(u_0,\theta)$. At a point where a positive
maximum $U(t)$ is attained, $u_x=0$ and $u_{xx}\le0$; this also
holds at an endpoint by the Neumann condition. Since
\[
|\epsilon u(v-u)|\le2|\epsilon|U(t)^2,
\]
estimates \eqref{estimate-v} and \eqref{condition-eq1-2} give
$u_t\le-(\xi-2|\epsilon|)U(t)^2+M_1$ there.
At a negative minimum, the same bound holds for $-u_t$.
The estimate applies to every point where $|u(t,\cdot)|$ attains
its maximum. If $U(t)=0$, then $u=v=0$ and
$|f(\theta\cdot t,x,0,0)|\le M_1$, so the same bound holds.
Consequently,
\[
D^+U(t)\le-a_1U(t)^2+M_1,
\qquad a_1:=\xi-2|\epsilon|>0,\quad 0<t<T_{\max}.
\]
Let $w$ solve
\[
w'=-a_1w^2+M_1,\qquad w(0)=U(0).
\]
This scalar solution is global, nonnegative, and converges to
$M_1^*$. Scalar comparison, with continuity at $t=0$, yields
\begin{equation}\label{chemotaxis-sup-comparison}
U(t)\le w(t)\le\max\{U(0),M_1^*\}=:M,
\qquad 0\le t<T_{\max}.
\end{equation}

By \eqref{chemotaxis-elliptic-bound}, the added drift coefficient
$-\epsilon(\mathcal Ru)_x$ is bounded by $|\epsilon|C_RM$, and
the added zeroth-order term is bounded by $2|\epsilon|M^2$.
Thus the full nonlinear term satisfies
\[
\bigl\|f(\theta\cdot t,\cdot,u,u_x)
-\epsilon\bigl[u_x(\mathcal Ru)_x+u(\mathcal Ru-u)\bigr]\bigr\|_\infty
\le C_M\bigl(1+\|u_x\|_\infty^q\bigr),
\qquad q=\max\{1,2-\eta\}<2.
\]
Together with \eqref{chemotaxis-sup-comparison}, this verifies
\eqref{parabolic-common-growth}. The common regularity argument
therefore gives $T_{\max}=\infty$.
Letting $t\to\infty$ in $U(t)\le w(t)$ now proves the assertion.
\end{proof}

\begin{proof}[Proof of Theorem~\ref{chemotaxis-thm}]
Fix $0<\bar\epsilon<\xi/2$ and choose
$M_b>\sqrt{M_1/(\xi-2\bar\epsilon)}$.
Comparison with $w'=-(\xi-2\bar\epsilon)w^2+M_1$ gives a common
positively invariant sup-norm ball of radius $M_b$ that absorbs
bounded initial sets, uniformly in $\theta$ and
$|\epsilon|\le\bar\epsilon$.
The growth estimate in the preceding proof is uniform on this
ball. Hence \eqref{parabolic-common-smoothing} yields
\[
\sup_{\substack{|\epsilon|\le\bar\epsilon,\ \theta\in\Theta\\
u_0\in X,\ \|u_0\|_\infty\le M_b,\ t\ge\tau}}
\|\phi^u_\epsilon(t,\cdot;u_0,\theta)\|_{X^\beta}
\le C_{\tau,M_b,\beta},
\qquad \alpha<\beta<1,\quad\tau>0.
\]
Here the subscript $\epsilon$ makes the parameter dependence
explicit. These estimates give a common relatively compact
absorbing set, precompactness of positive orbits, and compact
global attractors $\mathcal A^u_\epsilon\subset X\times\Theta$.

The perturbation is the map $G_\epsilon:X\to L^2(0,1)$ given by
\[
G_\epsilon(u)
=-\epsilon\bigl[u_x(\mathcal Ru)_x+u(\mathcal Ru-u)\bigr].
\]
Its derivative is
\[
\begin{aligned}
DG_\epsilon(u)h=-\epsilon\bigl[&h_x(\mathcal Ru)_x
+u_x(\mathcal Rh)_x\\
&+h(\mathcal Ru-u)+u(\mathcal Rh-h)\bigr].
\end{aligned}
\]
By \eqref{chemotaxis-elliptic-bound} and
$X\hookrightarrow C^1([0,1])$, for every bounded subset $B_0\subset X$,
\[
\sup_{u\in B_0}
\bigl(\|G_\epsilon(u)\|_{L^2}
+\|DG_\epsilon(u)\|_{L(X,L^2)}\bigr)
\le C_{B_0}|\epsilon|.
\]
The solution and variational estimates used in the proof of
Theorem~\ref{non-local-thm} therefore give
\[
\sup_{\substack{\frac12\le t\le1\\\theta\in\Theta}}
\|\phi^u_\epsilon(t,\cdot,\theta)-\phi(t,\cdot,\theta)\|_
{C^1(\mathcal W_0(\theta),X)}\longrightarrow0
\qquad(\epsilon\to0),
\]
which verifies \eqref{parabolic-c1} for sufficiently small $|\epsilon|$.

The common compact absorbing bound and finite-time convergence
also allow us to apply the upper-semicontinuity argument leading
to \eqref{asym-attractor2}. Thus
$\operatorname{dist}_{X\times\Theta}(\mathcal A^u_\epsilon,
\mathcal A)\to0$ as $\epsilon\to0$.
For sufficiently small $|\epsilon|$, the compact attractor
$\mathcal A^u_\epsilon$ lies in $\mathcal V_0$ and attracts every
positive orbit into a neighborhood contained in $\mathcal V_0$.
The orbit therefore eventually remains in $\mathcal V_0$, and
its $\omega$-limit set is contained there.
Proposition~\ref{parobolic-thm} completes the proof.
\end{proof}

\subsection{Competitive--cooperative tridiagonal systems}\label{cct}

We next consider an $n$-dimensional tridiagonal system and its
perturbation
\begin{equation}\label{tri-equation-g}
\dot x=f(t,x)+\epsilon g(t,x),
\end{equation}
where $n\ge2$ and $f,g:\R\times\R^n\to\R^n$ are $C^1$ functions.
We assume that $f$, $g$, and their first derivatives with respect
to $x$ are uniformly almost periodic in $t$, uniformly on bounded
subsets of $\R^n$, and that $g$ and $D_xg$ are bounded on
$\R\times\R^n$. The unperturbed system has the tridiagonal form
\begin{equation}\label{tri-equation-f1}
\begin{split}
\dot x_1
&=f_1(t,x_1,x_2),\\
\dot x_i
&=f_i(t,x_{i-1},x_i,x_{i+1}),
\qquad 2\le i\le n-1,\\
\dot x_n
&=f_n(t,x_{n-1},x_n).
\end{split}
\end{equation}

We assume that the two off-diagonal derivatives associated with
each neighboring pair have the same fixed sign and are uniformly
bounded away from zero. More precisely, there exist $\kappa_0>0$
and $\delta_i\in\{-1,1\}$, $1\le i\le n-1$, such that
\[
\delta_i\frac{\partial f_i}{\partial x_{i+1}}(t,x)\ge\kappa_0,
\qquad
\delta_i\frac{\partial f_{i+1}}{\partial x_i}(t,x)\ge\kappa_0
\]
for all $(t,x)\in\R\times\R^n$.
Such systems are called competitive--cooperative tridiagonal
systems. A diagonal change of signs in the state variables
transforms the system into a cooperative one; see
\cite{simith1991}. Applying the same transformation to the
perturbation preserves the stated regularity and boundedness
assumptions. We may therefore assume that
\begin{equation}\label{positive-jacobi-condition}
\frac{\partial f_i}{\partial x_{i+1}}(t,x)\ge\kappa_0,
\qquad
\frac{\partial f_{i+1}}{\partial x_i}(t,x)\ge\kappa_0,
\qquad 1\le i\le n-1,
\end{equation}
for all $(t,x)\in\R\times\R^n$.

The asymptotic behavior of bounded solutions of the unperturbed
system has been studied extensively. In the autonomous case,
bounded trajectories converge to equilibria \cite{Smillie1984};
in the periodic case, they approach periodic solutions
\cite{simith1991}. In the almost periodic case, the $\omega$-limit
set of a bounded orbit contains at most two minimal sets, each of
which is an almost automorphic extension of the base flow
\cite{wangyi2007,Fangchun2013}.
Our purpose is to establish the persistence of this qualitative
structure under small perturbations that need not be tridiagonal.

Let $\Theta=H(f,g)$ be the joint hull. As before, we use $f$ and
$g$ for the corresponding hull extensions. The unperturbed and
perturbed families over this common base are
\begin{equation}\label{tri-equation-induced}
\begin{split}
\dot x_1
&=f_1(\theta\cdot t,x_1,x_2),\\
\dot x_i
&=f_i(\theta\cdot t,x_{i-1},x_i,x_{i+1}),
\qquad 2\le i\le n-1,\\
\dot x_n
&=f_n(\theta\cdot t,x_{n-1},x_n),
\end{split}
\end{equation}
and
\begin{equation}\label{tri-equation-g1}
\begin{split}
\dot x_1
&=f_1(\theta\cdot t,x_1,x_2)
+\epsilon g_1(\theta\cdot t,x),\\
\dot x_i
&=f_i(\theta\cdot t,x_{i-1},x_i,x_{i+1})
+\epsilon g_i(\theta\cdot t,x),
\qquad 2\le i\le n-1,\\
\dot x_n
&=f_n(\theta\cdot t,x_{n-1},x_n)
+\epsilon g_n(\theta\cdot t,x).
\end{split}
\end{equation}
By continuity, \eqref{positive-jacobi-condition} extends to every
element of the hull:
\begin{equation}\label{positive-jacobi-condition2}
\frac{\partial f_i}{\partial x_{i+1}}(\theta\cdot t,x)\ge\kappa_0,
\qquad
\frac{\partial f_{i+1}}{\partial x_i}(\theta\cdot t,x)\ge\kappa_0,
\qquad 1\le i\le n-1,
\end{equation}
for all $(t,x,\theta)\in\R\times\R^n\times\Theta$.
These equations generate local skew-product semiflows
\[
\Lambda_t(x^0,\theta)
=\bigl(\phi(t,x^0,\theta),\theta\cdot t\bigr),
\qquad
\Lambda_t^g(x^0,\theta)
=\bigl(\phi^g(t,x^0,\theta),\theta\cdot t\bigr),
\]
defined on the respective maximal forward existence intervals.

We equip $\R^n$ with the maximum norm
\[
\|x\|:=\max_{1\le i\le n}|x_i|.
\]
For the unperturbed system, consider the following dissipativity
condition.

\medskip
\noindent{\bf(D-T)}
There exist $\delta,C>0$ such that
\[
x_i f_i(t,x)\le-\delta
\]
whenever $|x_i|\ge C$, every existing neighbor of $x_i$ has
absolute value at most $|x_i|$, and $t\in\R$.
Here $f_i(t,x)$ is understood as the component of the full vector
field, with the tridiagonal dependence specified in
\eqref{tri-equation-f1}.

Under {\bf(D-T)}, every solution of the unperturbed equation enters
the box
\[
\mathcal B:=\{x\in\R^n:\|x\|\le C\}
\]
in finite time and remains there; see \cite[Lemma~4.3]{YZ}.
Indeed, for $R(t)=\|\phi(t,x^0,\theta)\|$, the maximum-norm
estimate gives $D^+R(t)\le-\delta/R(t)$ whenever $R(t)\ge C$.
This also rules out finite-time blow-up and gives an entrance time
uniform in $\theta$ and in bounded sets of initial data.
Consequently, the induced skew-product semiflow admits a compact
global attractor $\mathcal A\subset\mathcal B\times\Theta$.
Strict inward motion on the boundary and invariance of the
attractor imply
\[
\mathcal A\subset\mathring{\mathcal B}\times\Theta.
\]

For the perturbed system, we impose the stronger condition below.

\medskip
\noindent{\bf(D-T-1)}
There exist $\delta,C>0$ such that
\[
f_i(t,x)\le-\delta
\quad\text{whenever }x_i\ge C,\quad
|x_j|\le x_i\ \text{for every existing neighbor }j\text{ of }i,
\]
and
\[
f_i(t,x)\ge\delta
\quad\text{whenever }x_i\le-C,\quad
|x_j|\le-x_i\ \text{for every existing neighbor }j\text{ of }i,
\]
uniformly in $t\in\R$.

Condition {\bf(D-T-1)} implies {\bf(D-T)}, possibly with a
different dissipativity constant, since
$x_if_i(t,x)\le-\delta|x_i|\le-\delta C$ under the stated
conditions. Both conditions extend to the joint hull by continuity.

Set
\[
M_g:=\sup_{(t,x)\in\R\times\R^n}\|g(t,x)\|.
\]
The same bound holds for the hull extensions.
We use the convention $\delta/M_g=+\infty$ when $M_g=0$ and assume
\begin{equation}\label{tri-small}
|\epsilon|<\frac{\delta}{M_g}.
\end{equation}

\begin{theorem}\label{Tr-c-thm}
Assume the regularity, boundedness, and tridiagonal hypotheses
above, together with {\bf(D-T-1)}.
There exists $\epsilon_1>0$ such that, whenever
$|\epsilon|<\epsilon_1$, every positive orbit of the skew-product
semiflow associated with \eqref{tri-equation-g1} is globally
defined and precompact in $\R^n\times\Theta$.
For each $(x^0,\theta)\in\R^n\times\Theta$, its $\omega$-limit set
$\omega^g(x^0,\theta)$ satisfies the three alternatives in
Theorem~\ref{perturbation-thm}. In particular, it contains at most
two minimal sets, each of which is an almost $1$-cover of the base
flow $\Theta$.
\end{theorem}

We first establish dissipativity of the perturbed system.

\begin{lemma}\label{dissipative-tri}
Assume {\bf(D-T-1)} and \eqref{tri-small}. Every solution
$\phi^g(t,x^0,\theta)$ of \eqref{tri-equation-g1} is defined for
all $t\ge0$, enters $\mathcal B$ in finite time, and remains there
afterwards.
\end{lemma}

\begin{proof}
Set
\[
R(t)=\|\phi^g(t,x^0,\theta)\|,
\qquad
\eta=\delta-|\epsilon|M_g>0.
\]
Whenever $R(t)\ge C$, define
\[
I(t):=\{i:|\phi_i^g(t,x^0,\theta)|=R(t)\}.
\]
Since $C>0$, every maximizing component is nonzero, and the upper
right Dini derivative satisfies
\[
D^+R(t)
=
\max_{i\in I(t)}
\left\{
\operatorname{sgn}\bigl(\phi_i^g(t,x^0,\theta)\bigr)
\frac{d}{dt}\phi_i^g(t,x^0,\theta)
\right\}.
\]
For $i\in I(t)$, all neighboring components have absolute value
at most $R(t)$. If $\phi_i^g(t,x^0,\theta)>0$, the first
inequality in {\bf(D-T-1)} gives
\[
\frac{d}{dt}\phi_i^g(t,x^0,\theta)
\le-\delta+|\epsilon|M_g=-\eta.
\]
If $\phi_i^g(t,x^0,\theta)<0$, the second inequality gives
\[
-\frac{d}{dt}\phi_i^g(t,x^0,\theta)
\le-\delta+|\epsilon|M_g=-\eta.
\]
Consequently,
\[
D^+R(t)\le-\eta
\qquad\text{whenever }R(t)\ge C.
\]

In particular, the vector field points strictly inward on every
active face of $\mathcal B$, so $\mathcal B$ is positively
invariant. Moreover,
\[
R(t)\le\max\{R(0),C\}
\]
throughout the maximal existence interval.
Since the base is compact and the vector field is continuous and
locally Lipschitz in the state variable, the continuation theorem
gives existence for all $t\ge0$.

If $R(0)>C$, integration of the Dini derivative inequality yields
\[
R(t)\le R(0)-\eta t
\]
as long as the trajectory remains outside $\mathcal B$.
Thus the solution enters $\mathcal B$ no later than
$(R(0)-C)/\eta$. If $R(0)\le C$, it starts in $\mathcal B$.
Positive invariance completes the proof.
\end{proof}

For every fixed $\bar\epsilon>0$ with $\bar\epsilon M_g<\delta$,
the entrance-time bound is uniform in $\theta$ and
$|\epsilon|\le\bar\epsilon$, for initial data in bounded sets.
Hence $\mathcal B\times\Theta$ is a common compact absorbing set
for this family of skew-product semiflows.

We next verify the NIC conditions for the unperturbed system.
Set
\[
\mathcal U_0=\R^n\times\Theta,
\qquad
K_0=\mathcal A.
\]
For $z=(x^1,x^2,\theta)\in\mathcal U_0^e$, let
\[
\psi(t):=\phi(t,x^1,\theta)-\phi(t,x^2,\theta).
\]
Then $\psi=(\psi_1,\dots,\psi_n)$ satisfies
\begin{equation}\label{linearity-tri-system}
\begin{split}
\dot\psi_1
&=a_{11}(t)\psi_1+a_{12}(t)\psi_2,\\
\dot\psi_i
&=a_{i,i-1}(t)\psi_{i-1}
+a_{ii}(t)\psi_i+a_{i,i+1}(t)\psi_{i+1},
\qquad 2\le i\le n-1,\\
\dot\psi_n
&=a_{n,n-1}(t)\psi_{n-1}+a_{nn}(t)\psi_n,
\end{split}
\end{equation}
where
\[
a_{ij}(t)
=
\int_0^1
\frac{\partial f_i}{\partial x_j}
\bigl(
\theta\cdot t,\,
s\phi(t,x^1,\theta)+(1-s)\phi(t,x^2,\theta)
\bigr)\,ds.
\]
The tridiagonal structure gives $a_{ij}(t)=0$ when $|i-j|>1$,
and \eqref{positive-jacobi-condition2} implies
\[
a_{i,i+1}(t)\ge\kappa_0,
\qquad
a_{i+1,i}(t)\ge\kappa_0,
\qquad t\ge0,\quad 1\le i\le n-1.
\]

Let $T(t,z)$ be the fundamental matrix of
\eqref{linearity-tri-system}, normalized by $T(0,z)=I_n$. Then
\[
T(t,z)(x^1-x^2)
=\phi(t,x^1,\theta)-\phi(t,x^2,\theta).
\]
Uniqueness and time translation of the coefficient equation give
the cocycle identity. Continuous dependence on the coefficients
gives {\bf(H1)}. Since the state space is finite-dimensional,
{\bf(H2)} holds automatically, and invertibility of the fundamental
matrix gives {\bf(H3)}.

Following \cite{Smillie1984}, define
\[
\mathcal R
=
\Bigl\{
x\in\R^n:
x_1\ne0,\ x_n\ne0,\ \text{and }
x_{i-1}x_{i+1}<0
\text{ whenever }x_i=0,\ 2\le i\le n-1
\Bigr\}.
\]
For $x\in\mathcal R$, let $\sigma(x)$ denote the number of sign
changes in the ordered sequence of its nonzero components.
Here $\sigma$ denotes the discrete sign-variation number.
The defining conditions on $\mathcal R$ make this number locally
constant on $\mathcal R$. Set
\[
C_0=\{0\},
\qquad
C_k:=\overline{\{x\in\mathcal R:\sigma(x)\le k-1\}},
\qquad 1\le k\le n,
\]
where the closure is taken in $\R^n$.
These sets form a nested family of solid $k$-cones, with
$C_n=\R^n$. The discrete sign-variation theorem for cooperative
tridiagonal systems gives the cone-preservation property in
{\bf(H4)}; see
\cite[Proposition~1.2]{simith1991} and
\cite[Lemma~2.1]{Fangchun2013}.

Define
\[
\pi:\R^n\to\R,\qquad \pi(x)=x_1,
\qquad
\mathcal H:=\pi^{-1}(0).
\]
Then $\mathcal H$ has codimension one. To verify {\bf(H5)},
suppose that
\[
T(t,z)v\in(C_i\setminus C_{i-1})\cap\mathcal H
\]
is nonzero and, contrary to the claim, that $v\in C_i$.
At regular times before $t$, the sign-variation number of
$T(s,z)v$ is at most $i-1$. Since its first component vanishes
at time $t$, the strict-drop property gives at most $i-2$ sign
changes at regular times immediately after $t$.
For $i\ge2$, taking these times down to $t$ and using closedness
of $C_{i-1}$ yields $T(t,z)v\in C_{i-1}$, a contradiction.
For $i=1$, a strict decrease from a sign-variation number already
equal to zero is impossible. Thus {\bf(H5)} holds with $d=1$.

Applying the abstract results gives the following tridiagonal
counterpart of Proposition~\ref{parobolic-thm}.

\begin{lemma}\label{Tri-thm}
There exist open bounded neighborhoods $\mathcal V_0$ and
$\mathcal W_0$ of $K_0$, with
\[
\overline{\mathcal V_0}\subset\mathcal W_0\subset\mathcal U_0,
\]
and a number $\epsilon_0>0$ with the following property.
For any skew-product semiflow
\[
\tilde\Lambda_t
=\bigl(\tilde\phi(t,\cdot,\theta),\theta\cdot t\bigr)
\]
on $\mathcal U_0$, with
$\tilde\phi(t,\cdot,\theta)\in
C^1(\mathcal W_0(\theta),\R^n)$ for $t\in[1/2,1]$, such that
\begin{equation}\label{Tri-c1}
\sup_{\substack{\frac12\le t\le1\\ \theta\in\Theta}}
\|\tilde\phi(t,\cdot,\theta)-\phi(t,\cdot,\theta)\|_
{C^1(\mathcal W_0(\theta),\R^n)}
<\epsilon_0,
\end{equation}
any minimal set of $\tilde\Lambda_t$ contained in $\mathcal V_0$
is an almost $1$-cover of the base flow $\Theta$.
Moreover, every precompact positive orbit that is eventually
contained in $\mathcal V_0$ and whose $\omega$-limit set is
contained in $\mathcal V_0$ satisfies the three alternatives in
Theorem~\ref{perturbation-thm}; in particular, its $\omega$-limit
set contains at most two minimal sets.
\end{lemma}

\begin{proof}[Proof of Theorem~\ref{Tr-c-thm}]
Fix $\bar\epsilon>0$ such that $\bar\epsilon M_g<\delta$.
By Lemma~\ref{dissipative-tri}, every positive orbit is globally
defined and bounded, hence precompact in $\R^n\times\Theta$.
Moreover, $\mathcal B\times\Theta$ is a common compact absorbing
set for $|\epsilon|\le\bar\epsilon$. Let
$\mathcal A^g_\epsilon$ denote the corresponding global attractor,
so that $\mathcal A^g_0=\mathcal A$.

The hull vector fields are $C^1$ in $x$, and their derivatives
depend continuously on the hull variable. Since $g$ and $D_xg$
are uniformly bounded, the perturbation $\epsilon g$ tends to
zero in the $C^1$ norm on bounded sets, uniformly in $\theta$.
The dissipative bound keeps solution segments with bounded initial
data in a common bounded set. Continuous dependence of solutions
and their variational equations therefore yields
\[
\sup_{\substack{\frac12\le t\le1\\ \theta\in\Theta}}
\|\phi^g_\epsilon(t,\cdot,\theta)-\phi(t,\cdot,\theta)\|_
{C^1(\mathcal W_0(\theta),\R^n)}
\longrightarrow0
\qquad(\epsilon\to0),
\]
where the subscript $\epsilon$ makes the parameter dependence
explicit. Thus \eqref{Tri-c1} holds for sufficiently small
$|\epsilon|$.

The common compact absorbing set and finite-time convergence of
the semiflows allow us to apply the upper-semicontinuity argument
used in the proof of Theorem~\ref{non-local-thm}. Consequently,
\[
\operatorname{dist}_{\R^n\times\Theta}
\bigl(\mathcal A^g_\epsilon,\mathcal A\bigr)
\longrightarrow0
\qquad(\epsilon\to0).
\]
Hence $\mathcal A^g_\epsilon\subset\mathcal V_0$ for all
sufficiently small $|\epsilon|$. Since this attractor is compact
and attracts every positive orbit, each such orbit eventually
remains in $\mathcal V_0$, and its $\omega$-limit set is contained
there.

Choose $\epsilon_1\le\bar\epsilon$ sufficiently small that both
\eqref{Tri-c1} and the preceding attractor inclusion hold for
$|\epsilon|<\epsilon_1$. Lemma~\ref{Tri-thm} completes the proof.
\end{proof}

The applications above illustrate how the abstract NIC framework
preserves the qualitative structure of $\omega$-limit sets under
small perturbations that need not retain the original local or
tridiagonal structure.

\section{Discussion and Further Remarks}

We conclude the discussion of the main results with remarks on their
scope and their relation to earlier work.

\medskip
\noindent\textbf{(1) The perturbation hypothesis and the NIC method.}
The fibrewise $C^1$ hypothesis is the uniform skew-product counterpart of
the assumption in Tere\v{s}\v{c}\'ak's semiflow theory \cite{Te}.
Lemma~\ref{small-pertur-lem} derives the one-step operator estimate from
this hypothesis, using the derivative bound for nearby points and the
$C^0$ bound for separated points. No additional perturbed difference
cocycle is needed. A $C^0$ estimate alone does not give the required
control of increments divided by the distance between the initial points.
In infinite dimensions, compactness and the regularizing image property
\eqref{regularizing-image} allow the cone chain to be truncated; for
scalar parabolic equations, this property follows from the Sturm
zero-number theorem.

Many systems with an NIC structure possess an integer-valued discrete
Lyapunov functional; see
\cite{BSZM,Fusco1987,Fusco1990,Gedeon,Smillie1984,
Ma-Sm,Ma-Se1,Ma-Se2,Ma-Se3,Matano,MS}.
Such a functional may depend on the precise form of the equation and
cease to be available after perturbation. The NIC method retains the
geometric and operator estimates underlying the restrictions on dynamics,
and thereby treats several classes of perturbations in a common framework.

\medskip
\noindent\textbf{(2) Relation to the periodic theory.}
Tere\v{s}\v{c}\'ak's work already includes both discrete-time and
continuous-time perturbation estimates
\cite[Propositions~5.1 and~6.1]{Te}. The extension developed here concerns
the recurrent dynamics over an almost periodic base. Without a common
return time, the cone estimates must be combined with information about
minimal extensions and the order in individual fibres. This combination
yields the almost $1$-covering property and the bound of two minimal sets
in Theorems~\ref{a-a-one} and~\ref{perturbation-thm}.

The finite-dimensional realization has a different scope: it requires
the all-time cone-layer hypothesis in Corollary~\ref{imbeding-d-plane}
and reduces only the fibre coordinate. It is not asserted for every
compact invariant set, and imposes no dimension bound on the base.

\medskip
\noindent\textbf{(3) Applications and possible extensions.}
Classical monotone-system results give generic convergence under
appropriate positivity assumptions; see \cite{Pol89} for autonomous
nonlocal equations and \cite{HP,Pol92} for time-periodic systems.
The nonlocal application here permits sign-changing coefficients, as
Example~\ref{example-quasiperiodic-nonlocal} illustrates. Its conclusions
hold for every positive orbit when the perturbation is sufficiently
small. The smallness condition matters: larger nonlocal perturbations
can generate complicated dynamics \cite{Fie-Pol}.

For periodic coefficients, the NIC theory gives convergence to periodic
dynamics under suitable dissipativity assumptions \cite{Te}. With
almost periodic forcing, our conclusions concern minimal extensions and
the structure of $\omega$-limit sets; they do not assert convergence to
a single almost periodic solution. The chemotaxis application likewise
allows general source terms satisfying {\bf(D)} and {\bf(P)}, rather
than requiring a logistic source or nonnegative initial data.
In all the applications of Section~\ref{applications-sec}, dissipativity,
compactness, and upper semicontinuity of attractors place the asymptotic
dynamics in the neighborhood where the abstract theory applies.

The parabolic applications use homogeneous Neumann boundary conditions.
The zero-number method also produces nested cones for homogeneous
Dirichlet or separated Robin conditions, but extending the applications
requires verifying transversality and the dissipative estimates in the
corresponding phase spaces.

\section{Proof of Proposition \ref{pertur-pro}}\label{proof-section}

Throughout this section we assume that {\bf(H1)--(H4)} are satisfied.

The proof adapts the cone estimates in
\cite[Propositions~5.1 and~6.1]{Te} to orbit pairs over the base flow.
We retain the quantitative bounds needed for the uniform perturbation
argument. The proof is organized as follows.

\begin{itemize}
\item[\rm(1)] We verify that the assumptions of Lemma~\ref{ES-lm} are
satisfied on the doubled invariant set $Z$, and hence obtain the
corresponding exponential separation bundles $V_z^i$ and
$\operatorname{Ann}(L_z^i)$.

\item[\rm(2)] We construct a finite family of perturbed cones
$\mathcal C_i$ near the original cone chain $\{C_i\}$, together with a
collection of quantitative constants that will be used later.

\item[\rm(3)] We prove Proposition~\ref{pertur-pro}(1) and (2) by combining a
one-step operator estimate derived from the fibrewise $C^1$ perturbation with
the invariance of the new cones and a contraction estimate outside the
terminal cone.

\item[\rm(4)] Assuming additionally {\bf(H5)}, we prove
Proposition~\ref{pertur-pro}(3) by introducing an intermediate family
$\mathcal W_i(s)$ describing the $i$-th cone layer and showing that
$\mathcal H$ cannot intersect this layer for $s\ll1$.
\end{itemize}

\subsection{Exponential separation along the doubled invariant set}

Define the doubled skew-product semiflow on $\mathcal U_0^e$ by
\[
\Lambda_t^e(x,y,\theta)
=
\bigl(
\phi(t,x,\theta),
\phi(t,y,\theta),
\theta\cdot t
\bigr),
\qquad t\ge0.
\]
Let
\[
Z=Z(K_0)
=
\left\{
(x,y,\theta):
(x,\theta),(y,\theta)\in K_0
\right\}.
\]
Since $K_0$ is a compact invariant set of $\Lambda_t$, the set $Z$ is a
compact invariant set of $\Lambda_t^e$. For
$z=(x,y,\theta)\in Z$, write
\[
z\cdot t
:=
\Lambda_t^e z
=
\bigl(
\phi(t,x,\theta),
\phi(t,y,\theta),
\theta\cdot t
\bigr).
\]

The restriction of $\Lambda_t$ to $K_0$ extends uniquely to a flow.
Indeed, invariance gives surjectivity of each positive-time map.
If two images coincide, invertibility of the base flow gives equality
of the initial base points, and \eqref{operator-T-eq1} together with
{\bf(H3)} gives equality of the state variables. Each time map is thus
a continuous bijection of the compact set $K_0$, hence a homeomorphism.
The induced maps on $Z$ likewise define a flow.

Let
\[
\{T(t,z):(t,z)\in[0,\infty)\times Z\}
\]
be the family of bounded linear operators defined in Section~2.2. By
{\bf(H1)},
\[
T(0,z)=I,
\]
and
\begin{equation}\label{difference-equation}
T(t,x,y,\theta)(x-y)
=
\phi(t,x,\theta)-\phi(t,y,\theta)
\end{equation}
for every $t\ge0$ and every $(x,y,\theta)\in Z$. Moreover,
\[
T(t+s,z)
=
T(s,z\cdot t)\circ T(t,z),
\qquad t,s\ge0.
\]
Hence
\[
(v,z)\longmapsto
\bigl(T(t,z)v,z\cdot t\bigr),
\qquad t\ge0,
\]
defines a linear vector bundle semiflow on $X\times Z$ over the compact
base flow $z\mapsto z\cdot t$.

By {\bf(H2)}, $T(t,z)$ is compact for every $t>0$ and $z\in Z$, while
{\bf(H4)} implies that it strongly preserves each cone $C_i$. Therefore,
for every admissible index $i$, all the assumptions of
Lemma~\ref{ES-lm} are satisfied with respect to the $k_i$-cone $C_i$.

For each admissible index $i$, let $V_z^i$ and
$\operatorname{Ann}(L_z^i)$ denote the invariant bundles of the
exponential separation along $Z$ associated with the cone $C_i$, as in
Definition~\ref{ES-define}. We next record the nesting properties of
these bundles.

For later use, fix $h\in[1/2,1]$. The bundle map
$(z,v)\mapsto(z\cdot h,T(h,z)v)$ satisfies assumption {\rm(H)} of
Tere\v{s}\v{c}\'ak: the base map is a homeomorphism, and
{\bf(H1)--(H4)} give continuity, compactness, nondegeneracy, and strong
cone preservation. The image-exhaustion property follows from
\eqref{regularizing-image} when $N=+\infty$, and from $C_N=X$ when
$N<+\infty$. By uniqueness of exponential separation, its leading and
dual bundles coincide with $V_z^i$ and $L_z^i$.

\begin{lemma}\label{invariant-space}
Assume that {\bf(H1)--(H4)} hold. Then, for every $z\in Z$ and every
$1\le i<N$,
\[
V_z^i\subset V_z^{i+1},
\qquad
L_z^i\subset L_z^{i+1}.
\]
Moreover,
\[
V_z^{i+1}
=
V_z^i\oplus
\bigl(
V_z^{i+1}\cap\operatorname{Ann}(L_z^i)
\bigr),
\qquad z\in Z,\quad 1\le i<N.
\]
For each $1\le i<N$, the family
\[
\bigcup_{z\in Z}
\left[
\bigl(
V_z^{i+1}\cap\operatorname{Ann}(L_z^i)
\bigr)\times\{z\}
\right]
\]
is a continuous finite-dimensional subbundle of $X\times Z$ of fibre
dimension $k_{i+1}-k_i$. Its unit-sphere bundle is compact, and it is
invariant under the linear vector bundle semiflow
\[
(v,z)\longmapsto\bigl(T(t,z)v,z\cdot t\bigr).
\]
\end{lemma}

\begin{proof}
Apply \cite[Lemma~5.2]{Te} to the time-one bundle map discussed above.
It gives the nesting, the direct-sum decomposition, and continuity of
the intersection bundle, whose fibre dimension is $k_{i+1}-k_i$.
Its unit-sphere bundle is compact because $Z$ is compact.
For every $t>0$, invariance of the exponential separations gives
\[
T(t,z)\bigl(V_z^{i+1}\cap\operatorname{Ann}(L_z^i)\bigr)
\subset V_{z\cdot t}^{i+1}\cap\operatorname{Ann}(L_{z\cdot t}^i).
\]
By {\bf(H3)}, the restriction is injective. Both spaces have dimension
$k_{i+1}-k_i$, so the inclusion is an equality.
\end{proof}

\subsection{Construction of the perturbed cone family \texorpdfstring{$\mathcal C_i$}{Ci}}

We first establish a uniform contraction estimate that will be used to
truncate the infinite cone chain. This estimate is a continuous-time
consequence of \cite[Lemma~5.1]{Te}.

\begin{lemma}\label{contract-lem}
Assume that {\bf(H1)--(H4)} hold and that $N=+\infty$. Then, for every
$0<\lambda_1<1$, there exists an integer $i_0\ge1$ such that
\begin{equation}\label{contract-equa1}
\|T(t,z)u\|
\le
\lambda_1^t\|u\|
\end{equation}
for every $(t,z)\in[\tfrac12,+\infty)\times Z$ and every
$u\in\operatorname{Ann}(L_z^{i_0})$.
\end{lemma}

\begin{proof}
Apply \cite[Lemma~5.1]{Te} to the time-$1/2$ bundle map, which satisfies
its hypotheses by the preceding discussion. Set
\[
\tau:=\max\left\{1,\sup_{\substack{0\le s\le1/2\\z\in Z}}
\|T(s,z)\|_{L(X,X)}\right\},\qquad \rho:=\lambda_1/\tau.
\]
Assumption {\bf(H1)} and $T(0,z)=I$ give $\tau<\infty$, so
$0<\rho<1$. The cited lemma provides $i_0\ge1$ such that
\begin{equation}\label{half-time-contraction}
\|T(\tfrac12,z)u\|\le\rho\|u\|,
\qquad z\in Z,\quad u\in\operatorname{Ann}(L_z^{i_0}).
\end{equation}
For $1/2\le t\le1$, the cocycle identity gives
\[
\|T(t,z)u\|
\le\|T(t-\tfrac12,z\cdot\tfrac12)\|\,
\|T(\tfrac12,z)u\|
\le\tau\rho\|u\|=\lambda_1\|u\|\le\lambda_1^t\|u\|.
\]
For $t>1$, choose $m=\lceil t\rceil$ and divide $t$ into $m$ equal
subintervals of length $t/m\in[1/2,1]$. Positive invariance of
$\operatorname{Ann}(L_z^{i_0})$ allows the preceding estimate to be
iterated, giving $\|T(t,z)u\|\le\lambda_1^t\|u\|$.
\end{proof}

By Lemma~\ref{ES-lm}, the bundles $\bigcup_{z\in Z}V_z^i\times\{z\}$ and $\bigcup_{z\in Z}L_z^i\times\{z\}$
are continuous, and there exist constants $M_i>0$ and $0<\gamma_i<1$
for every admissible index $i$, such that the corresponding exponential
separation estimates hold. Let $P_z^i$ denote the projection of $X$ onto $V_z^i$ along
$\operatorname{Ann}(L_z^i)$, and let
\[
Q_z^i:=I-P_z^i,\qquad z\in Z,\quad i\text{ admissible}.
\]
We also set
\[
V_z^0=\{0\},\qquad L_z^0=\{0\},\qquad P_z^0=0,\qquad Q_z^0=I
\]
for all $z\in Z$. We also use $C_z^0(s)=\{0\}$ for $s\ge0$.

If $N=\infty$, then for any fixed $0<\lambda_1<1$,
Lemma~\ref{contract-lem} yields $i_0\in\NN$ such that
\begin{equation}\label{contrac-condition}
P_z^{i_0}u=0
\quad\Longrightarrow\quad
\|T(t,z)u\|\le \lambda_1^t\|u\|,
\qquad t\in[\tfrac12,\infty).
\end{equation}
We then define
\[
N_0=
\begin{cases}
N, & N<\infty,\\[1mm]
i_0, & N=\infty.
\end{cases}
\]

For each $s\ge0$, $z\in Z$, and $1\le i\le N_0$, define
\[
C_z^i(s)=\{u\in X:\ \|Q_z^i u\|\le s\|P_z^i u\|\},
\qquad
\overline C_z^i(s)=C_z^i(s)\cap S,
\]
and
\[
D_z^i(s)=\{u\in X:\ \|P_z^i u\|\le s\|Q_z^i u\|\},
\qquad
\overline D_z^i(s)=D_z^i(s)\cap S,
\]
where $S=\{u\in X:\|u\|=1\}$. For $1\le i\le N_0$ and $s\ge0$, set
\[
\mathcal C_i(s)=\bigcup_{z\in Z}C_z^i(s),
\qquad
\mathcal D_i(s)=\bigcup_{z\in Z}D_z^i(s).
\]

The next lemma shows that, for $s>0$ sufficiently small, each
$\mathcal C_i(s)$ is itself a $k_i$-cone.

\begin{lemma}\label{small-cone-lem}
Assume that {\bf(H1)--(H4)} hold. Then there exists
$0<\delta_0\le1$ such that
\begin{equation}\label{inter-zero}
\mathcal C_i(s)\cap\mathcal D_i(s)=\{0\},
\qquad
0\le s\le\delta_0,\quad 1\le i\le N_0.
\end{equation}
Moreover, $\mathcal C_i(s)$ is a $k_i$-cone for every
$1\le i\le N_0$ and $0<s\le\delta_0$.
\end{lemma}

\begin{proof}
The continuity of the exponential separations implies that
$z\mapsto P_z^i$ and $z\mapsto Q_z^i$ are continuous in the operator
norm. Since $Z$ is compact,
\[
\kappa_i:=
\sup_{z\in Z}
\bigl(\|P_z^i\|_{L(X,X)}+\|Q_z^i\|_{L(X,X)}\bigr)
<+\infty,
\qquad 1\le i\le N_0.
\]

Suppose, to the contrary, that no $\delta_0>0$ satisfies
\eqref{inter-zero} simultaneously for all $1\le i\le N_0$. Then, for
every $n\ge1$, there exist $i_n\in\{1,\ldots,N_0\}$,
$s_n\in[0,1/n]$, and a unit vector
\[
u_n\in\mathcal C_{i_n}(s_n)\cap\mathcal D_{i_n}(s_n).
\]
Since $\mathcal C_i(s)$ and $\mathcal D_i(s)$ are increasing in $s$,
\[
u_n\in
\mathcal C_{i_n}(1/n)\cap\mathcal D_{i_n}(1/n).
\]
Passing to a subsequence, we may assume that $i_n=i_0$ for some fixed
$i_0\in\{1,\ldots,N_0\}$. Thus there exist $z_n^1,z_n^2\in Z$ such that
\begin{equation}\label{two-estimates}
\begin{aligned}
\|Q_{z_n^1}^{i_0}u_n\|
&\le \frac1n\|P_{z_n^1}^{i_0}u_n\|,\\
\|P_{z_n^2}^{i_0}u_n\|
&\le \frac1n\|Q_{z_n^2}^{i_0}u_n\|.
\end{aligned}
\end{equation}
Set
\[
v_n^j=P_{z_n^j}^{i_0}u_n,\qquad
w_n^j=Q_{z_n^j}^{i_0}u_n,
\qquad j=1,2.
\]
Then
\[
v_n^j\in V_{z_n^j}^{i_0},
\qquad
w_n^j\in\operatorname{Ann}(L_{z_n^j}^{i_0}).
\]
By the definition of $\kappa_{i_0}$ and \eqref{two-estimates},
\[
\|w_n^1\|\le\frac{\kappa_{i_0}}n,
\qquad
\|v_n^2\|\le\frac{\kappa_{i_0}}n,
\]
and hence
\[
w_n^1\to0,\qquad v_n^2\to0,\qquad \|v_n^1\|\to1.
\]

Since $Z$ is compact, after passing to a further subsequence we may assume
that $z_n^j\to z^j\in Z$ for $j=1,2$. The unit-sphere bundle of the
finite-dimensional continuous bundle $\{V_z^{i_0}\}_{z\in Z}$ is compact.
We may therefore assume further that
\[
v_n^1\to v\in V_{z^1}^{i_0},
\qquad \|v\|=1.
\]
Since $u_n=v_n^1+w_n^1$ and $w_n^1\to0$, we have $u_n\to v$. Moreover,
\[
w_n^2=u_n-v_n^2\to v.
\]
The continuity of the complementary bundle gives
$v\in\operatorname{Ann}(L_{z^2}^{i_0})$. Consequently,
\[
0\neq v
\in V_{z^1}^{i_0}\cap\operatorname{Ann}(L_{z^2}^{i_0})
\subset
C_{i_0}\cap\operatorname{Ann}(L_{z^2}^{i_0})
=\{0\},
\]
a contradiction. Hence \eqref{inter-zero} holds for some
$0<\delta_0\le1$.

It remains to prove that $\mathcal C_i(s)$ is a $k_i$-cone. Fix
$1\le i\le N_0$, $0<s\le\delta_0$, and $z_0\in Z$. By definition,
\[
\begin{aligned}
V_{z_0}^i
&\subset C_{z_0}^i(s)\subset\mathcal C_i(s),\\
\operatorname{Ann}(L_{z_0}^i)
&\subset D_{z_0}^i(s)\subset\mathcal D_i(s).
\end{aligned}
\]
It follows from \eqref{inter-zero} that
\[
\mathcal C_i(s)\cap\operatorname{Ann}(L_{z_0}^i)=\{0\}.
\]
Thus $\mathcal C_i(s)$ contains the $k_i$-dimensional subspace
$V_{z_0}^i$ and is transverse to the codimension-$k_i$ subspace
$\operatorname{Ann}(L_{z_0}^i)$. It is also invariant under
multiplication by nonzero real scalars.

We finally verify that $\mathcal C_i(s)$ is closed. Let
$u_n\in\mathcal C_i(s)$ with $u_n\to u$. Choose $z_n\in Z$ such that
\[
\|Q_{z_n}^iu_n\|\le s\|P_{z_n}^iu_n\|.
\]
After passing to a subsequence, assume that $z_n\to z\in Z$. By the
operator-norm continuity of the projections,
\[
P_{z_n}^iu_n\to P_z^iu,
\qquad
Q_{z_n}^iu_n\to Q_z^iu.
\]
Passing to the limit gives
\[
\|Q_z^iu\|\le s\|P_z^iu\|,
\]
so $u\in C_z^i(s)\subset\mathcal C_i(s)$. Hence
$\mathcal C_i(s)$ is closed and is therefore a $k_i$-cone.
\end{proof}

\subsubsection*{Choice of constants}

We now fix several constants and notations that will be used throughout the
rest of the section.

\begin{itemize}
\item Let
\[
c_P=\sup\{\|P_z^i\|:\ 1\le i\le N_0,\ z\in Z\},
\qquad
c_Q=\sup\{\|Q_z^i\|:\ 1\le i\le N_0,\ z\in Z\},
\]
and define
\[
\zeta=\max\left\{\sup_{(t,z)\in[\tfrac12,1]\times Z}\|T(t,z)\|,\,1\right\},
\qquad
r=\max\{c_Pc_Q,1\}.
\]

\item Let $\delta_0$ be as in \eqref{inter-zero}, and choose
$0<\delta<\delta_0$. When $N=\infty$, decrease $\delta$ so that
\begin{equation}\label{coeff-relation}
\lambda_0:=2\delta+\frac{\delta\zeta+\lambda_1^{1/2}}{1-\delta}<1.
\end{equation}
If {\bf(H5)} is assumed, choose $\delta$ from the outset so that
$2\delta/(1-\delta)<\delta_1$, where $\delta_1$ is supplied by
Lemma~\ref{plane-projection-lem}. That lemma depends only on the
unperturbed bundles and $\mathcal H$, so this choice precedes all
constants and cones defined in terms of $\delta$.
When $N<\infty$, set $\lambda_0=1/2$. In that case
$\mathcal C_{N_0}=X$, so the terminal-cone contraction assertion is vacuous.

\item Define
\[
M=\max\{M_i:\ 1\le i\le N_0\},
\qquad
\gamma=\min\{\gamma_i:\ 1\le i\le N_0\}\in(0,1),
\qquad
c=\frac1\delta.
\]
Then there exists $T_1\ge1$ such that
\begin{equation}\label{T1-choice}
cMe^{-\gamma t}\le (8r)^{-N_0}\delta
\qquad\text{for all }t\ge T_1.
\end{equation}

\item Put
\[
B=(\zeta+1)^{2T_1+2}.
\]
Choose $T_0>T_1+1$. When $N=\infty$, enlarge $T_0$ so that
\begin{equation}\label{large-const}
(\lambda_0-\delta)^{t-T_1}B\le \lambda_0^t
\qquad\text{for all }t\ge T_0.
\end{equation}
\end{itemize}

We now define the finite nested cone family
\[
\mathcal C_i:=\mathcal C_i\big((2r)^{i-N_0}\delta\big),
\qquad 1\le i\le N_0.
\]
If $1\le i<N_0$ and $u\in \mathcal C_i((2r)^{i-N_0}\delta)$, then there is
$z\in Z$ such that
\[
\|Q_z^i u\|\le (2r)^{i-N_0}\delta \|P_z^i u\|.
\]
By Lemma~\ref{invariant-space},
\[
Q_z^{i+1}u=Q_z^{i+1}Q_z^i u,
\qquad
P_z^iu=P_z^i P_z^{i+1}u.
\]
Hence
\[
\begin{split}
\|Q_z^{i+1}u\|
&=\|Q_z^{i+1}Q_z^i u\|
\le \|Q_z^{i+1}\|\,\|Q_z^i u\|\\
&\le (2r)^{i-N_0}\delta \|Q_z^{i+1}\|\,\|P_z^i u\|\\
&=(2r)^{i-N_0}\delta \|Q_z^{i+1}\|\,\|P_z^iP_z^{i+1}u\|\\
&\le (2r)^{i-N_0}\delta \|Q_z^{i+1}\|\,\|P_z^i\|\,\|P_z^{i+1}u\|\\
&\le r(2r)^{i-N_0}\delta \|P_z^{i+1}u\|\\
&=\frac12(2r)^{i+1-N_0}\delta \|P_z^{i+1}u\|.
\end{split}
\]
Thus every nonzero $u\in\mathcal C_i$ satisfies the defining inequality
for $C_z^{i+1}((2r)^{i+1-N_0}\delta)$ with a factor-two margin. By
continuity of $P_z^{i+1}$ and $Q_z^{i+1}$, a neighborhood of $u$ is
contained in this fixed fibre cone. Therefore
\[
u\in\operatorname{int}\mathcal C_{i+1}
\qquad (u\in\mathcal C_i\setminus\{0\}),
\]
and in particular
\[
\mathcal C_i\setminus\{0\}\subset\operatorname{int}\mathcal C_{i+1},
\qquad 1\le i<N_0.
\]

\subsection{Cone invariance under small perturbations}

We prove the eventual cone invariance required in
Proposition~\ref{pertur-pro}(1). For orbit segments that remain in the
perturbative neighborhood, suitable finite products of the auxiliary
one-step operators map $\mathcal C_i\setminus\{0\}$ into
$\operatorname{int}\mathcal C_i\setminus\{0\}$. Iteration then gives the
assertion for all times $t\ge T_0$.

It is immediate from \eqref{inter-zero} that
\begin{equation}\label{cone-parame-subset}
\mathcal C_i(s)\subset \mathcal C_i(\delta_0)
\subset \Bigl\{u\in X:\ \|Q_z^i u\|\le \frac1{\delta_0}\|P_z^i u\|\Bigr\}
\end{equation}
for all $0<s\le\delta_0$, $1\le i\le N_0$, and $z\in Z$.

By exponential separation,
\begin{equation}\label{exponential-inequa}
\|Q_z^i v\|\le c\|P_z^i v\|
\quad\Longrightarrow\quad
\|Q_{z\cdot t}^i T(t,z)v\|
\le (8r)^{-N_0}\delta\|P_{z\cdot t}^i T(t,z)v\|
\end{equation}
for all $t\ge T_1$, $z\in Z$, $v\in X$, and $1\le i\le N_0$. Combining
\eqref{exponential-inequa}, \eqref{cone-parame-subset}, and
\eqref{difference-equation}, we see that Proposition~\ref{pertur-pro}(1) is
true for the unperturbed semiflow $\phi$.

We first need a persistence property of the slices $\overline C_z^i(s)$ under
small perturbations.

\begin{lemma}\label{cone-nested}
Assume that $0\le s_1<s_2$ and $1\le i\le N_0$. Let
\[
\varrho_0=
\min\left\{
\frac{s_2-s_1}{2(1+c_Q+c_Ps_1)(1+s_1)},
\frac{1}{2c_P(1+s_1)}
\right\}.
\]
Then for any $z\in Z$ and $0<\varrho\le\varrho_0$,
\[
B(u,\varrho)\cap S\subset \overline C_z^i(s_2)
\qquad\text{for all }u\in \overline C_z^i(s_1),
\]
where $B(u,\varrho)=\{u_0\in X:\ \|u-u_0\|<\varrho\}$.
\end{lemma}

\begin{proof}
Fix $u_1\in \overline C_z^i(s_1)$ and choose $u_2\in S$ with
$\|u_2-u_1\|<\varrho_0$. Then
\begin{equation}\label{cone-inequa}
\begin{split}
\|Q_z^i u_2\|
&=\|Q_z^i(u_2-u_1)+Q_z^i u_1\|\\
&\le c_Q\varrho_0+s_1\|P_z^i u_1\|\\
&\le c_Q\varrho_0+s_1\bigl(\|P_z^i(u_1-u_2)\|+\|P_z^i u_2\|\bigr)\\
&\le (c_Q+s_1 c_P)\varrho_0+s_1\|P_z^i u_2\|\\
&\le \frac{s_2-s_1}{2(1+s_1)}+s_1\|P_z^i u_2\|.
\end{split}
\end{equation}
On the other hand,
\[
\|P_z^i u_2\|
\ge \|P_z^i u_1\|-\|P_z^i(u_2-u_1)\|
\ge \frac{1}{1+s_1}-c_P\varrho_0
\ge \frac{1}{2(1+s_1)}.
\]
Combining this with \eqref{cone-inequa} yields
\[
\|Q_z^i u_2\|\le s_2\|P_z^i u_2\|,
\]
hence $u_2\in \overline C_z^i(s_2)$.
\end{proof}

We next introduce the perturbed skew-product semiflow
\begin{equation}\label{pertur-skew-product-semiflow}
\tilde{\Lambda}_t(x,\theta)
=
\bigl(\tilde{\phi}(t,x,\theta),\theta\cdot t\bigr),
\qquad t\geq 0,
\end{equation}
and its natural extension to pairs,
\begin{equation}\label{pertur-skew-product-semiflow-e}
\tilde{\Lambda}_t^e(x,y,\theta)
=
\bigl(
\tilde{\phi}(t,x,\theta),
\tilde{\phi}(t,y,\theta),
\theta\cdot t
\bigr).
\end{equation}

No linear difference cocycle is assumed for the perturbed semiflow.
Instead, we construct auxiliary one-step linear operators that represent
differences of perturbed trajectories. For
$\mu\in[\tfrac12,1]$, define
\[
H_{\mu,\theta}(x)
:=
\tilde{\phi}(\mu,x,\theta)-\phi(\mu,x,\theta).
\]
Thus $H_{\mu,\theta}$ measures the difference between the perturbed and
unperturbed time-$\mu$ fibre maps.

Let $\tilde{z}=(x,y,\theta)\in\mathcal U_0^e$ with $x\neq y$.
By the Hahn--Banach theorem, there exists
$\ell_{x,y}\in X^*$ such that
\[
\|\ell_{x,y}\|=1,
\qquad
\ell_{x,y}(x-y)=\|x-y\|.
\]
Define
\begin{equation}\label{auxiliary-difference-operator}
\widehat{T}(\mu,\tilde{z})v
:=
T(\mu,\tilde{z})v
+
\frac{\ell_{x,y}(v)}{\|x-y\|}
\bigl(
H_{\mu,\theta}(x)-H_{\mu,\theta}(y)
\bigr),
\qquad v\in X.
\end{equation}
For $x=y$, set
\[
\widehat{T}(\mu,x,x,\theta)
:=
T(\mu,x,x,\theta).
\]
It follows from \eqref{auxiliary-difference-operator} that
\begin{align}
\widehat{T}(\mu,x,y,\theta)(x-y)
&=
T(\mu,x,y,\theta)(x-y)
+H_{\mu,\theta}(x)-H_{\mu,\theta}(y)
\notag\\
&=
\tilde{\phi}(\mu,x,\theta)
-\tilde{\phi}(\mu,y,\theta).
\label{one-step-difference-identity}
\end{align}
For $x\ne y$, we also have
\begin{equation}\label{rank-one-correction-estimate}
\bigl\|
\widehat{T}(\mu,\tilde{z})-T(\mu,\tilde{z})
\bigr\|_{L(X,X)}
\leq
\frac{
\|H_{\mu,\theta}(x)-H_{\mu,\theta}(y)\|
}{
\|x-y\|
}.
\end{equation}
For $x=y$, the two operators agree by definition, so their difference
has norm zero.
The uniform estimate needed below will be derived from the fibrewise
$C^1$-perturbation assumption after the neighborhoods
$\mathcal V_0$ and $\mathcal W_0$ have been chosen appropriately;
see Lemma~\ref{small-pertur-lem}.

If the segment $[x,y]$ is contained in the relevant fibre domain,
one may instead define
\[
\widehat{T}(\mu,\tilde{z})
=
T(\mu,\tilde{z})
+
\int_0^1
D_xH_{\mu,\theta}\bigl(y+s(x-y)\bigr)\,ds.
\]
Both constructions satisfy \eqref{one-step-difference-identity},
although the resulting operators need not coincide. The latter is the
averaged-derivative construction used in \cite[(5.31)]{Te}.

Let
\[
\boldsymbol\mu=(\mu_1,\ldots,\mu_q),
\qquad
\mu_j\in[\tfrac12,1],
\]
and set
\[
s_0=0,
\qquad
s_j=\mu_1+\cdots+\mu_j,
\qquad 1\leq j\leq q.
\]
For an initial perturbed pair
$\tilde{z}_0=(x_0,y_0,\theta_0)$, put
\[
\tilde{z}_j
=
\tilde{\Lambda}_{s_j}^e\tilde{z}_0,
\qquad 0\leq j\leq q,
\]
and suppose that
$\tilde{z}_j\in\mathcal U_0^e$ for $0\leq j\leq q-1$.
Define
\begin{equation}\label{perturbed-product}
\widehat{\mathcal T}_{\boldsymbol\mu}(\tilde{z}_0)
:=
\widehat{T}(\mu_q,\tilde{z}_{q-1})
\circ\cdots\circ
\widehat{T}(\mu_1,\tilde{z}_0).
\end{equation}
Repeated application of \eqref{one-step-difference-identity} gives
\begin{equation}\label{perturbed-product-difference}
\widehat{\mathcal T}_{\boldsymbol\mu}(\tilde{z}_0)(x_0-y_0)
=
\tilde{\phi}(s_q,x_0,\theta_0)
-
\tilde{\phi}(s_q,y_0,\theta_0).
\end{equation}

No continuity or cocycle property of $\widehat{T}$ is required.
The argument below uses only the one-step identity
\eqref{one-step-difference-identity} and the one-step estimate
\eqref{pertur-opera}.

For comparison, let $z_0\in Z$ and put
\[
z_j=z_0\cdot s_j,
\qquad 0\leq j\leq q.
\]
The cocycle identity for the unperturbed difference operators gives
\[
T(\mu_q,z_{q-1})
\circ\cdots\circ
T(\mu_1,z_0)
=
T(s_q,z_0).
\]
Thus it remains to compare
$\widehat{\mathcal T}_{\boldsymbol\mu}(\tilde{z}_0)$
with $T(s_q,z_0)$.

The following lemma is the key perturbative cone-invariance statement.

\begin{lemma}\label{equiv-lem-i}
There exists $\epsilon_1>0$ with the following property.
Let $z_0\in Z$, let $\tilde{z}_0$ be an initial perturbed pair,
and let
\[
\boldsymbol\mu=(\mu_1,\ldots,\mu_q),
\qquad
\mu_j\in[\tfrac12,1],
\qquad
t:=s_q\in[T_1,2T_0].
\]
Put
\[
z_j=z_0\cdot s_j,
\qquad
\tilde{z}_j
=
\tilde{\Lambda}_{s_j}^e\tilde{z}_0,
\qquad 0\leq j\leq q,
\]
and assume that
$\tilde{z}_j\in\mathcal U_0^e$ for $0\leq j\leq q-1$.
Suppose that
\begin{equation}\label{pertur-opera}
\bigl\|
\widehat{T}(\mu_j,\tilde{z}_{j-1})
-
T(\mu_j,z_{j-1})
\bigr\|_{L(X,X)}
<
\epsilon_1,
\qquad 1\leq j\leq q.
\end{equation}
Then, for every $1\leq i\leq N_0$,
\begin{equation}\label{In-cone}
\widehat{\mathcal T}_{\boldsymbol\mu}(\tilde{z}_0)
\bigl(C_{z_0}^i(c)\setminus\{0\}\bigr)
\subset
C_{z_0\cdot t}^i\bigl((4r)^{-N_0}\delta\bigr)
\setminus\{0\}.
\end{equation}
Consequently,
\begin{equation}\label{In-cone1}
\widehat{\mathcal T}_{\boldsymbol\mu}(\tilde{z}_0)
\bigl(\mathcal C_i\setminus\{0\}\bigr)
\subset
\operatorname{int}(\mathcal C_i)\setminus\{0\},
\qquad 1\leq i\leq N_0.
\end{equation}
\end{lemma}

\begin{proof}
Put $a=(8r)^{-N_0}\delta$ and $b=(4r)^{-N_0}\delta$.
By \eqref{exponential-inequa} and {\bf(H3)},
\begin{equation}\label{unperturbed-cone-compression}
T(t,z_0)\bigl(C_{z_0}^i(c)\setminus\{0\}\bigr)
\subset C_{z_0\cdot t}^i(a)\setminus\{0\}.
\end{equation}
Lemma~\ref{cone-nested} gives $\varrho>0$, independent of $z\in Z$
and $1\le i\le N_0$, such that
\begin{equation}\label{uniform-angular-margin}
B(u,\varrho)\cap S\subset\overline C_z^i(b),
\qquad u\in\overline C_z^i(a).
\end{equation}

Compactness of the unit-sphere bundles of $V_z^i$, continuity of $T$,
and {\bf(H3)} give a uniform $m>0$ with
$\|T(\tau,z)v\|\ge m\|v\|$ for $v\in V_z^i$,
$T_1\le\tau\le2T_0$, and $1\le i\le N_0$.
Invariance of the splitting gives
$P_{z\cdot\tau}^iT(\tau,z)=T(\tau,z)P_z^i$.
Since $\|P_z^iu\|\ge(1+c)^{-1}$ for
$u\in\overline C_z^i(c)$, we obtain
\begin{equation}\label{uniform-lower-image}
\|T(\tau,z)u\|
\ge c_P^{-1}\|T(\tau,z)P_z^iu\|
\ge\frac{m}{c_P(1+c)}=:m_0>0.
\end{equation}

Write $A_j=T(\mu_j,z_{j-1})$ and
$B_j=\widehat T(\mu_j,\tilde z_{j-1})$.
For $\epsilon_1\le1$, \eqref{pertur-opera} and the bound $\zeta$
fixed above give $\|A_j\|\le\zeta$ and $\|B_j\|\le\zeta+1$.
Also $q\le q_*:=\lceil4T_0\rceil$, since $q/2\le s_q\le2T_0$.
The telescoping identity
\[
B_q\cdots B_1-A_q\cdots A_1
=\sum_{k=1}^q B_q\cdots B_{k+1}(B_k-A_k)A_{k-1}\cdots A_1
\]
and the unperturbed cocycle property imply
\begin{equation}\label{finite-product-close}
\|\widehat{\mathcal T}_{\boldsymbol\mu}(\tilde z_0)-T(t,z_0)\|_{L(X,X)}
\le C_{T_0}\epsilon_1,
\qquad C_{T_0}:=q_*(\zeta+1)^{q_*-1}.
\end{equation}
Choose $0<\epsilon_1<\min\{1,\delta\}$ so that
\begin{equation}\label{epsilon1-choice}
C_{T_0}\epsilon_1
<\min\{m_0/2,\varrho m_0/4\}.
\end{equation}

By homogeneity, take $u\in C_{z_0}^i(c)$ with $\|u\|=1$, and write
$Au=T(t,z_0)u$ and
$Bu=\widehat{\mathcal T}_{\boldsymbol\mu}(\tilde z_0)u$.
The preceding estimates give $\|Au\|\ge m_0$ and
$\|Bu\|>m_0/2$, and hence
\[
\left\|\frac{Bu}{\|Bu\|}-\frac{Au}{\|Au\|}\right\|
\le\frac{4\|Bu-Au\|}{m_0}
\le\frac{4C_{T_0}\epsilon_1}{m_0}<\varrho.
\]
By \eqref{unperturbed-cone-compression}, the normalized unperturbed
image lies in $\overline C_{z_0\cdot t}^i(a)$.
Equation~\eqref{uniform-angular-margin} therefore places the normalized
perturbed image in $\overline C_{z_0\cdot t}^i(b)$, proving
\eqref{In-cone}. Finally, \eqref{cone-parame-subset} and
Lemma~\ref{cone-nested} give
\begin{equation}\label{equivalent-iv}
\mathcal C_i\subset C_z^i(c),
\qquad C_z^i\bigl((4r)^{-N_0}\delta\bigr)\setminus\{0\}
\subset\operatorname{int}(\mathcal C_i),
\end{equation}
for every $z\in Z$ and $1\le i\le N_0$, yielding \eqref{In-cone1}.
\end{proof}

The next lemma derives the factor estimate \eqref{pertur-opera} from the
single $C^1$ perturbation hypothesis.

\begin{lemma}\label{small-pertur-lem}
There exist open bounded neighborhoods $\mathcal V_0$ and $\mathcal W_0$
of $K_0$, with
\[
\overline{\mathcal V_0}
\subset
\mathcal W_0
\subset
\mathcal U_0,
\]
and constants $\epsilon_0>0$ and $\eta>0$ with the following property.
Assume \eqref{perturb-c1} and set
\[
\mathcal V_0^e
:=
\left\{
(x,y,\theta):
x,y\in\mathcal V_0(\theta)
\right\}.
\]
Let $\tilde z_0\in\mathcal V_0^e$, and suppose that its perturbed orbit
segment of length $s_q\leq2T_0$ remains in $\mathcal V_0^e$. Let
\[
\boldsymbol\mu=(\mu_1,\ldots,\mu_q),
\qquad
\mu_j\in[\tfrac12,1],
\]
be a partition of this segment, and put
\[
s_0=0,
\qquad
s_j=\mu_1+\cdots+\mu_j,
\qquad
\tilde z_j
=
\tilde\Lambda_{s_j}^e\tilde z_0,
\quad 1\leq j\leq q.
\]
Then there exists $z_0\in Z$ such that, with
\[
z_j=z_0\cdot s_j,
\qquad 0\leq j\leq q,
\]
one has
\begin{equation}\label{finite-shadowing}
d_e(\tilde z_j,z_j)<\eta,
\qquad 0\leq j\leq q,
\end{equation}
and
\begin{equation}\label{factor-estimate-small-perturbation}
\bigl\|
\widehat T(\mu_j,\tilde z_{j-1})
-
T(\mu_j,z_{j-1})
\bigr\|_{L(X,X)}
<
\epsilon_1,
\qquad 1\leq j\leq q.
\end{equation}
In particular, \eqref{pertur-opera} holds for every $1\leq j\leq q$.
\end{lemma}

\begin{proof}
Fix a compatible metric $d_\Theta$ on $\Theta$ and use
\[
d_e\bigl((x,y,\theta),(x',y',\theta')\bigr)
=\|x-x'\|+\|y-y'\|+d_\Theta(\theta,\theta').
\]
Since $Z$ is compact and $\mathcal U_0^e$ is open, continuity of $T$
allows us to choose $\eta>0$ such that the $\eta$-neighborhood of $Z$
lies in $\mathcal U_0^e$ and
\begin{equation}\label{operator-uniform-continuity}
\begin{gathered}
z\in Z,\quad z'\in\mathcal U_0^e,\quad d_e(z',z)<\eta
\quad\Longrightarrow\quad
\|T(\mu,z')-T(\mu,z)\|_{L(X,X)}<\epsilon_1/2
\end{gathered}
\end{equation}
uniformly for $\mu\in[1/2,1]$.
Recall $q_*:=\lceil4T_0\rceil$, so that $q\le q_*$.
Continuity of $(\mu,z)\mapsto\Lambda_\mu^e z$ near the compact set
$[1/2,1]\times Z$ gives radii
$0<r_0<\cdots<r_{q_*}<\eta$, chosen recursively backwards, such that
\begin{equation}\label{finite-step-modulus}
d_e(w,z)<r_j,\quad z\in Z
\quad\Longrightarrow\quad
d_e(\Lambda_\mu^e w,\Lambda_\mu^e z)<r_{j+1}/2
\end{equation}
for $0\le j<q_*$, $w\in\mathcal U_0^e$, and $\mu\in[1/2,1]$.

Choose an open bounded neighborhood $\mathcal W_0$ of $K_0$ whose
closure lies in $\mathcal U_0$. Compactness provides a positive uniform
metric margin around $K_0$ inside $\mathcal W_0$. A sufficiently small
metric neighborhood $\mathcal V_0$ of $K_0$ then satisfies
$\overline{\mathcal V_0}\subset\mathcal W_0$, and there is $a_0>0$
such that
$[x,y]\subset\mathcal W_0(\theta)$ whenever
$x,y\in\mathcal V_0(\theta)$ and $\|x-y\|<a_0$.
We may also require $\mathcal V_0^e$ to lie in the $r_0$-neighborhood
of $Z$. Indeed, if this failed for arbitrarily small neighborhoods,
there would be pairs $(x_n,y_n,\theta_n)$ staying a positive distance
from $Z$ while both component points approached $K_0$. Compactness
would give limits $(x,\theta),(y,\theta)\in K_0$ with the same base
coordinate, contradicting $(x,y,\theta)\in Z$.

The $C^0$ part of \eqref{perturb-c1} gives
\begin{equation}\label{one-step-map-close}
d_e(\tilde\Lambda_\mu^e w,\Lambda_\mu^e w)
\le C_e\epsilon_0,\qquad C_e=2,
\end{equation}
for $w\in\mathcal V_0^e$ and $\mu\in[1/2,1]$.
Set $C_0:=\max\{1,2/a_0\}$ and choose $\epsilon_0>0$ such that
\begin{equation}\label{epsilon0-choice-small-perturbation}
C_e\epsilon_0<\frac12\min_{1\le j\le q_*}r_j,
\qquad C_0\epsilon_0<\epsilon_1/2.
\end{equation}
For $\tilde z_0\in\mathcal V_0^e$, choose $z_0\in Z$ with
$d_e(\tilde z_0,z_0)<r_0$. Inductively,
\eqref{finite-step-modulus} and \eqref{one-step-map-close} give
\[
d_e(\tilde z_{j+1},z_{j+1})
< C_e\epsilon_0+r_{j+1}/2<r_{j+1},\qquad 0\le j<q,
\]
where the orbit-segment hypothesis ensures
$\tilde z_j\in\mathcal V_0^e$. Since $r_j<\eta$, this proves
\eqref{finite-shadowing}.

For the factor estimate, fix $\mu$ and $\theta$ and write
$H=H_{\mu,\theta}$. If $x,y\in\mathcal V_0(\theta)$ are less than
$a_0$ apart, the mean-value formula on
$[x,y]\subset\mathcal W_0(\theta)$ bounds
$\|H(x)-H(y)\|$ by $\epsilon_0\|x-y\|$.
Otherwise, the $C^0$ bound gives
$\|H(x)-H(y)\|\le2\epsilon_0
\le(2\epsilon_0/a_0)\|x-y\|$. Thus
\begin{equation}\label{H-uniform-Lipschitz-estimate}
\|H(x)-H(y)\|\le C_0\epsilon_0\|x-y\|,
\qquad x,y\in\mathcal V_0(\theta).
\end{equation}
The rank-one estimate \eqref{rank-one-correction-estimate} implies
\begin{equation}\label{derived-operator-closeness}
\|\widehat T(\mu,x,y,\theta)-T(\mu,x,y,\theta)\|_{L(X,X)}
\le C_0\epsilon_0.
\end{equation}
This also holds at $x=y$, where the two operators agree.
Combining it with \eqref{operator-uniform-continuity} and
\eqref{finite-shadowing}, we obtain
\[
\begin{aligned}
&\|\widehat T(\mu_j,\tilde z_{j-1})-T(\mu_j,z_{j-1})\|_{L(X,X)}\\
&\quad\le
\|\widehat T(\mu_j,\tilde z_{j-1})-T(\mu_j,\tilde z_{j-1})\|_{L(X,X)}
+\|T(\mu_j,\tilde z_{j-1})-T(\mu_j,z_{j-1})\|_{L(X,X)}\\
&\quad<C_0\epsilon_0+\epsilon_1/2<\epsilon_1,
\qquad 1\le j\le q.
\end{aligned}
\]
This proves \eqref{factor-estimate-small-perturbation}.
\end{proof}
\subsection{Contraction outside the terminal cone}

We now prepare the proof of Proposition~\ref{pertur-pro}(2).

\begin{lemma}\label{outside-contrac}
Assume $N=\infty$, and let $\delta$ and $\lambda_0$ be as chosen in
\eqref{coeff-relation}. Let $u\in X$ and $z\in Z$ satisfy
\[
\|P_z^{N_0}u\|\le \delta \|Q_z^{N_0}u\|.
\]
Then
\[
\|T(t,z)u\|\le (\lambda_0-2\delta)\|u\|
\qquad\text{for all }t\in[\tfrac12,1].
\]
\end{lemma}

\begin{proof}
Write
\[
u=P_z^{N_0}u+Q_z^{N_0}u=:v+w,
\qquad \|v\|\le \delta\|w\|.
\]
By \eqref{contrac-condition} and \eqref{coeff-relation},
\[
\begin{split}
\|T(t,z)(v+w)\|
&\le \|T(t,z)\|\,\|v\|+\lambda_1^t\|w\|\\
&\le \delta\|T(t,z)\|\,\|w\|+\lambda_1^{1/2}\|w\|\\
&\le \frac{\|T(t,z)\|\delta+\lambda_1^{1/2}}{1-\delta}(1-\delta)\|w\|\\
&\le (\lambda_0-2\delta)(1-\delta)\|w\|\\
&\le (\lambda_0-2\delta)\|v+w\|.
\end{split}
\]
\end{proof}

\begin{lemma}\label{equiv-lem-iii}
Assume $N=\infty$. Let $z_0\in Z$, let $\tilde z_0$ be an initial perturbed
pair, let $u_0\in X$, and let
$\boldsymbol\mu=(\mu_1,\ldots,\mu_q)$,
$\mu_j\in[\tfrac12,1]$, be a partition of
$t_0\in[T_0,2T_0]$ such that $s_p=t_0-T_1$ for some $p$, and suppose
that \eqref{pertur-opera} holds. For
$1\le j\le q$, set
\[
u_j:=\widehat T(\mu_j,\tilde z_{j-1})\cdots
\widehat T(\mu_1,\tilde z_0)u_0.
\]
Assume that $u_0,u_q\notin\mathcal C_{N_0}$.
Then
\begin{equation}\label{equivalent-iii}
\|u_q\|\le \lambda_0^{t_0}\|u_0\|.
\end{equation}
\end{lemma}

\begin{proof}
We first show that
\begin{equation}\label{outside-contrac-eq}
u_j\notin C_{z_j}^{N_0}(c),
\qquad 0\le j\le p.
\end{equation}
Indeed, if $u_j\in C_{z_j}^{N_0}(c)$ for some $j\le p$, then the remaining
time $t_0-s_j$ belongs to $[T_1,2T_0]$. Applying
Lemma~\ref{equiv-lem-i} to the remaining factors would give
$u_q\in\mathcal C_{N_0}$, a contradiction.

Since $c=1/\delta$, \eqref{outside-contrac-eq} implies
\[
\|P_{z_j}^{N_0}u_j\|<\delta\|Q_{z_j}^{N_0}u_j\|,
\qquad0\le j\le p.
\]
Using Lemma~\ref{outside-contrac}, \eqref{pertur-opera}, and
$\epsilon_1<\delta$, we therefore obtain, for $1\le j\le p$,
\[
\begin{split}
\|u_j\|
&\le \|T(\mu_j,z_{j-1})u_{j-1}\|
   +\epsilon_1\|u_{j-1}\|\\
&\le(\lambda_0-\delta)\|u_{j-1}\|.
\end{split}
\]
Because the number of factors in the first part is at least $t_0-T_1$,
\[
\|u_p\|\le(\lambda_0-\delta)^{t_0-T_1}\|u_0\|.
\]
The last interval has length $T_1$ and contains at most $2T_1+2$ factors.
Their norms are bounded by $\zeta+1$, so its product has norm at most $B$.
Consequently, \eqref{large-const} gives
\[
\|u_q\|\le B(\lambda_0-\delta)^{t_0-T_1}\|u_0\|
\le\lambda_0^{t_0}\|u_0\|.
\]
This proves \eqref{equivalent-iii}.
\end{proof}

\subsection{Proof of Proposition \ref{pertur-pro}(1) and (2)}

We now complete the proof of the first two parts of
Proposition~\ref{pertur-pro}. The key point is that the assumed closeness of
the time maps in the fibrewise $C^1$ topology implies the one-step operator
estimate \eqref{pertur-opera}. Lemmas~\ref{equiv-lem-i} and
\ref{equiv-lem-iii} can therefore be applied to products along orbit
segments inside $\mathcal V_0$.

\begin{proof}[Proof of Proposition~\ref{pertur-pro}(1) and (2)]
Let $x,y\in\mathcal V_0(\theta)$ be distinct and assume
that
\[
\tilde\phi(s,x,\theta),\ \tilde\phi(s,y,\theta)\in \mathcal V_0(\theta\cdot s)
\quad\text{for all }0\le s\le t.
\]

\medskip
\noindent\textbf{Proof of part {\rm(1)}.}
Assume $x-y\in \mathcal C_i\setminus\{0\}$ for some $1\le i\le N_0$ and
$t\ge T_0$. Put $q=\lfloor t/T_0\rfloor$ and $\ell=t/q$. Then
$q\ge1$, $T_0\le\ell\le2T_0$, and $t=q\ell$. On each of the $q$
orbit segments of length $\ell$, choose a partition with step lengths in
$[1/2,1]$. Lemma~\ref{small-pertur-lem} gives \eqref{pertur-opera}, while
\eqref{perturbed-product-difference} identifies the resulting product with
the difference at the end of the segment. Lemma~\ref{equiv-lem-i}, applied
successively, therefore gives
\[
\tilde\phi(t,x,\theta)-\tilde\phi(t,y,\theta)
\in \operatorname{int}(\mathcal C_i)\setminus\{0\}.
\]
This proves Proposition~\ref{pertur-pro}(1).

\medskip
\noindent\textbf{Proof of part {\rm(2)}.}
Assume now that
\[
\tilde\phi(t,x,\theta)-\tilde\phi(t,y,\theta)\notin \mathcal C_{N_0}
\qquad\text{for some }t\ge T_0.
\]
If $N<\infty$, then $N_0=N$ and $\mathcal C_{N_0}=X$, so the
premise is impossible. Assume henceforth that $N=\infty$. Put
$q=\lfloor t/T_0\rfloor$, $\ell=t/q$, and let
\[
u_j:=\tilde\phi(j\ell,x,\theta)-\tilde\phi(j\ell,y,\theta),
\qquad0\le j\le q.
\]
Thus
$T_0\le\ell\le2T_0$. Part $(1)$, applied to every remaining orbit
segment, shows that $u_j\notin\mathcal C_{N_0}$ for all
$0\le j\le q$; otherwise the terminal vector would belong to
$\operatorname{int}(\mathcal C_{N_0})$. On each block, choose a partition
of $\ell-T_1$ followed by a partition of $T_1$, all with step lengths in
$[1/2,1]$. Lemmas~\ref{small-pertur-lem} and \ref{equiv-lem-iii}, together
with \eqref{perturbed-product-difference}, yield
\[
\|u_q\|\le \lambda_0^\ell\|u_{q-1}\|
\le\cdots\le \lambda_0^{q\ell}\|u_0\|
=\lambda_0^t\|u_0\|.
\]
Since $u_q=\tilde\phi(t,x,\theta)-\tilde\phi(t,y,\theta)$ and $u_0=x-y$,
this is the desired estimate.
This proves Proposition~\ref{pertur-pro}(2).
\end{proof}

\subsection{Transversality and proof of Proposition \ref{pertur-pro}(3)}

In this subsection we also assume that {\bf(H5)} holds. The goal is to show
that if an orbit difference stays in the same cone layer at times $0$ and
$2T_0$, then it cannot lie in the codimension-$d$ subspace $\mathcal H$
at the intermediate time $T_0$.

For any $z\in Z$, $1\le i\le N_0$, and $s\ge0$, define
\[
W_z^i(s)
:=\left\{u\in X:\ \|Q_z^i u+P_z^{i-1}u\|
\le s\|Q_z^{i-1}P_z^i u\|\right\}
\]
and
\[
\mathcal W_i(s):=\bigcup_{z\in Z}W_z^i(s),
\qquad s\ge0,\ 1\le i\le N_0.
\]
These sets describe a thin tubular neighborhood of the $i$-th intermediate
bundle
\[
W_z^i(0)=V_z^i\cap \operatorname{Ann}(L_z^{i-1})
\subset (C_i\setminus C_{i-1})\cup\{0\}.
\]
Moreover, if $0\le s_1<s_2$, then
$W_z^i(s_1)\subset W_z^i(s_2)$ and
$\mathcal W_i(s_1)\subset \mathcal W_i(s_2)$.

\begin{lemma}\label{plane-projection-lem}
Assume that {\bf(H1)--(H5)} are satisfied. Then there exists $\delta_1>0$ such that
\begin{equation}\label{plane-inter-zero}
\mathcal H\cap \mathcal W_i(\delta_1)=\{0\}
\end{equation}
and
\begin{equation}\label{projection-subset}
I_s^i:=
\bigcup_{z\in Z}\Bigl\{
u\in X:\ \|P_z^{i-1}u\|\le s\|Q_z^{i-1}u\|,
\ \|Q_z^i u\|\le s\|P_z^i u\|
\Bigr\}
\subset \mathcal W_i\Bigl(\frac{2s}{1-s}\Bigr)
\end{equation}
for all $0\le s<1$ and $1\le i\le N_0$.
\end{lemma}

\begin{proof}
We first prove \eqref{plane-inter-zero}. The unit-sphere bundle associated with
\[
\bigcup_{z\in Z}\bigl[
\bigl(V_z^i\cap \operatorname{Ann}(L_z^{i-1})\bigr)\times\{z\}\bigr]
\]
is compact, and the underlying bundle is invariant. Moreover,
\[
V_z^i\cap \operatorname{Ann}(L_z^{i-1})
\subset (C_i\setminus C_{i-1})\cup\{0\},
\]
assumption {\bf(H5)} implies that
\[
\mathcal H\cap \mathcal W_i(0)=\{0\}.
\]
Indeed, otherwise there would exist $z\in Z$ and a nonzero
$u\in \mathcal H\cap W_z^i(0)\subset (C_i\setminus C_{i-1})\cap \mathcal H$.
For any $t>0$, invariance and invertibility on the finite-dimensional bundle
give $z_{-t}\in Z$ and
$v\in V_{z_{-t}}^i\cap\operatorname{Ann}(L_{z_{-t}}^{i-1})\subset C_i$
such that $T(t,z_{-t})v=u$.  This contradicts {\bf(H5)}.  Compactness of the
unit-sphere bundle and closedness of $\mathcal H$ now give a uniform angular
separation; by continuity this separation persists on
$\mathcal W_i(\delta_1)$ for some $\delta_1>0$.

We now prove \eqref{projection-subset}. Let $u\in I_s^i$. Then $u$ can be
written as
\[
u=w_1+w_2+v,
\]
where
\[
w_1\in \operatorname{Ann}(L_z^i),\qquad
w_2\in V_z^i\cap \operatorname{Ann}(L_z^{i-1}),\qquad
v\in V_z^{i-1}.
\]
By the definition of $I_s^i$,
\[
\|v\|\le s\|w_1+w_2\|\le s\|w_1\|+s\|w_2\|,
\qquad
\|w_1\|\le s\|w_2+v\|\le s\|w_2\|+s\|v\|.
\]
Hence
\[
\|w_1+v\|\le \|w_1\|+\|v\|
\le \frac{2s}{1-s}\|w_2\|.
\]
Therefore
\[
u\in W_z^i\Bigl(\frac{2s}{1-s}\Bigr),
\]
which proves \eqref{projection-subset}.
\end{proof}

Recall that, when {\bf(H5)} holds, the choice of constants above ensures
\[
0<\delta<\min\left\{\delta_0,\frac{\delta_1}{2+\delta_1}\right\}.
\]
We can now prove the final part of Proposition~\ref{pertur-pro}.

\begin{proof}[Proof of Proposition~\ref{pertur-pro}(3)]
Let
\[
\Delta(t):=\tilde\phi(t,x,\theta)-\tilde\phi(t,y,\theta),
\qquad u=\Delta(0)=x-y,
\]
and assume that
\[
u,\ \Delta(2T_0)\in \mathcal C_i\setminus \mathcal C_{i-1}
\]
for some $1\le i\le N_0$. It suffices to show that
\[
\Delta(T_0)\notin \mathcal H.
\]

Choose a partition of $[0,2T_0]$ by steps in $[1/2,1]$ having $T_0$ as a
grid point. Lemma~\ref{small-pertur-lem} supplies $z_0\in Z$ and the factor
estimates \eqref{pertur-opera} along this partition.

By \eqref{cone-parame-subset}, we have
\begin{equation}\label{more-gener-supp}
u\in C_{z_0}^i(c),
\qquad
\Delta(2T_0)
\in C_{z_0\cdot 2T_0}^i(c)\setminus
C_{z_0\cdot 2T_0}^{\,i-1}\bigl((4r)^{-N_0}\delta\bigr).
\end{equation}
The vector $\Delta(T_0)$ is nonzero, since otherwise the semiflow property
would imply $\Delta(2T_0)=0$.
If $i=1$, it therefore lies outside $C_{z_0\cdot T_0}^0(c)=\{0\}$.
For $i\ge2$, if $\Delta(T_0)\in C_{z_0\cdot T_0}^{\,i-1}(c)$, then applying
Lemma~\ref{equiv-lem-i} to the factors over the second interval of length
$T_0$ would imply
\[
\Delta(2T_0)\in
C_{z_0\cdot 2T_0}^{\,i-1}\bigl((4r)^{-N_0}\delta\bigr),
\]
contradicting \eqref{more-gener-supp}.
Thus
\begin{equation}\label{middle-not-lower-cone}
\Delta(T_0)\notin C_{z_0\cdot T_0}^{\,i-1}(c).
\end{equation}

On the other hand, since $u\in C_{z_0}^i(c)$, applying
Lemma~\ref{equiv-lem-i} to the first interval gives
\[
\Delta(T_0)\in
C_{z_0\cdot T_0}^{\,i}\bigl((4r)^{-N_0}\delta\bigr).
\]
Together with \eqref{middle-not-lower-cone} and $c=1/\delta$, this implies
\[
\Delta(T_0)\in I_\delta^i.
\]
Hence, by \eqref{projection-subset},
\[
\Delta(T_0)\in \mathcal W_i\Bigl(\frac{2\delta}{1-\delta}\Bigr).
\]
Our choice of $\delta$ guarantees
\[
\frac{2\delta}{1-\delta}\le \delta_1.
\]
Therefore
\[
\Delta(T_0)\in \mathcal W_i(\delta_1).
\]
If $\Delta(T_0)\in \mathcal H$, then
\[
\Delta(T_0)\in \mathcal H\cap \mathcal W_i(\delta_1),
\]
which contradicts \eqref{plane-inter-zero}. Hence
\[
\Delta(T_0)\notin \mathcal H.
\]
This completes the proof of Proposition~\ref{pertur-pro}(3).
\end{proof}

\section*{Acknowledgement}
The author was partially supported by the National Natural Science Foundation of China 
under Grant Nos. 12671213 and 12331006. The author is grateful to Professor Wenxian Shen 
of Auburn University and Professor Yi Wang of the University of Science and Technology 
of China for their constructive comments and stimulating discussions.

\section*{Declaration of interests}

The author declares that there is no conflict of interest. The manuscript has no associated data.

\end{document}